\magnification=\magstep1
\input amstex
\input epsf
\documentstyle{amsppt}
%\nologo

\pagewidth{6.0truein}
\pageheight{8.5truein}
\define\bSym{\text{\bf Sym}}
\define\Sym{\operatorname{Sym}}
\define\Ker{\operatorname{Ker}}
\define\card{\operatorname{card}}
\define\dom{\operatorname{dom}}
\define\ov#1{\overline {#1}}
\NoBlackBoxes
\topmatter
\thanks
Vladimir Retakh was partially supported by NSF, NSA, and Max-Planck-Institut
f\"ur Mathematik (Bonn). 
\endthanks
\keywords
quasidetermininats, noncommutative algebra,
symmetric functions
\endkeywords
\subjclass Primary
05E05, 15A15, 16W30 
\endsubjclass
\title
Quasideterminants
\endtitle
\leftheadtext{I. Gelfand, S. Gelfand,V. Retakh, R. L. Wilson}
\author
Israel Gelfand, Sergei Gelfand,\\ Vladimir Retakh, Robert Lee Wilson
\endauthor

\address
I.~G.: Department of Mathematics, Rutgers University,
Piscataway, NJ 08854
\endaddress
\email
igelfand\@math.rutgers.edu
\endemail
\address
S.~G.: American Mathematical Society, 201, Charles Street,
Providence, RI 02904
\endaddress
\email
sxg\@ams.org
\endemail
\address
V.~R.: Department of Mathematics, Rutgers University,
Piscataway, NJ 08854
\endaddress
\email
vretakh\@math.rutgers.edu
\endemail
\address
R.~W.: Department of Mathematics, Rutgers University,
Piscataway, NJ 08854
\endaddress
\email
rwilson\@math.rutgers.edu
\endemail

\dedicatory
To the memory of Gian-Carlo Rota 
\enddedicatory
\endtopmatter

\abstract

The determinant is a main organizing tool in commutative
linear algebra. In this review we present a theory of
the quasideterminants defined for matrices over a division ring.

\endabstract 

\document

{\centerline {\bf Contents}}
\medskip

0. Introduction

1. General theory and main identities

2. Important example: quaternionic quasideterminants

3. Noncommutative determinants

4. Noncommutative Pl\"ucker and flag coordinates 

5. Factorization of Vandermonde Quasideterminants and the Vi\`ete Theorem

6. Noncommutative symmetric functions

7. Universal quadratic algebras associated with pseudo-roots of
noncommutative polynomails and noncommutative differential polynomials

8. Noncommutative traces, determinants and cyclic vectors

9. Some applications

\medskip

\head 0. Introduction\endhead

The ubiquitous notion of a determinant has a long history, both visible and
invisible. The determinant has been a main organizing tool in commutative
linear algebra and we cannot accept the point of view of a modern textbook
\cite {FIS} that ``determinants ... are of much less importance than they once
were".

Attempts to define a determinant for matrices with noncommutative entries 
started more than 150 years ago and include several great names.  For many
years the most famous examples of matrices of noncommutative objects were
quaternionic matrices and block matrices. It is not suprising that the first
noncommutative determinants or similar notions were defined for such
structures.

A. Cayley \cite {C} was the first to define, in 1845, the determinant of a
matrix with noncommutative entries. He mentioned that for a quaternionic matrix
$A=\left (\matrix a_{11}&a_{12}\\a_{21}&a_{22}\endmatrix \right )$  the
expressions $a_{11}a_{22}-a_{12}a_{21}$ and $a_{11}a_{22}- a_{21}a_{12}$ are
different and suggested choosing one of them as the determinant of the matrix
$A$. The analog of this construction for $3\times 3$-matrices was also proposed
in \cite {C} and later developed in \cite {J}. This ``naive" approach is now
known to work for quantum determinants and some other cases. Different forms of
quaternionic determinants were considered later by E. Study \cite {St}, E.H.
Moore \cite {Mo} and F. Dyson \cite {Dy}.

There were no direct ``determinantal" attacks on block matrices (excluding
evident cases) but important insights were given by G. Frobenius \cite {Fr} and
I. Schur \cite {Schur}.

A theory of determinants of matrices with general noncommutative entries was in
fact originated by J.H.M. Wedderburn in 1913. In  \cite {W} he constructed a
theory of noncommutative continued fractions or, in modern terms,
``determinants" of noncommutative  Jacobi matrices. 

In 1926-1928 A. Heyting \cite {H} and A. Richardson \cite {Ri, Ri1} suggested
analogs of a determinant for matrices over division rings. Heyting is known as
a founder of intuitionist logic and Richardson as a creator of the
Littlewood-Richardson rule. Heyting tried to construct a noncommutative
projective geometry. As a  computational tool, he introduced the ``designant"
of a noncommutative matrix. The designant of a $2\times 2$-matrix $A=(a_{ij})$
is defined as $a_{11}-a_{12}a_{22}^{-1}a_{21}$. The designant of an $n\times
n$-matrix is defined then by a complicated inductive procedure. The inductive
procedures used by Richardson were even more complicated. It is important to
mention that determinants of Heyting and Richardson in general are rational
functions  (and not polynomials!) in matrix entries. 

The idea to have non-polynomial determinants was strongly criticized by O. Ore
\cite {O}. In \cite {O} he defined a polynomial determinant for matrices over
an imporatnt class of noncommutative rings (now known as Ore rings).

The most famous and widely used noncommutaive determinant is the Dieudonne
determinant. It was defined for matrices over a division ring $R$ by J.
Dieudonne in 1943 \cite {D}. His idea was to consider determinants with values
in $R^*/[R^*, R^*]$ where $R^*$ is the monoid of invertible elements in $R$.
The properties of Dieudonne determinants are close to those of commutative
ones, but, evidently, Dieudonne determinants cannot be used for solving systems
of linear equations.

An interesting generalization of commutative determinants belongs to F. Berezin
\cite {B, Le}. He defined determinants for matrices over so called
super-commutative algebras. In particular, Berezin also understood that it is
impossible to avoid rational functions in matrix entries in his definition.   

Other famous examples of noncommutative determinants developed for different
special cases are:  quantum determinants \cite {KS, Ma}, Capelli determinants
\cite {We}, determinants introduced by Cartier-Foata \cite {CF, F} and
Birman-Williams \cite {BW}, etc. As we explain later (using another universal
notion, that of quasideterminants) these determinants and the determinants of
Dieudonne, Study, Moore, etc., are related to each other much more than one
would expect.

The notion of quasideterminants for matrices over a noncommutative division
ring was introduced in \cite {GR, GR1, GR2}. Quasideterminants are defined in
the ``most noncommutative case", namely, for matrices over free division rings.
We believe that quasideterminants should be one of main organizing tools in
noncommutative algebra giving them the same role determinants play in
commutative algebra.  The quasideterminant is not an analog of the commutative
determinant but rather of a ratio of the determinant of an $n\times n$-matricx
to the determinant of an $(n-1)\times (n-1)$-submatrix.

The main property of quasideterminants is a ``heredity principle'': let $A$ be
a square matrix over a division ring and $(A_{ij})$  a block decomposition of
$A$ (into submatrices of $A$). Consider the $A_{ij}$'s as elements of a matrix
$X$. Then the quasideterminant of the matrix X will be a matrix  $B$, and
(under natural assumptions) the quasideterminant of $B$ is equal to a suitable
quasideterminant of $A$. Since determinants of block matrices are not defined,
there is no analog of this principle for ordinary (commutative) determinants.

Quasideterminants have been effective in many areas including noncommutative
symmetric functions \cite {GKLLRT}, noncommutative integrable systems \cite
{RS, EGR, EGR1}, quantum algebras and Yangians \cite {GR, GR1, GR2, KL, Mol,
Mol1, MolR}, and so on \cite {P, Sch, RSh, RRV}. Quasideterminants and related
quasi-Pl\"ucker coordinates are also important in various approches to
noncommutative algebraic geometry (e.g., \cite {K, KR, SvB})

Many areas of noncommutative mathematics (Ore rings, rings of differential
operators, theory of factors, ``quantum mathematics'', Clifford algebras, etc)
were developed separately from each other. Our approach shows an advantage of
working with totally noncommutative variables (over free rings and division
rings). It leads us to a large variety of results, and their specialization to
different noncommutative areas implies known theorems with additional
information.

The price one pays for this is a huge number of inversions in rational
noncommutative expressions. The minimal number of successive inversions
required to express an element is called the height of this element. This
invariant (inversion height) reflects the  ``degree of noncommutativity'' and
it is of a great interest by itself.

Our experience shows that in dealing with noncommutative objects one should not
imitate the classical commutative mathematics, but follow ``the way it is''
starting with basics. In this paper we concentrate on two problems:
noncommutative Pl\"ucker coordinates (as a background of a noncommutative
geometry) and the noncommutative Bezout and Vi\`ete theorems (as a background
of noncommutative algebra). We apply the obtained results to the theory of
noncommutative symmetric functions started in \cite{GKLLRT}.

We have already said that the universal notion of a determinant has a long
history, both visible and  invisible.  The visible history of determinants
comes from the fact that they are constructed from another class of universal
objects: matrices. 

The invisible history of determinants is related with the Heredity principle
for matrices: matrices can be viewed as matrices with matrix entries (block
matrices) and some matrix properties come from the corresponding properties of
block matrices. In some cases, when the matrix entries of the block matrix
commute, the  determinant of a matrix can be computed in terms of  the
determinants of its blocks, but in general it is not possible:  the determinant
of a matrix with matrix entries is not defined because the entries do not
commute. In other words, the determinant does not satisfy the Heredity
principle.

Quasideterminants are defined for matrices over division rings and satisfy the
Heredity Principle. Their definition can be specialized for matrices over a
ring (including noncommutative rings) and can  be connected with  different
"famous" determinants. This reflects another general principle: in many cases
noncommutative algebra can be made  simpler and more natural than commutative
algebra.

The survey describes the first 10 years of development of this very active
area, and we hope that future work will bring many new interesting resutls. 

The paper is organized as follows. In Section 1 a definition of
quasideterminants is given and the main properties of quasideterminants
(including the Heredity principle) are described.

In Section 2 we discuss an important example: quasideterminants of quaternionic
matrices. These quasideterminants can be written as polynomials with real
coefficients in the matrix entries and their quaternionic conjugates.

As we already mentioned, mathematics knows a lot of different versions of
noncommutative determinants. In Section 3 we give a general definition of
determinants of noncommutative matrices (in general, there are many
determinants of a fixed matrix) and show how to obtain some well-known
noncommutative determinants as specializations of our definition.

In Section 4 we introduce noncommutative versions of Pl\"ucker and flag
coordinates for rectangular matrices over division rings.

In Section 5 we discuss two related classical problems for noncommutative
polynomials in one variable: how to factorize a polynomial into products of
linear polynomials and how to express the coefficients of a polynomial in terms
of its roots.

This results obtained in Section 5 lead us to a theory of noncommutative
symmetric functions (Section 6) and to universal quadratic algebras associated
with so-called pseudo-roots of noncommutative polynomails and noncommutative
differential polynomials (Section 7).

In Section 8 we present another approach to the theory of noncommutative
determinants, traces, etc., and relate it to the results presented in Sections
3 and 5.

Some applications to noncommutative continued fractions, characteristic
functions of graphs, noncommutative orthogonal polynomials and integrable
systems are given in Section 9.

We thank the referee for careful reading of the manuscript.

\head 1. General theory and main identities 
\endhead

\subhead  1.1.  The division ring of rational functions in free
variables\endsubhead
Throughout the paper we will work with rings of fractions of various
noncommutative rings. There are several ways to define rings of fractions in
the noncommutative case. We will use the approach developed by Amitsur,
Bergman and P. M. Cohn (for a detailed exposition see, e.g., \cite{Co}). The
advantage of this approach is that it is constructive; its disadvantage is that
it does not look very natural.

First, we define the free division ring generated by a finite set.
Let $X=\{x_1,\dots,x_m\}$ be a finite set. 
Define $\Cal F(X)$ as the free algebra generated by $m+2$ 
elements $0,1,x_1,\dots , x_m$, unary operations $a\mapsto -a$,
$a\mapsto a^{-1}$, and binary operations $+$ and $\times $, so that $\Cal F(X)$
contains such elements as $(x - x)^{-1}$ and even $0^{-1}$. 
No commutativity, associativity, distributivity, or other conditions are
imposed, so, that, in particular, elements $(x_1+x_2)\times x_3$ and $x_1\times
x_3 +x_2\times x_3$ are distinct, and three elements 
$$
(-x_1)^{-1}, \quad (0-x_1)^{-1}, \quad -x_1^{-1}\times 1^{-1},
$$
are also distinct. Elements of $\Cal F(X)$ are called {\it formulas\/} or
rational expressions over $X$. 

Denote by $\Cal P(X)$ the subset of $\Cal F(X)$ consisting of formulas without
division, i.e., without operation $(\ )^{-1}$.
                  
Now let $R$ be a $\Bbb Q$-algebra. By a partial homomorphism of
$\Cal F(X)$ to $R$ we mean a pair $(G,\beta)$ consisting of
a subset $G\subset \Cal F(X)$ and a map $\beta:G\to R$ such that
 
\roster
\item "(i)" $0,1 \in G$ and $\beta (0)=0$, $\beta (1)=1$,

\item "(ii)" if $a_1=-b$, $a_2=b+c$, $a_3=b\times c$ are elements in $G$
then $b, c \in G$ and $\beta (a_1)=-\beta (b)$,
$\beta (a_2)=\beta (b)+\beta (c)$, $\beta (a_3)=\beta (b)\beta (c)$.

\item"(iii)"  Let $b\in G$ and let $\beta (b)$ be invertible in $R$. Then
$b^{-1}\in G$ and $\beta (b^{-1}) = (\beta(b))^{-1}$.
\endroster

Let again $R$ be a $\Bbb Q$-algebra and $\alpha:X\to R$ an arbitrary map.
We say that a partial homomorphism ($G,\beta)$ of $\Cal F(X)$ to $R$ is an
extension of a map $\alpha$ if, in addition to (i)--(iii), the following
condition is satisfied. 

\roster
\item"(iv)" For $i=1,\dots,m$ we have $x_i\in G$ and $\beta(x_i)=\alpha(x_i)$. 
\endroster

Clearly, for an arbitrary $\alpha$, conditions (i), (ii), (iv) determine a
natural extension $(\Cal P(X), \alpha_{\Cal P})$, and for any other extension
$(G,\beta)$ we have $G\supset \Cal P(X)$, $\beta\mid_{\Cal P(X)} = \alpha_{\Cal
P}$. 

For two extension $(G_1,\beta_1)$ and $(G_2,\beta_2)$ of $\alpha$
define their intersection $(G,\beta)$ as follows:
$$
\gather
\text{$G$ is the set of all $a\in G_1\cap G_2$ such that
$\beta_1(a)=\beta_2(a)$,}\\
\text{$\beta(a)=\beta_1(a)=\beta_2(a)$ for $a\in G$.}
\endgather
$$
Clearly, $(G,\beta)$ is again an extension of $\alpha$. Therefore, the
intersection of all  extensions of $\alpha$ is again an extension. We call it
the {\it canonical extension\/} of  $\alpha$ and denote by
$(G_0,\overline\alpha$), or simply $\overline\alpha$.

Since each $\alpha:X\to R$ admits at least one extension (for example, $(\Cal
P( X), \alpha_\Cal P)$), the definition of the canonical extension makes sense.

If $(G_0,\overline\alpha)$ is the canonical extension of $\alpha$
and $a\in G_0$ we say that  $\alpha$ {\it can be evaluated at $a$.}

\smallskip

Let $D$ be a division ring over $\Bbb Q$. Denote by $M(X,D)$ the set of 
all maps
$\alpha : X\to D$. Clearly, $M(X,D)$ is isomorphic to $D^m$, where
$m=\card X$.

\definition{Definition 1.1.1} (i) The {\it domain\/} $\dom a$ of an element
$a\in \Cal F(X)$ is subset of $M(X,D) $
consisting of the maps $\alpha:X\to D$
such that $\alpha$ can be evaluated at $a$.

(ii) An element $a\in \Cal F(X)$ is called {\it
nondegenerate\/} if $\dom a\ne \emptyset$, and {\it degenerate\/} otherwise.

(iii) Two elements $a_1,a_2\in \Cal F(X)$ are called {\it equivalent\/} if
they are both nondegenerate and $\overline\alpha(a_1)=\overline\alpha(a_2)$ for
all $\alpha\in \dom a_1\cap \dom a_2$.
\enddefinition

For example, for $x\in X$, the elements $x-x\in \Cal F(X)$ is nondegenerate and
equivalent  to $0\in \Cal F(X)$ whereas $(x-x)^{-1}$ is degenerate. Another
example:  for $x\in X$,  the element $a_1=(1-x)^{-1}+(1-x^{-1})^{-1}$ is
equivalent to $a_2=1$.

\proclaim{Theorem 1.1.2 \rm (\cite {Co}, Section 7.2)}
{\rm(i)} If $a_1,a_2\in \Cal F(X)$ are both nondegenerate, then
$\dom_1\cap\dom a_2\ne \emptyset$.

 {\rm(ii)} Assume, in addition, that $D$ is a division ring with the
center $\Bbb Q$.
Then the equivalence classes of elements in $\Cal F(X)$
form a division ring, called $F_D(X)$.

 {\rm(iii)} If division rings $D_1$ and $D_2$ with center $\Bbb Q$ are
infinite dimensional over  $\Bbb Q$, then the projections $\Cal F(X)\to
F_{D_1}$ and $\Cal F(X)\to F_{D_2}$ induce an isomorphism $F_{D_1}\sim
F_{D_2}$.
\endproclaim

Part (iii) of Theorem 1.1.2 allows us to identify the division rings
$F_D(X)$ for all division rings $D$ infinite-dimensional over $\Bbb Q$.
We denote this division ring by $F(X)$ and called it
the {\it free division ring\/} generated by $X$.
For example, if $X$ consists of one element $x$, then
$F(X)=\Bbb Q(x)$ is the field of rational functions over $\Bbb Q$ in one
variable.

Elements $f\in F(X)$ are called (noncommutative) {\it rational functions\/}  in
variables $x\in X$, and  any element $a\in \Cal F(X)$ in the equivalence class
$f$ is called a  {\it rational expression\/} of the function $f$.

\remark{Remark}Similar results hold if $\Bbb Q$ is replaced by an arbitrary
field  $k$ of characteristic 0 (for example, by $\Bbb C$).
\endremark

The next proposition shows that for an arbitrary $\Bbb Q$-algebra $R$
evaluations of a map $\alpha:X\to R$ on two equivalent elements coincide.

\proclaim{Proposition 1.1.3} Let $R$ be a $\Bbb Q$-algebra,  $\alpha:X\to R$ a
map, and $(G_0, \overline\alpha)$ the canonical extension of $\alpha$. If
$a_1,a_2\in \Cal F(X)$ are equivalent and both lie in $G_0\subset \Cal F(X)$,
then $ \overline\alpha(a_1)= \overline\alpha(a_2)$.
\endproclaim

\definition{Definition 1.1.4} Let $\alpha:X\to R$ and let $(G_0,
\overline\alpha)$ be the canonical extension of $\alpha$. We say that $\alpha$ can be
evaluated at $f\in F(X)$  if there exists $a\in \Cal F(X)$ in the equivalence
class of $f$ such that $a\in G_0$, and in this case we define the value of
$\alpha$ at $f$ by the formula $\alpha(f) = \overline\alpha (a)$.
\enddefinition

Proposition 1.1.3 shows that this definition makes sense.
\smallskip
Cohn has shown (see \cite {Co}, Section 7.2) that the division ring $F(X)$
can be characterized by a universality property as follows.

\proclaim{Theorem 1.1.5} There exists a division ring $F(X)$ over $\Bbb Q$
and a monomorphism of algebras $\theta:\Bbb Q\langle X\rangle
\to F(X)$ with the following property.

If $D$ is an arbitrary division ring and
$$
\varphi : \Bbb Q\langle X\rangle \rightarrow D,
$$
a homomorphism, then there is a unique pair $(R,\psi)$ consisting of a
subring $R\subset F(X)$ containing
$\theta(\Bbb Q\langle X\rangle)$
and a homomorphism
$$
\psi : R\rightarrow D
$$
such that $\varphi = \psi\theta$ and
if $a\in R$ and $\psi (a)\neq 0$, then $a^{-1}\in R$.

The pair $(F(X),\theta)$ is determined uniquely up to a unique isomorphism.
\endproclaim

To conclude this subsection, we recall the definition of inversion height
(see, for example, \cite {Re}).

\definition{Definition 1.1.6} (i)  The inversion height of a
formula  $a\in\Cal F(X)$ is the maximal number of nested inversions in $a$.

(ii) The inversion height of an element $f\in F(X)$ is the smallest inversion
height of a formula in the equivalence class $f$.
\enddefinition

\example{Examples} (i) The inversion height of a polynomial in
generators $x\in X$ equals zero.

(ii) The inversion height of the ratio
of two polynomials $PQ^{-1}$ equals $0$ if $P$ is right divisible by
$Q$ (i.e., there exists a polynomial $R$ such that $P=RQ$),
and $1$ otherwise.

In the next two examples, let $x,y,z$ be three different elements of $X$.

(iii)  Consider the elements  $a_1,a_2\in \Cal F(X)$ given by the formulas
$a_1=(1-x)^{-1}+(1-x^{-1})^{-1}$ and
$a_2=x^{-1}+x^{-1}(z^{-1}y^{-1}-x^{-1})^{-1}x^{-1}$. Let $f_1$ and $f_2$ be
the corresponding elements in $F(X)$. Then the inversion height of $a_1$ and
$a_2$ equals $2$. On the other hand, in $F(X)$ we have $a_1=1$ and
$a_2=(x-yz)^{-1}$. Hence, the height of $f_1$ equals $0$ and the height of
$f_2$ equals $1$.

(iv) The height of the element $f\in F(X)$ given by the formula
$(x-yw^{-1}z)^{-1}$ equals $2$.
\endexample

\subhead 1.2. Definition of quasideterminants\endsubhead

Let $I,J$ be two finite sets of the same cardinality $n$ and $\Cal X$ be the
set of  $n^2$ elements $x_{ij}$, $1\leq i,j \leq n$. Denote by $F(\Cal X)$ the
free division ring generated by $\Cal X$ (see 1.1). Let $X$ be the 
$n\times
n$-matrix over $F(\Cal X)$  with rows indexed by elements of $I$, columns
indexed by elements of $J$, and the $(i,j)$-th entry equal to $x_{ij}\in
F(\Cal X)$.

\proclaim{Proposition 1.2.1} The matrix $X$ is inverible over $F(\Cal X)$.
\endproclaim

\demo{Proof } The proof is by inducion in $n$. Let us assume, for simplicity,
that $I=J=\{1,\dots,n\}$.

For $n=1$, $X=(x_{11})$ and the inverse matrix $Y=X^{-1}$ is $Y=(y_{11})$,
where $y_{11}=(x_{11})^{-1}$.

Let $n\geq 2$. Represent $X=(x_{ij})$ as a $2\times 2$ block matrix according
to the decompositions $\{1,\dots,n\}=\{1,\dots,n-1\}\cup\{n\}$ of $I$ and $J$,
$$
X=\pmatrix X_{11} & X_{12} \\ X_{21} & X_{22} \endpmatrix
$$
so that $X_{11}$, $X_{12}$, $ X_{21}$, $X_{22}$ are matrices of order
$(n-1)\times(n-1)$,  $(n-1)\times 1$, $1\times(n-1)$, and $1\times 1$
respectively. Then one can directly verify that the matrix $Y$ given in the
same block decomposition
$$
Y=\pmatrix Y_{11} & Y_{12} \\ Y_{21} & Y_{22} \endpmatrix
$$
by the formulas
$$
\align
Y_{11}&=(X_{11} - X_{12}X_{22}^{-1}X_{21})^{-1},\\
Y_{12}&=-X_{11}^{-1}X_{12} ( X_{22} - X_{21}X_{11}^{-1}X_{12})^{-1},\\
Y_{21}&=-X_{22}^{-1}X_{21} (X_{11} - X_{12}X_{22}^{-1}X_{21})^{-1},\\
Y_{22}&=(X_{22} - X_{21}X_{11}^{-1}X_{12})^{-1},
\endalign
$$
is the inverse to $X$.
\enddemo

Let $I,J$ be as in 1.2.1 and let $Y$ be the matrix inverse to
$X$, as in Proposition 1.2.1. Notice that each entry $y_{ij}$ of
$Y$ is a nonzero element of the division ring $F(\Cal X)$.

\definition{Definition 1.2.2 \rm(Quasideterminant of a matrix with formal
entries)} For $i\in I$, $j\in J$ the ($i,j)$-th
quasideterminant  $|X|_{ij}$ of $X$ is the element of $F(\Cal X)$ defined by
the formula
$$
|X|_{ij} = (y_{ji})^{-1}
$$
where $Y=X^{-1}=(y_{ij})$, see Proposition 1.2.1.
\enddefinition

From the proof of Proposition 1.2.1, we obtain the following
recurrence relations for $|X|_{ij}$.

First of all, if $n=1$, so that $I=\{i\}$, $J=\{j\}$, then $|X|_{ij}=x_{ij}$.

Next, let $n\geq 2$ and let $X^{ij}$ be the $(n-1)\times(n-1)$-matrix obtained
from $X$ by deleting the $i$-th row and the $j$-th column. Then
$$
|X|_{ij} = x_{ij} -\sum
x_{ii'}(|X^{ij}|_{j'i'})^{-1} x_{j'j}.\tag 1.2.1
$$
Here the sum is taken over $i'\in I\smallsetminus\{i\}$, $j'\in
J\smallsetminus\{j\}$.

\remark {Remark} In part (ii) of Definition 1.2.1, $X^{ij}$ is the matrix
with formal entries $x_{i'j'}$ indexed by elements $i'\in
I\smallsetminus\{i\}$, $j'\in J\smallsetminus\{j\}$, and
$(|X^{ij}|_{i'j'})^{-1}$ is the inverse of the quasideterminant
$|X^{ij}|_{i'j'}$ in the corresponding free division ring
$F( \Cal X')\subset F(\Cal X)$, where $\Cal X' = \{x_{i'j'},\,i'\in
I\smallsetminus \{i\},\,j'\in J\smallsetminus \{j\}\}$.
\endremark

\example{Examples 1.2.3} (a) For the $2\times 2$-matrix $X=(x_{ij})$, $i,j
=1,2$, there are four quasideterminants:
$$
\matrix
|X|_{11} = x_{11} - x_{12}\cdot x_{22}^{-1}\cdot x_{21},\quad
&|X|_{12}=x_{12}-x_{11}\cdot x_{21}^{-1} \cdot x_{22},\\
|X|_{21} = x_{21} - x_{22}\cdot x_{12}^{-1} \cdot x_{11},\quad &
|X|_{22}=x_{22}-x_{21}\cdot x_{11}^{-1} \cdot x_{12}.
\endmatrix
$$

(b) For the $3\times 3$-matrix $X=(x_{ij})$, $i,j=1,2,3$, there are 9
quasideterminants. One of them is
$$
\align
|X|_{11}=x_{11}&-x_{12}(x_{22}-x_{23}x_{33}^{-1}x_{32})^{-1}x_{21}
-x_{12}(x_{32}-x_{33}\cdot x_{23}^{-1} x_{22})^{-1} x_{31}
\\
& -x_{13}(x_{23}-x_{22}x_{32}^{-1}x_{33})^{-1}x_{21}
-x_{13}(x_{33}-x_{32}\cdot x_{22}^{-1}x_{23})^{-1}x_{31}.
\endalign
$$
\endexample

The action of the product of symmetric groups $S_n\times S_n$ on
$I\times J$, $|I|=|J|=n$, induces the action of $S_n\times S_n$ on the
the set of variables $\{a_{ij}\}$, $i\in I$, $j\in J$, and the corresponding
action on the free division ring $F(\Cal X)$. We denote this latter action
by $f\mapsto (\sigma, \tau)f$, $\sigma,\tau\in S_n$.

The following proposition shows that the definition of the quasideterminant is
compatible with this action.

\proclaim{Proposition 1.2.4} For $(\sigma,\tau)\in S_n\times S_n$ we have
$(\sigma,  \tau)\big( |X|_{ij}\big)=|X|_{ \sigma (i) \tau (j)}$.
\endproclaim

In particular, the stabilizer subgroup of $ |X|_{ij}$ under the action of
$S_n\times S_n$ is isomorphic to $S_{n-1}\times S_{n-1}$.

Proposition 1.2.4 shows that in the definition of the quasideterminant, we do
not need to require $I$ and $J$ to be ordered or a bijective correspondence
between $I$ and $J$ to be given.

We go now to the definition of quasideterminants over a ring $R$ with unit.
Let $A=(a_{ij})$, $i\in I$, $j\in J$, be a matrix over $R$. Such a matrix determines
the map $\alpha_A:\Cal X\to R$, $\Cal X = \{x_{ij}\}$,
given by the formula $\alpha_A(x_{ij})=a_{ij}$.

\definition{Definition 1.2.5 \rm(Quasideterminant of a matrix over a ring)}
Let $i\in I$, $j\in J$, and the formal quasideterminant $|X|_{ij}\in F(\Cal X)$
can be evaluated at $\alpha_A$ in the sense of Definition 1.1.4. Then we say
that the $(ij)$-th quasideterminant $|A|_{ij}$ of $A$ exists and is equal to
$\alpha_A(|X|_{ij})$.
Otherwise, we say that $|A|_{ij}$ does not exist.
\enddefinition

According to this definition, the quasideterminant $|A|_{ij}$ of a matrix $A$
over $R$ is an element of $R$.

According to Definition 1.2.2 and Proposition 1.2.4, the quasideteminant
$|A|_{ij}$ can be computed as  follows. Denote by $r_i^j$ the row submatrix of
length $n-1$ obtained from $i$-th row of $A$ by deleting the element $a_{ij}$,
and by $c_j^i$ the  column submatrix of height $n-1$ obtained from $j$-th
column of $A$ by deleting the element $a_{ij}$.

\proclaim{Proposition 1.2.6} Let $|I|, |J|>1$ and assume that the
$(n-1)\times(n-1)$-matrix  $A^{ij}$ is invertible over $R$. Then
$$
|A|_{ij}=a_{ij}-r_i^j(A^{ij})^{-1}c_j^i. \tag 1.2.2
$$
\endproclaim

\remark{Remark} For a generic matrix $A$, to find the quasideterminant
$|A|_{ij}$, one should take  the formula to $|X|_{ij}$, substitute
$x_{ij}\mapsto a_{ij}$, and verify that all inversions exist in $R$.  However,
in special cases (for example, when some of the entries of $A$ equal zero), one
might need to  replace the formula for the quasideterminant by an equivalent
formula and only then to substitute  $x_{ij}\mapsto a_{ij}$. Here is an
example.

Let
$$
A=\pmatrix a_{11}&a_{12}&a_{13}\\
 a_{21}&a_{22}&a_{23}\\
 0&a_{32}&a_{33}\endpmatrix
$$
where $a_{21}$ and $a_{32}$ are invertible in $R$.
The quasideterminant $|A|_{13}$ cannot be defined using formula
(1.2.1) since the rational expression
$a_{12}(a_{22}-a_{21}a_{31}^{-1}a_{32})^{-1}a_{23}$ is not defined.  However,
if we replace this expression in formula (1.2.1) by the equivalent expression
$a_{12}a_{32}^{-1}a_{31}(a_{22}a_{32}^{-1}a_{31}-a_{21})^{-1}a_{23}$, the new
formula is defined for the matrix $A$ and the corresponding rational function
given the quasideterminant $|A|_{13}$

Let us note also that the since the submatrix
$$
A^{13}=\pmatrix a_{21}&a_{22}\\
 0&a_{32}\endpmatrix
$$
is invertible, the quasideterminant $|A|_{13}$ can be defined using formula (1.2.2).
\endremark

Sometimes it is convenient to adopt a more graphic notation for the quasideterminant
by boxing the element $a_{ij}$. For $A=(a_{ij})$, $i,j=1,\dots , n$, we write
$$
|A|_{pq}=\left |\matrix
a_{11}&\dots &a_{1q}&\dots &a_{1n}\\
      &\dots &      &\dots &       \\
a_{p1}&\dots &\boxed {a_{pq}}&\dots &a_{pn}\\
      &\dots &      &\dots &       \\
a_{n1}&\dots &a_{nq}&\dots &a_{nn}\endmatrix \right |.
$$

If $A$ is a generic $n\times n$-matrix (in the sense that all square
submatrices of $A$ are invertible), then there exist $n^2$ quasideterminants of
$A$. However, a non-generic matrix may have $k$ quaisdeterminants, where $0\leq
k \leq n^2$. Example 1.2.3(a) shows that each of the quasideterminants
$|A|_{11}$, $|A|_{12}$, $|A|_{21}$, $|A|_{22}$ of a $2\times 2$-matrix $A$ is
defined whenever the corresponding element $a_{22}$, $a_{21}$, $a_{12}$,
$a_{11}$ is invertible.

\remark{Remark} The definition of the quasidereminant can be generalized to
define $|A|_{ij}$ for a matrix $A=(a_{ij})$ in which each $a_{ij}$ is an
invertible morphism $V_j\to V_i$ in an additive category $C$ and the matrix
$A^{pq}$ of morphisms is invertible. In this case the quasideterminant
$|A|_{pq}$ is a morphism from the object $V_q$ to the object $V_p$.
\endremark

The next example shows that the notion of a quasideterminant is not a
generalization of a determinant over a commutative ring, but rather a
generalization of a ratio of two determinants.

\example {Example} If the elements $a_{ij}$ of the matrix $A$
commute, then
$$
|A|_{pq} = (-1)^{p+q} \frac{\det A}{\det A^{pq}}.
$$
\endexample

We will show in Section 3 that similar expressions for
quasideterminants can be given for quantum matrices, quaternionic
matrices, Capelli matrices and other cases listed in the
Introduction.

In general quasideterminants are not polynomials in their entries,
but (non-commutative) rational functions.
The following theorem was conjectured by I. Gelfand and Retakh, and
proved by Reutenauer \cite {Re} in a slightly different form.

\proclaim{Theorem 1.2.7}  Quasideterminants of the $n \times n$-matrix
$X=(x_{ij})$ with formal entries have the inversion height $n-1$ over the free
division ring $F(\Cal X)$, $\Cal X=\{x_{ij}\}$.
\endproclaim

In the commutative case determinants are finite sums of monomials with
appropriate coefficients. As is shown in \cite {GR1, GR2}, in the
noncommutative  case quasideterminants of a matrix $X=(x_{ij})$ with formal
entries $x_{ij}$ can be identified with formal power series in the matrix
entries or their inverse. A simple example of this type is described below.

Let $X=(x_{ij})$, $i,j=1,\dots , n$, be a matrix with formal entries. Denote by
$E_n$ the identity matrix of order $n$  and by $\Gamma _n$ the complete
oriented graph with vertices $\{1,2, \dots , n\}$, with the arrow from $i$ to
$j$ labeled by $x_{ij}$. A path $p: i\rightarrow k_1\rightarrow k_2\rightarrow
\dots \rightarrow k_t\rightarrow j$ is labeled by the word
$w=x_{ik_1}x_{k_1k_2}x_{k_2k_3}\dots x_{k_tj}$.

Denote by $P_{ij}$ the set of words labelling paths going from $i$
to $j$, i.e. the set of words of the form
$w=x_{ik_1}x_{k_1k_2}x_{k_2k_3}\dots x_{k_tj}$. A simple path is
a path $p$ such that $k_s\neq i,j$ for every $s$. Denote by
$P'_{ij}$ the set of words labelling simple paths from $i$ to $j$.

Let $R$ be the ring of formal power series in $x_{ij}$ over a field. From
\cite{Co}, Section 4,  it follows that there is a canonical embedding of $R$ in
a division ring $D$ such that the image of $R$ generates $D$. We identify
$R$ with its image in $D$.

\proclaim{Proposition 1.2.8} Let $i, j$ be two distinct integers between $1$
and $n$. The
rational functions
$|I_n-X|_{ii}$, $|I_n-X|_{ij}^{-1}$ are defined in $D$ and
$$
|I_n-X|_{ii}=1-\sum _{w\in P'_{ii}}w, \qquad
|I_n-X|_{ij}^{-1}=\sum _{w\in P_{ij}}w.
$$
\endproclaim

\example {Example} For $n=2$,
$$|I_2-X|_{11}=1-x_{11}-\sum _{p\geq 0}x_{12}x_{22}^px_{21}.$$
\endexample

For some matrices of special form over a ring, quasideterminants can be expressed as
polynomials in the entries of the matrix.
The next proposition shows that this holds, in particular, for the so-called
almost triangular matrices. Such matrices play am important role in many
papers, including \cite{DS, Ko, Gi}.

\proclaim{Proposition 1.2.9} The following quasideterminant is a polynomial in
its entries:
$$
\left| \matrix
a_{11} & a_{12} & a_{13} & \dots & \boxed{a_{1n}}\\
-1         & a_{22} & a_{23} & \dots & a_{2n}\\
0           & -1        & a_{33}  & \dots & a_{3n}\\
             &             & \dots      &         &           \\
0           &             & \dots      & -1     & a_{nn}
\endmatrix \right |
= a_{1n} + \sum _{1\leq j_1< j_2< \dots <j_k <n}
a_{1j_1}a_{j_1+1,j_2}a_{j_2+1,j_3}\dots a_{j_k+1,n}.
$$
\endproclaim

\remark{Remark} Denote the expression on right-hand side by $P(A)$. Note that
$(-1)^{n-1}P(A)$ equals to the determinant of the almost upper-triangular
matrix over a commutative ring. For non-commutative almost upper triangular
matrices, Givental \cite{Gi} (and others) defined the determinant as $(-1)^{n-1}P(A)$.
\endremark

\example{Example} For $n=3$ we have
$$
P(A) = a_{13}+a_{11}a_{23} +a_{12}a_{33} +a_{11}a_{22}a_{33}.
$$
\endexample

\subhead 1.3. Transformation properties of quasideterminants\endsubhead
Let $A=(a_{ij})$ be a square matrix of order $n$ over a ring $R$.

(i)  The quasideterminant $|A|_{pq}$ does not depend on permutations
of rows and columns in the matrix $A$ that do not involve the $p$-th 
row and the $q$-th column. This follows from Proposition 1.2.3.

(ii)  {\it The multiplication of rows and columns.}  Let
the matrix $B=(b_{ij})$ be obtained from the matrix $A$ by
multiplying the $i$-th row by $\lambda\in R$  from
the left, i.e., $b_{ij}=\lambda a_{ij}$ and $b_{kj}=a_{kj}$ for
$k\neq i$. Then
$$
 |B|_{kj}=\cases \lambda |A|_{ij} \quad&\text{ if } k = i,\\
|A|_{kj} \quad&\text{ if } k \neq i \text{ and } \lambda \text{ is
invertible.}\endcases
$$

Let the matrix $C=(c_{ij})$ be obtained from the matrix A by
multiplying the $j$-th column by $\mu\in R$ from the
right, i.e. $c_{ij}=a_{ij}\mu $ and $c_{il}=a_{il}$ for all $i$
and $l\neq j$. Then
$$
|C|_{i\ell}=\cases |A|_{ij} \mu \quad&\text{ if } l = j,\\
|A|_{i\ell} \quad&\text{ if } l \neq j \text{ and } \mu \text{ is
invertible.}\endcases
$$

(iii) {\it The addition of rows and columns.} Let the matrix
$B$ be obtained from $A$ by replacing the $k$-th row of $A$ with the 
sum of the $k$-th and $l$-th rows, i.e., $b_{kj}=a_{kj}+a_{lj}$,
$b_{ij}=a_{ij}$ for $i\neq k$. Then           
$$
|A|_{ij} = |B|_{ij}, \qquad  i=1, \dots k-1, k+1,\dots n,
\quad j=1, \dots, n.
$$

Let the matrix
$C$ be obtained from $A$ by replacing the $k$-th column of $A$ with the 
sum of the $k$-th and $l$-th columns, i.e., $c_{ik}=a_{ik}+a_{il}$,
$c_{ij}=a_{ij}$ for $j\neq k$. Then           
$$
|A|_{ij}= |C|_{ij} , \qquad i=1,\dots, n,\quad,\dots ,\ell -1,\ell +1,\dots n.
$$

\subhead  1.4. General properties of quasideterminants
\endsubhead

\comment
\subhead Quasideterminants and inverse matrices \rm(\cite
{GR})\endsubhead 

\proclaim{Theorem 1.4.1} If $A=(a_{ij})$ is a generic matrix and
$B=(b_{ij})=A^{-1}$, then $b_{ji}=|A|_{ij}^{-1}$ for
all $i, j$.
\endproclaim
\endcomment

\subsubhead 1.4.1 Two involutions \rm(see \cite {GR4}) \endsubsubhead
For a square matrix $A=(a_{ij})$ over a ring $R$, denote by $IA=A^{-1}$ the
inverse matrix (if it exists), 
and by $HA=(a^{-1}_{ji})$ the Hadamard inverse matrix (a;so if it exists). It
is evident that if $IA$ exists, then $I^2A=A$, and if $HA$ exists, then 
$H^2A=A$.

Let $A^{-1}=(b_{ij})$. According to Theorem 1.2.1,
$b_{ij}=|A|^{-1}_{ji}$. This formula can be rewritten in the following form.

\proclaim{Theorem 1.4.1} For a square matrix A over a ring $R$,
$$
HI(A)=(|A|_{ij})\tag 1.4.1
$$
provided that all quasideterminants $|A|_{ij}$ exist.
\endproclaim

\subsubhead 1.4.2. Homological relations \rm (see \cite{GR}) \endsubsubhead Let
$X=(x_{ij})$ be a square matrix of order $n$ with formal entries. For 
$1\leq k,l \leq n$ let $X^{kl}$ be the submatrix of order $n-1$ of the matrix
$X$ obtained by deleting the $k$-th row and the $l$-th column.
Quasideterminants of the matrix $X$ and the submatrices are connected by the
following {\it homological relations.}

\proclaim{Theorem 1.4.2}
{\rm(i)}  Row homological relations:
$$
-|A|_{ij} \cdot |A^{i\ell}|^{-1}_{sj} = |A|_{i\ell}\cdot
|A^{ij}|^{-1}_{s\ell},\qquad  s\neq i
$$

{\rm(ii)}  Column homological relations:
$$
-|A^{kj}|^{-1}_{it} \cdot |A|_{ij} = |A^{ij}|^{-1}_{kt}\cdot |A|_{kj},
\qquad  t\neq j
$$
\endproclaim

The same relations hold for matrices over a ring $R$ provided the corresponding
quasideterminants exist and are invertible.

A consequence of homological relations is that the ratio of two
quasideterminants of an $n\times n$ matrix (each being a rational function of
inversion height $n-1$) actually equals a ration of two rational functions each 
having inversion height $< n-1$.

\subsubhead 1.4.3. Heredity \endsubsubhead
Let $A=(a_{ij})$ be an $n\times n$ matrix over a ring $R$, and let 
$$
A=\pmatrix A_{11} & \dots & A_{1s}\\
 {}&{} &{}\\
A_{s1} &\dots &A_{ss}\endpmatrix
\tag{1.4.2}
$$
be a block decomposition of $A$, where each $A_{pq}$ is a $k_p\times l_q$
matrix, $k_1+\dots+k_s=l_1+\dots+l_s=n$. Let us choose $p'$ and $q'$ such that
$k_{p'}=l_{q'}$, so that  $A_{p'q'}$ is a square matrix.

Let also $X=(x_{pq})$ be a matrix with formal variables and $|X|_{p'q'}$ be the
${p'q'}$-quasideterminant of $X$. In the formula for $|X|_{p'q'}$ as a rational
function in variable $x_{pq}$ we can substitute each variable $x_{pq}$ with the
corresponding matrix $A_{pq}$, obtaining a rational expression $F(A_{pq})$. Let
us note that  all matrix operations in this rational expression formally make
sense, i.e., in each addition, the orders of summands coincide, in each
multiplication, the number of columns of the first multiplier equals the number
of rows of the second multiplier, and each matrix that has to be inverted is a
square matrix. Let us assume that all matrices in this rational expression for 
that need to be inverted, are indeed invertible over $R$. Computing
$F(A_{pq})$, we obtain an $k_{p'}\times l_{q'}$ matrix over $R$, whose rows are
naturally numbered by indices
$$
i=k_1+\dots+k_{p'-1}+1,\,\dots,\, k_1+\dots+k_{p'}
\tag{1.4.3}
$$
and columns are numbered by indices
$$
j=l_1+\dots+l_{q'-1}+1,\,\dots,\, l_1+\dots+l_{q'}. 
\tag{1.4.4}
$$
We denote this matrix by $|X|_{p'q'}(A)$.

Let us note that under our assumptions, $k_{p'} = l_{q'}$, so that 
$|X|_{p'q'}(A)$ is a square matrix over $R$.

\proclaim{Theorem 1.4.3} Let the index $i$ lies in the range \rom{(1.4.3)} and
the index $j$ lies in the range \rom{(1.4.4)}. Let as assume that the matrix  
$|X|_{p'q'}(A)$ is defined. Then each of the quasideterminants $|A|_{ij}$ and 
$||X|_{p'q'}(A)|_{ij}$ exist if and only of the other exists, and in this case
$$
|A|_{ij}=||X|_{p'q'}(A)|_{ij}.
\tag{1.4.5}
$$
\endproclaim

\example{Example 1} Let in (1.4.2) $s=2$, $p'=q'=1$ and $k_1=l+1=1$. Then
formula (S0) becomes the inductive definition of the quasideterminant
$|A|_{ij}$ (see Definition 1.2.5).
\endexample

\example{Example 2} Let
$$
A=\pmatrix
a_{11}&a_{12}&a_{13}&a_{14}\\
a_{21}&a_{22}&a_{23}&a_{24}\\
a_{31}&a_{32}&a_{33}&a_{34}\\
a_{41}&a_{42}&a_{43}&a_{44} \endpmatrix.
$$
Take the decomposition $ A= \pmatrix  A_{11}&A_{12}\\A_{21}&A_{22}
\endpmatrix $ of $A$ into four $2\times 2$ matrices, so that $A_{11}=\pmatrix
a_{11}&a_{12}\\a_{21}&a_{22}\endpmatrix $, $A_{12}=
\pmatrix a_{13}&a_{14}\\a_{23}&a_{24}\endpmatrix $,
$A_{21}=\pmatrix a_{31}&a_{32}\\a_{41}&a_{42}\endpmatrix$, 
$A_{22}= \pmatrix a_{33}&a_{34}\\a_{43}&a_{44}\endpmatrix $. Let us use formula
(1.4.5) to find the quasideterminant $|A|_{13}$. 
We have
$$
\align
|X|_{12}(A)&=A_{12}-A_{11}A_{21}^{-1}A_{22}\\
&=
\pmatrix a_{13}&a_{14}\\a_{23}&a_{24}\endpmatrix
-\pmatrix a_{11}&a_{12}\\a_{21}&a_{22}\endpmatrix
{\pmatrix a_{31}&a_{32}\\a_{41}&a_{42}\endpmatrix}^{-1}
\pmatrix a_{33}&a_{34}\\a_{43}&a_{44}\endpmatrix
\\
&= 
\pmatrix a_{13}-\dots &a_{14}-\dots
\\a_{23}-\dots &a_{24}-\dots \endpmatrix. 
\endalign
$$
Denote the matrix in the right-hand side of this formula by
$\pmatrix c_{13}&c_{14}\\c_{23}&c_{24}\endpmatrix$. Then 
$$
|A|_{13}=\left| 
\matrix c_{13}&c_{14}\\c_{23}&c_{24}\endmatrix
\right|_{13},
$$
or, in other notation,
$$
|A|_{13}=\left| \matrix \boxed {c_{13}}
&c_{14}\\c_{23}&c_{24}\endmatrix \right| .
$$
\endexample

\subsubhead 1.4.4. A generalization of the homological relations\endsubsubhead 
Homological relations admit the following
generalization. For a matrix $A=(a_{ij})$, $i\in I$, $j\in J$, and two subsets
$L\subset I$, $M\subset J$ denote by $A^{L,M}$ the submatrix of the matrix $A$
obtained  by deleting the rows with the indexes $\ell \in L$ and the columns
with the indexes $m\in M$. Let $A$ be a square matrix, $L= (\ell_1,\dots,
\ell_k), M=(m_0,\dots,m_k)$. Set $M_i = M\smallsetminus\{m_i\}$, $i=0,\dots,
k$.

\proclaim{Theorem 1.4.4 \rm [GR1, GR2]}  For $p\notin L$ we have
$$
\align
\sum^k_{i=0}
|A^{L,M_i}|_{pm_i}\cdot|A|^{-1}_{\ell m_i}&= \delta_{p\ell},\\
\sum^k_{i=0}
|A|^{-1}_{m_i\ell }\cdot|A^{M_i,L}|_{m_ip} &= \delta_{\ell p},
\endalign
$$
provided the corresponding quasideterminants are defined and the matrices
$|A|^{-1}_{m_i\ell}$, $|A|^{-1}_{\ell m_i}$ are invertible over $R$.
\endproclaim

\subsubhead 1.4.5. Quasideterminants and Kronecker tensor products
\endsubsubhead Let $A=(a_{ij})$, $B=(b_{\alpha  \beta})$ be matrices over a
ring $R$. Denote by $C=A\otimes B$ the Kronecker tensor product, i.e., the
matrix with entries numbered by indices $(i\alpha,j\beta)$, and with the 
$(i\alpha,j\beta)$-th entry equal to $c_{i\alpha,j\beta}=a_{ij} b_{\alpha 
\beta}$.

\proclaim{Proposition 1.4.5} If quasideterminants $|A|_{ij}$ and
$|B|_{\alpha \beta}$ are defined, then the quasideterminant
$|A\otimes B|_{i\alpha, j \beta}$ is defined and
$$
|A\otimes B|_{i \alpha,j \beta}=|A|_{ij}|B|_{\alpha \beta}.
$$
\endproclaim

Note that in the commutative case the corresponding identity determinants is
different. Namely, if $A$ is a $m\times m$-matrix and $B$ is a $n\times
n$-matrix over a commutative ring, then $\det (A\otimes B)=(\det A)^n(\det
B)^m$.

\subsubhead 1.4.6. Quasideterminants and matrix rank \endsubsubhead
Let $A=(a_{ij})$ be a matrix over a division ring.

\proclaim{Proposition 1.4.6}  If the quasideterminant $|A|_{ij}$ is
defined, then the following statements are equivalent.
\roster
\item"{(i)}" $|A|_{ij} = 0,$
\item"{(ii)}"  the $i$-th row of the matrix $A$ is a left linear
combination of the other rows of $A$;
\item"{(iii)}"  the $j$-th column of the matrix $A$ is a right linear
combination of the other columns of $A$.
\endroster
\endproclaim

\example {Example} Let $i,j=1,2$ and $|A|_{11}=0$, i.e., 
$a_{11}-a_{12}a_{22}^{-1}a_{21}=0$. Therefore,
$a_{11}=\lambda a_{21}$, where $\lambda =a_{12}a_{22}^{-1}$.
Since $a_{12}=(a_{12}a_{22}^{-1})a_{22}$, the first row of
$A$ is proportional to the second row.
\endexample

There exists the notion of linear dependence for elements of a
(right or left) vector space over a division ring. 
So there exists the notion of the row rank (the dimension of the left
vector space spanned by the rows) and the notion of the
column rank (the dimension of the right vector space spanned by
the columns) and these
ranks are equal \cite {Ja, Co}. This also follows from
Proposition 1.4.6.

By definition, an  $r$-quasiminor of a square matrix $A$ is a
quasideterminant of an $r\times r$-submatrix of $A$.

\proclaim{Proposition 1.4.7}  The rank of the matrix $A$ over a division
algebra is $\geq r$ if and only if at least one 
$r$-quasiminor of the matrix $A$ is defined and is not equal to zero.
\endproclaim

\subhead 1.5. Basic identities
\endsubhead

\subsubhead 1.5.1. Row and column decomposition \endsubsubhead
The following result is an analogue of the classical expansion of a
determinant by a row or a column.

\proclaim{Proposition 1.5.1} Let $A$ be a matrix over a ring $R$. 
For each $k\neq p$ and each $\ell \neq q$ 
we have
$$ 
\align
|A|_{pq}&=a_{pq}-\sum_{j\neq q} a_{pj}(|A^{pq}|_{kj})^{-1}
|A^{pj}|_{kq},\\
|A|_{pq}&=a_{pq}-\sum_{i\neq p} |A^{iq}|_{pi}(|A^{pq}|_{i\ell})^{-1}
a_{iq},
\endalign
$$
provided all terms in right-hand sides of these expressions are defined.
\endproclaim

As it was pointed out in \cite {KL}, Propostiion 1.5.1 immediately follows from 
the homological relations (Theorem 1.4.2).

\subsubhead 1.5.2. Sylvester's identity\endsubsubhead 
  Let $A=(A_{ij}),
i,j=1,\dots, n$, be a matrix over a ring $R$ and 
$A_0=(a_{ij}), i,j=1, \dots, k$, a
submatrix of $A$ that is invertible over $R$.  For $p,q=k+1,\dots, n$ set
$$
c_{pq} = \vmatrix &&& a_{1q}\\
&A_0& & \vdots\\
&&&a_{kq}\\
a_{p1}&\dots& a_{pk} &a_{pq}\endvmatrix_{pq} .
$$
These quasidetrminants are defined because matrix $A_0$ is invertible.

Consider the $(n-k)\times (n-k)$ matrix
$$
C=(c_{pq}), \quad p,q = k+1,\dots, n.
$$
The submatrix $A_0$ is called the {\it pivot\/} for the matrix $C$.

\proclaim{Theorem 1.5.2 \rm(see \cite{GR})} For $i,j = k+1,\dots, n$,
$$
|A|_{ij} = |C|_{ij}
$$
\endproclaim

The commutative version of Theorem 1.5.2 is the following Sylvester's theorem.
 
\proclaim{Theorem 1.5.3}  Let $A=(a_{ij}),
i,j=1,\dots, n$, be a matrix over a commutative ring.  Suppose
that the submatrix $A_0=(a_{ij}), i,j=1,\dots, k$, of $A$ is
invertible.  For $p,q = k+1,\dots, n$ set
$$
\aligned
\tilde b_{pq} &= \det\pmatrix && & & a_{1q}\\
& A_0&& & \vdots\\
&&& & a_{kq}\\
&a_{p1}&\dots& a_{pk} &a_{pq}\endpmatrix ,\\
\tilde B&=(\tilde b_{pq}),\quad p,q= k+1,\dots, n.
\endaligned
$$
Then
$$
\det A={\det\tilde B\over (\det A_0)^{n-k-1}}.
$$
\endproclaim

\remark{Remark 1} A quasideterminant of an $n\times n$-matrix $A$ is
equal to the corresponding quasideterminant of a $2\times 2$-matrix
consisting of $(n-1)\times (n-1)$-quasiminors of the matrix $A$,
or to the quasideterminant of an $(n-1)\times(n-1)$-matrix consisting of
$2\times 2$-quasiminors of the matrix $A$. One can use any of these
procedures for an inductive definition of quasideterminants. In fact, 
Heyting \cite {H} essentially defined the quasideterminants $|A|_{nn}$ 
for matrices $A=(a_{ij})$, $i,j=1,\dots , n$, in this way.
\endremark

\remark{Remark 2} Theorem 1.5.2 can be generalized to
the case where $A_0$ is a square submatrix of $A$ formed by some (not
necessarily consecutive and not necessarily the same) rows and columns of $A$.
In particular, in the case where $A_0=(a_{ij}), i,j=2,\dots, n-1$, Theorem
1.5.2 is an analogue of a well-known
commutative identity which is called the ``Lewis Carroll
identity'' (see, for example, \cite {Ho}).
\endremark

\subsubhead 1.5.3. Inversion for quasiminors\endsubsubhead  The following theorem
was formulated in \cite {GR}. For a matrix $A=(a_{ij})$, $i\in I$, $j\in J$, over a ring $A$ and subsets 
$P\subset I$, $Q\subset J$ denote by $A_{PQ}$ the submatrix 
$$
A_{PQ} = (a_{\alpha\beta}),\qquad \alpha\in P,\quad \beta\in Q.
$$

Let $|I|= |J|$ and $B=A^{-1} = (b_{rs})$.  Suppose that $|P| = |Q|$.

\proclaim{Theorem 1.5.4}  Let $k\notin P, \ell\notin Q$.  Then
$$
|A_{P\cup\{k\}, Q\cup\{\ell\}}|_{k\ell}\cdot |B_{I\setminus P,
J\setminus Q}|_{\ell k} = 1.$$
\endproclaim
Set $P=I\smallsetminus\{k\}, Q= J\smallsetminus\{\ell\}$.  Then
this theorem leads to the already mentioned identity
$$
|A|_{k\ell} \cdot b_{\ell k} = 1.
$$

\example {Example} Theorem 1.5.4 implies the following identity for principal 
quasiminors.
Let $A=(a_{ij})$, $i,j=1,\dots , n$ be an invertible matrix over $R$ and 
 $B=(b_{ij})=A^{-1}$. For a fixed $k$, $1\leq k\leq n$, 
set $A_{(k)}=(a_{ij})$, $i,j=1,\dots , k$ and  $B^{(k)}=(b_{ij})$, 
$i,j=k,\dots , n$. Then
$$
|A_{(k)}|_{kk}\cdot |B^{(k)}|_{kk}=1.
$$
\endexample

\subsubhead 1.5.4. Multiplicative properties of quasideterminants
\endsubsubhead 
Let $X =(x_{pq}), Y=(y_{rs})$ be $n\times n $-matrices. The following statement
follows directly from from Definition 1.2.2.

\proclaim{Theorem 1.5.5} We have 
$$|XY|^{-1}_{ij} = \sum^n_{p=1}|Y|^{-1}_{pj} |X|^{-1}_{ip}.
$$
\endproclaim

\subhead 1.6. Noncommutative linear algebra\endsubhead
In this section we use quasideterminants to noncommutative generalizations of
basic theorems about systems of linear equations (see \cite {GR, GR1, GR2}).

\subsubhead 1.6.1. Solutions of systems of linear equations\endsubsubhead

\proclaim{Theorem 1.6.1} Let $A=(a_{ij})$ be an $n\times n$ matrix over a 
ring $R$. Assume that all the quasideterminants $|A|_{ij}$ are defined and 
invertible. Then
$$
\cases
a_{11} x_1+\dots + a_{1n} x_n = \xi_1\\
\qquad\qquad \dots \\
a_{n1} x_1+\dots + a_{nn} x_n = \xi_n
\endcases
$$
for some $x_i\in R$ if and only if
$$
x_i = \sum^n_{j=1} |A|^{-1}_{ji} \xi_j.\qquad i=1,\dots, n.
$$
\endproclaim

\subsubhead 1.6.2. Cramer's rule \endsubsubhead  Let $A_{\ell}(\xi)$ be the
$n\times n$-matrix obtained by replacing the $\ell$-th column of the
matrix $A$ with the column $(\xi_1,\dots, \xi_n)$.

\proclaim{Theorem 1.6.2}  In notation of Theorem 1.6.1, if the
quasideterminants $|A|_{ij}$ and
$|A_j(\xi)|_{ij}$ are defined, then
$$
|A|_{ij} x_j = |A_j(\xi) |_{ij}.
$$
\endproclaim

\subsubhead 1.6.3. Cayley--Hamilton theorem\endsubsubhead  
Let $A=(a_{ij})$, $i,j=1,\dots, n$, be a matrix over a ring $R$. 
Denote by $E_n$ the identity matrix of order $n$. 

Let $t$ be a formal variable. Set $f_{ij} = |tE_n-  A|_{ij}$ for $1\leq i,j\leq
n$. Then $f_{ij}(t)$ is a rational function in $t$. Define  the matrix function
$\tilde f_{ij}(t)$ by replacing  in $f_{ij}(t)$ each element $a_{ij}$ with the
matrix $ \tilde a_{ij} = a_{ij} E_n$ of order $n$ and the variable $t$ by the
matrix $A$. The functions $f_{ij}(t)$ are called the characteristic functions
of the matrix $A$. 

The following theorem was stated in \cite{GR1, GR2}. 

\proclaim{Theorem 1.6.3}  $\tilde f_{ij}(A) = 0$ for all $i,j=1,\dots, n$.
\endproclaim

\head 2. Important example: quaternionic quasideterminants
\endhead

As an example, we compute here quasideterminants  of quaternionic
matrices.

\subhead 2.1. Norms of quaternionic matrices
\endsubhead 
Let $\Bbb H$ be the algebra of quaternions. Algebra $\Bbb H$
is an algebra over the field of real numbers $\Bbb R$ with
generators ${\bold i}, {\bold j}, {\bold k}$ such that
${\bold i}^2={\bold j}^2={\bold k}^2=-1$ and
${\bold i}{\bold j}={\bold k}$,
${\bold j}{\bold k}={\bold i}$,
${\bold k}{\bold i}={\bold j}$.
It follows from the definition that
${\bold i}{\bold j}+{\bold j}{\bold i}=0$,
${\bold i}{\bold k}+{\bold k}{\bold i}=0$,
${\bold j}{\bold k}+{\bold k}{\bold j}=0$.

Algebra $\Bbb H$ posseses a standard anti-involution
$a\mapsto \bar a$: if
$a=x+y{\bold i}+z{\bold j}+t{\bold k}$,
$x,y,z,t\in \Bbb R$, then
$\bar a=x-y{\bold i}-z{\bold j}-t{\bold k}$. It
follows that $a\bar a=x^2+y^2+z^2+t^2$.
The multiplicative functional $\nu: \Bbb H \to
\Bbb R_{\geq 0}$ where $\nu (a)=a\bar a$ is called
the norm of $a$. One can see that
$a^{-1}=\frac {\bar a}{\nu (a)}$ for $a\neq 0$.

We will need the following generalization of the norm
$\nu $ to quaternionic matrices.
Let $M(n, \Bbb H)$ be the $\Bbb R$-algebra of quaternionic
matrices of order $n$. There exists a unique multiplicative functional
$\nu : M(n, \Bbb H)\rightarrow \Bbb R_{\geq 0}$ such that

(i) $\nu (A)=0$ if and only if the matrix $A$ is non-invertible,

(ii) If $A '$ is obtained from $A $ by adding a left-multiple of
a row to another row or a right-multiple of a column to another
column, then $\nu (A ')=\nu (A)$.

(iii) $\nu (E_n)=1$ where $E_n$ is the identity matrix of order $n$.

The number $\nu(A)$ is called the {\it norm\/} of the quaternionic matrix $A$.

For a quaternionic matrix $A=(a_{ij})$, $i,j=1,\dots , n$, denote by $A^*=(\bar a_{ji})$
the conjugate matrix. It is known that $\nu (A)$ coincides
with the Dieudonne determinant of $A$ and with the Moore determinant
of $AA^*$ (see \cite {As} and Subsections 3.2--3.4 below). The norm
$\nu (A)$ is a real number and it is equal to an alternating sum of monomials
of  order $2n$ in the $a_{ij}$ and $\bar a_{ij}$. An expression for
$\nu (A)$ is given by Theorem 2.1.2 below.

Let $A=(a_{ij})$, $i,j=1,\dots , n$, be a quaternionic matrix. Let
$I=\{i_1,\dots , i_k\}$, $J=\{j_1,\dots , j_k\}$ be two ordered sets of natural numbers
such that all $i_p$ and all $j_p$ are distinct. Set
$$z_{I,J}=a_{i_1j_1}\bar a_{i_2j_1}a_{i_2j_2}\dots a_{i_kj_k}\bar a_{i_1j_k}.$$

Denote by $\mu _i(A)$ the sum of all $z_{I,J}(A)$ such that $i_1=i$. One can easily see
that $\mu _i(A)$ is a real number since with each monomial  $z_{I,J}$ 
it contains 
the conjugate monomial $\overline z_{I,J} = z_{I',J'}$, where
$I'=\{i_1,i_k,i_{k-1}\dots,i_2\}$, $J=\{j_k,j_{k-1},\dots , j_1\}$. 

\proclaim{Proposition 2.1.1} The sum $\mu _i(A)$ does not depend on $i$.
\endproclaim

\example{Example} For  $n=1$ the statement is obvious. For $n=2$ we have
$$
\align
\mu _1(A)&=a_{11}\bar a_{21}a_{22}\bar a_{12}+a_{12}\bar a_{22}a_{21}\bar
a_{11},\\ 
\mu _2(A)&=a_{22}\bar a_{12}a_{11}\bar a_{21}+a_{21}\bar a_{11}a_{12}\bar
a_{22}.
\endalign
$$ 
Note that for two arbitrary quaternions $x,y$ we have $xy+\bar y\bar x=2\Re
(xy)=2\Re (yx)=yx+\bar x\bar y$, where $\Re (a)$ is the real part of the
quaternion $a$.  By setting
 $x=a_{11}\bar a_{21}$, $y=a_{22}\bar a_{12}$ one see that $\mu _1(A)=\mu
_2(A)$. 
\endexample

Proposition 2.1.1 shows that we may omit the index $i$ in $\mu _i(A)$ and denote it by
$\mu (A)$. 

Let $A=(a_{ij})$, $i,j=1,\dots , n$ be a matrix. We call an (unordered)
 set of square submatrices $\{A_1,\dots , A_s\}$
where $A _p=(a_{ij})$,  $i\in I_p$, $j\in J_p$  a {\it complete set} if  $I_p\cap I_q=
J_p\cap J_q=\emptyset $ for all $p\neq q$ and $\cup _pI_p=\cup _pJ_p=\{1,\dots
, n\}$. 

\proclaim{Theorem 2.1.2} Let $A=(a_{ij})$, $i,j=1,\dots , n$ be a quaternionic
matrix. Then 
$$
\nu (A)= \sum  (-1)^{k_1+\dots+k_p-p}\mu (A_1)\dots  \mu (A_p),
$$
where the sum is taken over all complete sets $(A_1,\dots , A_p)$ of
submatrices of $A$, $k_i$ is the order of the matrix $A_i$. 
\endproclaim

\example {Example} For $n=2$ we have
$$\nu (A)= \nu (a_{11})\nu (a_{22})+ \nu (a_{12})\nu (a_{21}) -(a_{11}\bar
a_{21}a_{22}\bar a_{12}+a_{12}\bar a_{22}a_{21}\bar a_{11}).$$ 
\endexample

\proclaim{Corollary 2.1.3} Let $A$ be a square quaternionic matrix. Fix an
arbitrary $i\in \{1,\dots , n\}$. Then
$$
\nu (A)=\sum (-1)^{k(B_1)-1}\nu (B_1)\mu (B_2) 
$$
where the sum is taken over all complete sets of submatrices $(B_1, B_2)$ such
that $B_2$ contains an element from the $i$-th row, $k(B_1)$ the order of $B_1$, and
$\nu (B_1)=1$ if $B_2=A$. 
\endproclaim

\subhead 2.2. Quasideterminants of quaternionic matrices
\endsubhead
This section contains results from \cite {GRW1}.

Let $A=(a_{ij})$ ,  $i,j=1,\dots , n$, be a quaternionic matrix.
Let $I=\{i_1,\dots , i_k\}$ and $J=\{j_1,\dots , j_k\}$ be two ordered sets
of natural numbers $1\leq
i_1, i_2,\dots , i_k\leq n$ and $1\leq j_1,j_2,\dots , j_k\leq n$
such that all $i_p$ are distinct and all $j_p$ are distinct.
For $k=1$ set  $m_{I,J}(A)=a_{i_1j_1}$. For $k\geq 2$ set
$$
m_{I,J}(A)=a_{i_1j_2}\bar a_{i_2j_2}a_{i_2j_3}\bar a_{i_3j_3}
a_{i_3j_4}\dots \bar a_{i_kj_k}a_{i_kj_1}.
$$

If the matrix $A$ is Hermitian, i.e., $a_{ji}=\bar a_{ij}$ for all
$i,j$, then
$$
m_{I,J}(A)=a_{i_1j_2}a_{j_2i_2}a_{i_2j_3} a_{j_3i_3}
a_{i_3j_4}\dots a_{j_ki_k}a_{i_kj_1}.
$$

To a quaternionic  matrix $A=(a_{pq})$, $p,q=1,\dots , n$, and to
a fixed row index $i$ and a column index $j$ we associate a
polynomial in $a_{pq}$, $\bar a_{pq}$, which we call the
$(i,j)$-th double permanent of $A$. 

\proclaim {Definition 2.2.1}
The $(i,j)$-th double permanent of $A$ is the sum
$$
\pi _{ij}(A)=\sum  m_{I,J}(A),
$$
taken over all orderings $I=\{i_1,\dots , i_n\}$,  $J=\{j_1,\dots, j_n\}$ of
$\{1,\dots , n\}$ such that $i_1=i$ and  $j_1=j$ . 
\endproclaim

\example{Example} For $n=2$
$$
\pi _{11}(A)=a_{12}\bar a_{22}a_{21}.
$$
For $n=3$
$$
\pi _{11}(A)=a_{12}\bar a_{32}a_{33}\bar a_{23}a_{21}+ a_{12}\bar
a_{22}a_{23}\bar a_{33}a_{31}+ a_{13}\bar a_{33}a_{32}\bar
a_{22}a_{21}+ a_{13}\bar a_{23}a_{22}\bar a_{32}a_{31}.
$$
\endexample

For a submatrix $B$ of $A$ denote by $B^c$ 
the matrix obtained from $A$ by deleting all
rows and columns containing elements from $B$. If $B$ is a
$k\times k$-matrix, then $B^c$ is a $(n-k)\times (n-k)$-matrix.
$B^c$ is called the complementary submatrix of $B$.

Quasideterminants of a matrix $A=(a_{ij})$ are rational 
functions of elements $a_{ij}$. Therefore, 
for a quaternionic matrix $A$, its quasideterminants are 
polynomials in $a_{ij}$ and their conjugates, with coefficients that 
are rational functions of  $a_{ij}$ always taking rational values.
The following theorem gives expressions for these polynomials. 

\proclaim
{Theorem 2.2.2} If the quasideterminant $|A|_{ij}$ of a quaternionic
matrix is defined, then
$$
\nu ( A^{ij}) |A|_{ij}= \sum (-1)^{k(B)-1}\nu ( B^c)\pi _{ij}(B)
\tag 2.2.1
$$
where the sum is taken over all square submatrices $B$ of $A$ containing
$a_{ij}$, $k(B)$ is the order of $B$, and we set $\nu (B^c)=1$ for $B=A$.
\endproclaim

Recall that according to Proposition 1.2.6 the quasideterminant 
$|A|_{ij}$ is defined if the matrix
$A^{ij}$ is invertible. In this case $\nu (A^{ij})$ is invertible,
so that formula  (2.2.1) indeed gives an expression for $|A|_{ij}$.

The right-hand side in (2.2.1) is a linear combination  with real coefficients
of monomials of lengths $1,3,\dots , 2n-1$ in $a_{ij}$ and $\bar a_{ij}$. The
number $\mu (n)$  of such monomials for a matrix of order $n$ is $\mu
(n)=1+(n-1)^2\mu (n-1)$.

\example{Example} For $n=2$
$$
\nu ( a_{22}) |A|_{11}=\nu ( a_{22})a_{11}-a_{12}\bar
a_{22}a_{21}.
$$
For $n=3$
$$
\align 
\nu (A^{11})|A|_{11}={}&\nu (A^{11})a_{11}
-\nu (a_{33})a_{12}\bar a_{22}a_{21}- \nu
(a_{23})a_{12}\bar a_{32}a_{31}-
\\
&-\nu (a_{32})a_{13}\bar a_{23}a_{21}
-\nu (a_{22})a_{13}\bar a_{33}a_{31}+
a_{12}\bar a_{32}a_{33}\bar a_{23}a_{21}+\\
&+ a_{12}\bar a_{22}a_{23}\bar a_{33}a_{31}
+ a_{13}\bar a_{33}a_{32}\bar a_{22}a_{21}
+ a_{13}\bar a_{23}a_{22}\bar a_{32}a_{31}.
\endalign
$$
\endexample

The example shows how to simplify the general formula for
quasideterminants of matrix of order $3$ (see subsection 1.2) for
quaternionic matrices.

The  following
theorem, which is similar to Corollary 2.1.3, shows that the
coefficients in formula (2.2.1) are uniquely defined.
\proclaim {Theorem 2.2.3} Let quasideterminants $|A|_{ij}$ of quaternionic
matrices 
are given by the formula
$$
\xi (A^{ij})|A|_{ij}=\sum (-1)^{k(B)-1}\xi (B^c)\pi _{ij}(B)
$$
and all coefficients $\xi (C)$ depend of submatrix $C$ only, then
$\xi (C)=\nu (C)$ for all square matrix $C$.
\endproclaim

\example {Example} For $n=2$ set $a_{11}=0$. Then $\xi
(a_{22})a_{12}a_{22}^{-1}a_{21}=a_{12}\bar a_{22}a_{21}$. This
implies that $\xi (a_{22})=\bar a_{22}a_{22}=\nu (a_{22})$.
\endexample
\head 3. Noncommutative determinants \endhead

Noncommutative determinants were defined in different and, sometimes, not
related situations. In this  section we present some results from \cite {GR,
GR1, GR2, GRW1} describing a universal approach to noncommutative
determinants and norms of noncommutative matrices based on the notion of
quasideterminants.

\subhead 3.1. Noncommutative determinants as products of
quasiminors \endsubhead
Let $A=(a_{ij})$, $i,j=1,\dots , n$, be a matrix over a
division ring $R$ such that all square submatrices of $A$
are invertible. 
For $\{i_1,\dots , i_k\},\{j_1,\dots , j_k\}\subset\{1,\dots , n\}$
define $A^{i_1\dots
i_k, j_1\dots j_k}$ to be the submatrix of $A$ obtained by deleting
rows with indices $i_1,\dots , i_k$ and columns with indices
$j_1,\dots , j_k$. Next, for any orderings $I=(i_1,\dots , i_n)$,
$J=(j_1,\dots , j_n)$ of $\{1,\dots , n\}$  set
$$
D_{I,J}(A)=|A|_{i_1 j_1}|A^{i_1 j_1}|_{i_2 j_2}|A^{i_1i_2,
j_1j_2}|_{i_3 j_3} \dots a_{i_n j_n}.
$$

In the commutative case $D_{I,J}(A)$ is, up the the sign, the determinant of
$A$. When $A$ is a quantum matrix  $D_{I,J}(A)$ differs from the quantum
determinant of $A$ by a factor depending on $q$ \cite {GR, GR1, KL}. The same
is true for some other  noncommutative algebras. This suggests to call
$D_{I,J}(A)$ the $(I,J)$-predeterminants of $A$. From the ``categorical point of
view" the expressions $D_{I,\tilde I}(A)$ where $I=(i_1, i_2,\dots , i_n)$, $\tilde I=(i_2,
i_3,\dots , i_n, i_1)$ are particularly important. We denote
$D_I(A)=D_{I, \tilde I}(A)$.  It is also convenient to have the basic predeterminant
$$
\Delta (A)=D_{\{12\dots n\},\{23\dots n1\}}.
\tag3.1.1
$$

We use the homological relations for quasideterminants to compare
different $D_{I,J}$. Here we restrict ourselves to elementary transformations
of  $I$ and $J$.

Let $I=(i_1, \dots , i_p, i_{p+1},\dots , i_n)$ and   $J=(j_1, \dots , j_p,
j_{p+1},\dots , j_n)$. Set  $I'=(i_1, \dots , i_{p+1}, i_p,\dots , i_n)$, 
$J'=(j_1, \dots , j_{p+1}, j_p,\dots , j_n)$.
Set also
$$
\align
X&=|A|_{i_1, j_1}|A^{i_1, j_1}|_{i_2, j_2}\dots |A^{i_1\dots
i_{p-2}, j_1,\dots , j_{p-2}}|_{i_{p-1}, j_{p-1}},\\
Y&=|A^{i_1\dots i_{p+1}, j_1,\dots , j_{p+1}}|_{i_{p+2}, j_{p+2}}
\dots a_{i_n, j_n},\\
u&=|A^{i_1\dots i_p, j_1,\dots , j_p}|_{i_{p+1}, j_{p+1}},\\
w_1&=|A^{i_1\dots i_{p-1}i_{p+1}, j_1,\dots , j_p}|_{i_p,
j_{p+1}},\\ 
w_2&=|A^{i_1,\dots i_p, j_1,\dots ,j_{p-1}}|_{i_{p+1}, j_{p+1}}.
\endalign
$$

\proclaim {Proposition 3.1.1} We have
$$
\align
D_{I, J'}&=-D_{I, J}Y^{-1}u^{-1}w_2^{-1}uw_2Y,
\\
D_{I', J}&=-Xuw_1^{-1}X^{-1}D_{I, J}Y^{-1}u^{-1}w_1Y.
\endalign
$$
\endproclaim

Let  $C$ be a commutative ring with a unit and
$f: R\rightarrow C$ be a multiplicative map, i.e. $f(ab)=f(a)f(b)$ for all
$a,b\in R$. 

Let $I=(i_1,\dots , i_n)$, $J=(j_1,\dots , j_n)$ be any
orderings of $(1,\dots , n)$. For an element $\sigma $ from the
symmetric group of $n$-th order set $\sigma (I)=
(\sigma (i_1),\dots , \sigma (i_n))$. Let $p(\sigma )$ be the
parity of $\sigma $.

Proposition 3.1.1  immediately implies the following theorem.

\proclaim {Theorem 3.1.2} In notations of Section {\rm 3.1}
we have
$$
f(D_{I,J}(A))=f(-1)^{p(\sigma _1)+p(\sigma _2)}f(D_{\sigma (I),
\sigma (J)}(A)).
$$
\endproclaim

It follows that $f(D_{I,J}(A)$ is uniquely defined up to a
power of $f(-1)$. We call $f(D_{1\dots n, 1\dots n})(A))$
the $f$-determinant $A$ and denote it by $fD(A)$. Note that if $f$ is
a homomorphism then $f$-determinant $fD(A)$ equals to the usual determinant of
the commutative matrix $f(A)$. 

\proclaim {Corollary 3.1.3} We have
$$fD(AB)=fD(A)\cdot fD(B).$$
\endproclaim

When $R$ is the algebra of quaternions and $f(a)=\nu (a)=a\bar a$, or,
in other words, $f$ is the quaternionic norm, then one can see that $fD(a)$ is
the matrix quaternionic norm $\nu (A)$ (see Section 2.1).

In Theorems 3.1.4--3.1.6 we present formulas for determinants of triangular and
almost triangular matrices. A matrix 
$A=(a_{ij})$, $i,j=1,\dots , n$, is called an upper-triangular matrix
if $a_{ij}=0$ for $i>j$. An upper-triangular matrix $A$ is called a generic
upper-triangular matrix if every square submatrix $A$ consisting of the rows 
$i_1\le i_2\le \dots \le i_k$ and the columns $j_1\le j_2\le \dots \le j_k$ 
such that $i_1\le j_1$, $i_2\le j_2$, \dots, $i_k\le j_k$, is invertible.

\proclaim{Theorem 3.1.4}
Let $A=(a_{ij})$, $i,j=1,\dots , n$, be a generic upper-triangular matrix.
The determinants $D_{i_1i_2\dots i_n}(A)$ are defined if and only if
$i_1=n$. In this case
$$
\multline
D_{ni_2\dots i_{n-1}}(A)=a_{nn}\cdot |A^{n,i_2}|^{-1}_{i_2n}\cdot
a_{i_2i_2}\cdot |A^{n,i_2}|_{i_2n}\cdot |A^{ni_2,i_2i_3}|^{-1}_{i_3n}\cdot
a_{i_3i_3}|A^{ni_2,i_2i_3}|_{i_3n}\cdot \dots 
\\
\cdot |A^{ni_2i_3\dots i_{n-1},i_2i_3\dots i_n}|^{-1}_{i_nn}\cdot
a_{i_ni_n}\cdot |A^{ni_2i_3\dots i_{n-1},i_2i_3\dots i_n}|_{i_nn}.
\endmultline
$$
\endproclaim

In particular,
$$
D_{n,n-1\dots 2,1}(A)=a_{nn}a_{n-1,n}^{-1}a_{n-1,n-1}a_{n-1,n}\dots
a_{1n}^{-1}a_{11}a_{1n}.
$$

A matrix $A=(a_{ij})$, $i,j=1,\dots , n$, is called an almost
upper-triangular matrix if $a_{ij}=0$ for $i>j+1$. An almost
upper-triangular matrix $A$ is called a
Frobenius matrix
if $a_{ij}=0$ for all $j\neq n$ and $i\neq j+1$, and  $a_{j+1j}=1$ for
$j=1,\dots , n-1$.

\proclaim{Theorem 3.1.5} If $A$ is invertible upper-triangular matrix, then
$$
D_{1,n,n-1\dots 2}(A)=|A|_{1n}a_{n,n-1}a_{n-1,n-2}\dots a_{21}.
$$
\endproclaim

By Proposition 1.2.7, the determinant $D_{1,n,n-1\dots 2}(A)$ of an 
uppen-triangular matrix $A$ is polynomial in $a_{ij}$.

Let $p(I)$ be the parity of the ordering $I=(i_1,\dots , i_n)$.

\proclaim{Theorem 3.1.6} If $A$ is a Frobenius matrix and the determinant
$D_I(A)$ is defined, then $D_{I}(A)=(-1)^{p(I)+1}a_{1n}$.
\endproclaim    

Now let $R$ be a division ring, $R^*=R\smallsetminus \{0\}$ the monoid of
invertible elements in $R$ and $\pi:R^*\rightarrow R^*/[R^*,R^*]$ the
canonical homomorphism. 
To the abelian group $R^*/[R^*,R^*]$ we adjoin the zero element $0$ with
obvious multiplication, and denote the obtained semi-group by $\tilde R$.
Extend $\pi $ to a map $R \rightarrow \tilde  R$ by setting
$\pi (0)=0$.
 
We recall here the classical notion of the Dieudonne determinant
(see \cite {D, A}). There exists a unique homomorphism    
$$ 
\det : M_n(R)\rightarrow \tilde R 
$$ 
such that

(i) $\det A'=\tilde \mu \det A $ for any matrix $A'$ obtained from $A\in M_n(R)$
by multiplying one row of $A$ from the left by $\mu $;

(ii) $\det A''=\det A $ for any matrix $A''$ obtained from $A$ by
adding one row to another;

(iii) $\det(E_n)=1$ for the identity matrix $E_n$.

The homomorphism $\det $ is called the Dieudonne determinant.  

It is known that $\det A=0$ if $\text{rank}(A)<n$ (see \cite{A}, Chapter 4).
The next proposition gives a construction of the Dieudonne determinant in the
case where $\text{rank}(A)=n$.

\proclaim{Proposition 3.1.7}
Let $A$ be an $n\times n$-matrix over a division ring $R$. If
$\text{rank}(A)=n$, then 

{\rm (i)} There exist orderings $I$ and $J$ of $\{1,\dots,n\}$ such that
$D_{I,J}(A)$ is defined. 

{\rm (ii)} If  $D_{I,J}(A)$ is defined, then the Diedonne determinant is given
by the formula $\det A = p(I)p(J) \pi(D_{I,J}(A)) $, where $p(I)$ is the parity
of the ordering $I$. 
\endproclaim

Note that in \cite {Dr} Draxl introduced the Dieudonne
predeterminant, denoted $\delta \epsilon \tau $. For a generic
matrix $A$ over a division ring there exists the Gauss
decomposition $A=UDL$ where $U, D, L$ are upper-unipotent,
diagonal, and lower-unipotent matrices. Then Draxl $\delta
\epsilon \tau  (A)$ is defined as the product
of diagonal elements in $D$ from top to the bottom.
For nongeneric matrices Draxl used the Bruhat decomposition
instead of the Gauss decomposition.

\proclaim{Proposition 3.1.8} $\delta \epsilon \tau  (A)=\Delta (A)$, where
$\Delta(A)$ is given by {\rm(3.1.1)}.
\endproclaim

\demo{Proof \rm(for a generic $A$)} Let $y_1, \dots , y_n$ be the diagonal
elements in $D$ from top to the bottom. As shown in \cite {GR1, GR2} (see also
4.9), $y_k=|A^{12\dots k-1, 12\dots k-1}|_{kk}$. Then $\delta \epsilon \tau 
(A)=y_1y_2\dots y_n=\Delta (A)$. \qed
\enddemo

Below we consider below special examples of noncommutative determinants.

\subhead 3.2. Dieudonne determinant for quaternions \endsubhead
Let $A=(a_{ij})$, $i,j=1,\dots , n$, be a quaternionic
matrix. If $A$ is not inverstible, then the Dieudonne determinant of $A$ equals
zero. By Proposition 3.1.7, if $A$ is invertible, there exist  orderings
$I=(i_1,\dots , i_n)$, $J=(j_1,\dots , j_n)$ of $\{1,\dots , n\}$ such that the
following expressions are defined: 
$$
D_{I,J}(A)=|A|_{i_1 j_1}|A^{i_1 j_1}|_{i_2, j_2}|A^{i_1i_2,
j_1j_2}|_{i_3 j_3} \dots a_{i_n j_n}.
$$
By Theorem 2.2.2, $D_{I,J}(A)$ can be expressed as a polynomial in $a_{ij}$ and 
$\overline{a_{ij}}$ with real coefficients.

In the quaternionic case the Dieudonne determinant $D$ coincides with the map
$$
\det: M_n(\Bbb H)\rightarrow \Bbb R_{\geq 0}
$$
(see \cite {As}).

The following proposition generalizes a result in \cite {VP}.

\proclaim {Proposition 3.2.1} In the quaternionic case
for each $I,J$ we have 
$$
\det A = \nu(D_{I,J}(A))^{1/2}
$$ 
{\rm(}the positive square root{\rm)}.
\endproclaim
The proof of  Proposition 3.2.1 follows from the homological relations
for quasideterminants.

\subhead 3.3. Moore determinants of Herimitian quaternionic
matrices \endsubhead
A quaternionic matrix $A=(a_{ij})$, $i,j=1,\dots , n$, is
called Hermitian if $a_{ji}=\bar a_{ji}$ for all $i,j$. It follows that
all diagonal elements of $A$ are real numbers and that the
submatrices $A^{11}$,   $A^{12,12}$, $\dots $ are Hermitian.

The notion of determinant for Hermitian quaternionic matrices was introduced by
E.~M.~Moore in 1922 \cite {M, MB}. Here is the original definition.

Let $A=(a_{ij})$, $i,j=1,\dots , n$, be a matrix over a ring. Let
$\sigma $ be a permutation of $\{1,\dots , n\}$. Write $\sigma $
as a product of disjoint cycles.
Since disjoint cycles commute, we may write
$$
\sigma =(k_{11}\dots k_{1j_1})(k_{21}\dots k_{2j_2})\dots
(k_{m1}\dots k_{mj_m})
$$
where for each $i$, we have $k_{i1}<k_{ij}$ for all $j > 1$, and
$k_{11} > k_{21}> \dots > k_{m1}$. This expression is unique. Let
$p(\sigma )$ be the parity of $\sigma $. The Moore determinant $M(A)$
is defined as follows:
$$
M(A)=\sum _{\sigma \in S_n}p(\sigma)a_{ k_{11}, k_{12}}\dots  a_{
k_{1j_1}, k_{11}} a_{ k_{21}, k_{22}}\dots a_{ k_{mj_m}, k_{m1}}.
\tag 3.3.1
$$
(There are equaivalent formulations of this definition; e.g.,
one can require $k_{i1}>k_{ij}$ for all $j>1$.)
If $A$ is Hermitian quaternionic matrix then $M(A)$ is a real
number. Moore determinants have nice features and are widely used
(see, for example, {\cite {Al, As, Dy1}).

We will show (Theorem 3.3.2) that determinants of Hermitian quaternionic
matrices can be obtained using our general approach. First we prove that for a
quaternionic Hermitian matrix $A$, the determinants $D_{I,I'}(A)$ coincide up
to a sign.

Recall that $\Delta (A)=D_{I,I'}(A)$ for $I=\{1,\dots , n\}$ and
that $\Delta (A)$ is a pre-Dieudonne determinant in the sense of
\cite {Dr}. If $A$ is Hermitian, then $\Delta (A)$ is a product of
real numbers and, therefore, $\Delta (A)$ is real.

\proclaim{Proposition 3.3.1} Let $p(I)$ be the parity of the
ordering $I$. Then $\Delta (A)=p(I)p(J)D_{I,J}(A)$.
\endproclaim
The proof follows from homological relations for
quasideterminants.

\proclaim {Theorem 3.3.2} Let $A$ be a Hermitian quaternionic
matrix. Then $\Delta (A)=M(A)$ {\rm(}see {\rm (3.3.1))}.
\endproclaim
\demo {Proof} We use the noncommutative Sylvester formula for quasideterminants
(Theorem 1.5.2).

For $i,j=2,\dots , n$ define a Hermitian matrix  $B_{ij}$ by the formula
$$
B_{ij}=\pmatrix
a_{11} & a_{1j}\\
a_{i1} & a_{ij}
\endpmatrix.
$$
Let $b_{ij}=M(B_{ij})$  and
$c_{ij}=|B_{ij}|_{11}$.

Note that $B=(b_{ij})$ and $C=(c_{ij})$  also are Hermitian
matrices. It follows from (3.3.1) that  $M(A)=a_{nn}^{2-n}M(B)$.
Note, that $M(B)=a_{nn}^{n-1}M(C)$, therefore,  $M(A)=a_{nn}M(C)$.

By Theorem 1.5.2 ,
$|A|_{11}=|C|_{11}$, $|A^{11}|_{22}=|C^{11}|_{22}$, $\dots $. So,
$$
\multline
|A^{11}|_{22}|A^{11}|_{22}\dots |A^{12\dots n-1, 12\dots
n-1}|_{n-1, n-1}\\
= |C^{11}|_{22}|C^{11}|_{22}\dots |C^{12\dots
n-1, 12\dots n-1}|_{n-1, n-1}.
\endmultline
$$
The product on the left-hand side equals  $\Delta (A)a_{nn}^{-1}$ and the
product on right-hand side equals
 $\Delta (C)$, so  $\Delta (A)= \Delta (C)a_{nn}=M(A)$. \qed
\enddemo

\subhead 3.4. Moore determinants and norms of quaternionic
matrices \endsubhead
\proclaim{Proposition 3.4.1} For generic matrices  $A, B$ we have
$$
\nu (A)=\Delta (A)\Delta (A^*)=\Delta (AA^*).
$$
\endproclaim

Since $AA^*$ is a Hermitian matrix, one has the following

\proclaim{Corollary 3.4.2} $\nu (A)=M(AA^*)$.
\endproclaim

\subhead 3.5. Study determinants \endsubhead
An embedding of the field of complex numbers $\Bbb C$ into $\Bbb H$
is defined by an image of $\bold i\in \Bbb C$. Chose the embedding given by
$x+y\bold i \mapsto x+y\bold i + 0\bold j + 0\bold k $, where
$x,y\in \Bbb R$ and identify $\Bbb C$ with its image in $\Bbb H$.
Then any quaternion $a$ can be uniquely written as $a=\alpha + \bold j\beta $
where $\alpha, \beta \in \Bbb C$.

Let $M(n, F)$ be the algebra of matrices of order $n$ over a field $F$.
Define a homomorphism $\theta: \Bbb H\rightarrow M(2, \Bbb C)$ by setting
$$
\theta (a)=\left (\matrix
\alpha & -\bar \beta \\
\beta & \bar \alpha \endmatrix \right ).
$$
For $A=(a_{ij})\in M(n, \Bbb H)$, set $\theta _n (A)=(\theta
(a_{ij}))$.
This extends $\theta $ to homomorhism of matrix algebras
$$\theta _n: M(n, \Bbb H) \rightarrow M(2n, \Bbb C).$$ 

In 1920, Study \cite{S} defined a determinant $S(A)$ of  a
quaternionic matrix $A$ of order $n$ by setting $S(A)=\det \theta
_n (A)$. Here $\det $ is the standard determinant of a complex
matrix. The following proposition is well known (see \cite {As}).

\proclaim{Proposition 3.5.1} For any quaternionic matrix $A$
$$S(A)=M(AA^*).$$
\endproclaim

The proof in \cite {As} was based on properties of eigenvalues of
quaternionic matrices. Our proof based on Sylvester's identity
and homological relations actually shows that $S(A)=\nu (A)$ for
a generic matrix $A$.

\subhead 3.6. Quantum determinants \endsubhead

Note, first of all, that quantum determinants and the Capelli 
determinants (to be discussed in Section 3.7) are not 
defined for all matrices over the corresponding algebras. 
For this reason, they are not actual determinants, but, rather, 
``determinant-like'' expressions. However, using the 
traditional terminology, we will talk about quantum and 
Capelli determinants.

We say that $A=(a_{ij}), i,j=1,\dots, n$, is a {\it quantum matrix\/} if, 
for some central invertible element $q\in F$, the elements $a_{ij}$ satisfy the
following commutation relations:
$$
\aligned
a_{ik}a_{il}&=q^{-1}a_{il}a_{ik}\ \  {\text {for}}\   k<l,\\
a_{ik}a_{jk}&=q^{-1}a_{jk}a_{ik}\ \  {\text {for}} \  i<j,\\
a_{il}a_{jk}&=a_{jk}a_{il}\ \ {\text {for}}\  i<j, k<l,\\
a_{ik}a_{jl}-a_{jl}a_{ik}&=(q^{-1}-q)a_{il}a_{jk} \ \ 
{\text {for}}\   i<j,  k<l.
\endaligned
\tag 3.6.1
$$

Denote by $\Cal A(n,q)$ the algebra with generators $(a_{ij})$,
$i,j=1,\dots , n$, satisfying relations (3.6.1). The center of this algebra is
the one-dimensional subspace generated by the so called {\it quantum
determinant\/} of $A$.

The quantum determinant $\det_q A$ is defined as follows:
$$
\det {}_q A=\sum _{\sigma \in S_n}(-q)^{-l(\sigma )}
a_{1\sigma (1)}a_{2\sigma (2)}\dots a_{n\sigma (n)},
$$
where $l(\sigma )$ is the number of inversions in $\sigma $.

If $A$ is a quantum matrix, then any square submatrix of $A$
also is a quantum matrix with the same $q$.

Note that the algebra $\Cal A(n,q)$ admits the ring of fractions.

\proclaim{Theorem 3.6.1 \rm (\cite {GR,  KL})} In the  ring of fractions
of the algebra $\Cal A(n,q)$ we have
we have $$
\det {}_q A=
(-q)^{i-j}|A|_{ij}\cdot \det {}_q A^{ij}
=(-q)^{i-j}\det {}_q A^{ij}\cdot |A|_{ij}.
$$
\endproclaim

\proclaim{Corollary 3.6.2 \rm (\cite {GR, KL})} In the  ring of fractions
of the algebra $\Cal A(n,q)$ we have
$$
\det {}_q A=|A|_{11}|A^{11}|_{22}\dots a_{nn}
$$
and all factors on the right-hand side commute.
\endproclaim

An important generalization of this result for matrices
satisfying Faddeev--Reshe\-tikhin--Takhtadjan relations is given in \cite {ER}.

\subhead 3.7.  Capelli determinants \endsubhead
Let $X=(x_{ij})$, $i,j=1,\dots ,n$ be a matrix of formal
commuting variables and $X^T$ the transposed matrix. Let
$D=(\partial _{ij})$, $\partial _{ij}=\partial / \partial x_{ij}$,
be the matrix of the corresponding differential operators. Since each of the
matrices $X$, $D$ consists of commuting entries, $\det X$ and $\det D$ make
sense. Let us set $X^TD=(f_{ij})$, so that $f_{ij}=\sum _{k}x_{ki}\partial
/\partial x_{kj}$. 

Let $W$ be a diagonal matrix, $W={\text {diag}}(0,1,2,\dots , n)$.

By definition, the Capelli determinant $\det_{Cap} $ of $X^TD-W$ equals to
the sum
$$
\sum _{\sigma \in S_n}(-1)^{l(\sigma )} f_{\sigma (1)1}
(f_{\sigma (2)2}-\delta _{\sigma (2)2})\dots (f_{\sigma
(n)n}-(n-1)\delta _{\sigma (n)n}).
$$

The classical Capelli identity says that the sum is equal to
$\det X\det D$.

Set $Z=X^TD-I_n$.
It was shown in \cite {GR1, GR2} that the Capelli determinant can
be expressed as a product of quasideterminants. More precisely, 
let $\Cal D$ be the algebra of polynomial differential operators with variables
$x_{ij}$.

\proclaim{Theorem 3.7.1} In the ring of fractions of the algebra $\Cal D$ we 
have
$$ 
|Z|_{11}|Z^{11}|_{22}\dots z_{nn}=\det X \det D
$$
and all factors on the left-hand side commute. 
\endproclaim

It is known \cite {We} that the right-hand side in the theorem is equal to the Capelli
determinant. 

This theorem can also be interpreted in a different
way.

 Let $A=(e_{ij})$, $i,j=1,\dots n$ be the matrix of the
standard generators of the universal enveloping algebra
$U(gl_n)$. Recall that these generators satisfy the relations
$$[e_{ij}, e_{kl}]=\delta _{jk}e_{il}-\delta _{li}e_{kj}.$$

Let $E_n$ be the identity matrix of order $n$.
It is well known (see, for example, \cite {Ho}) that
coefficients of the polynomial in a central variable $t$
$$
\det(I_n+tA):=\sum _{\sigma \in S_n}(-1)^{l(\sigma )}
(\delta_{\sigma (1)1}+t e_{\sigma (1)1})\dots (\delta _{\sigma
(n)n}+t(e_{\sigma (n)n}-(n-1)\delta _{\sigma (n)n}))
$$
generate the center of $U(gl_n)$.

Theorem 3.7.1 can be reformulated in the following way \cite {GKLLRT}.

\proclaim {Theorem 3.7.2} $\det (I_n+tA)$ can be factored in the algebra of
formal power series in $t$ with coefficients in $U(gl_n)$:
$$
\align
\det(I_n+tA)&=(1+te_{11}) \left| \matrix
1+t(e_{11}-1) & te_{12}\\
te_{21}         & \boxed {1+t(e_{22}-1)}
\endmatrix \right |\cdot \dots
\\
&\quad \cdot 
\left| \matrix
1+t(e_{11}-n+1) & \dots & te_{1n}\\
\dots         & \dots & \dots\\
te_{n1}           &\dots  &\boxed{1+t(e_{nn}-n+1)}
\endmatrix \right |
\endalign
$$
and the factors on the right-hand side commute with each other.
\endproclaim

The above version is obtained by using the classical embedding of
$U(gl_n)$ into the Weyl algebra generated by $(x_{ij},
\partial/\partial x_{ij})$, $i,j=1,\dots , n$, where $e_{ij}$
corresponds to
$$f_{ij}=\sum _{k=1}^nx_{ki}\partial/\partial x_{kj}$$.

\subhead 3.8.  Berezinians \endsubhead
Let $p(k)$ be the parity of an integer $k$, i.e. $p(k)=0$ if $k$ is even
and $p(k)=1$ if $k$ is odd. A (commutative) super-ring over $R^0$ is a 
ring $R=R^0\oplus R^1$ such that 

(i) $\ a_ia_j\in R^{p(i+j)}$ for any $a_m\in R^m$, $m=0,1$,  

(ii) $\ ab=ba$ for any $a\in R^0$, $b\in R$, and $cd=-dc$ for
any $c,d\in R^1$.

Let 
$A=\pmatrix X&Y\\Z&T
\endpmatrix$ be an $(m+n)\times (m+n)$-block-matrix over a super-ring 
$R=R^0\oplus R^1$, where
$X$ is an $m\times m$-matrix over $R^0$, $T$ is an $n\times n$-matrix over $R^0$,
and $Y,Z$ are matrices over $R^1$. If  $T$ is an invertible matrix, then
$X-YT^{-1}Z$ is an invertible matrix over commutative ring $R^0$.
Super-determinant, or Berezinian, of $A$ is defined by the following formula:
$$ 
{\text{Ber}}\, A=\det (X-YT^{-1}Z)\det T^{-1}.
$$
Note that ${\text{Ber}}\, A \in R^0$.

\proclaim {Theorem 3.8.1} Let $R^0$ be a field. Set $J_k=\{1,2,\dots ,k\}$,
$k\leq m+n$ and $A^{(k)}=A^{J_k, J_k}$. Then ${\text {\rm Ber}}\, A$ is a product
of elements of $R^0$: 
$$
{\text {\rm Ber}}\, A=|A|_{11}|A^{(1)}|_{22}\dots |A^{(m-1)}|_{mm}
|A^{(m)}|_{m+1,m+1}^{-1}\dots |A^{(m+n-1)}|_{m+n,m+n}^{-1}.
$$
\endproclaim

\subhead 3.9.  Cartier-Foata determinants \endsubhead
Let $A=(a_{ij})$, $i,j=1,\dots , n$ be a matrix such that
the entries $a_{ij}$ and $a_{kl}$ commute when $i\neq k$.
In this case Cartier and Foata \cite {CF, F} defined a determinant of $A$ as
$$ 
\det {}_{CF} A=\sum _{\sigma \in S_n}
(-1)^{l(\sigma )}a_{1\sigma (1)}a_{2\sigma (2)}\dots a_{n\sigma(n)}.
$$

The order of factors in monomials $a_{1\sigma (1)}a_{2\sigma (2)}\dots
a_{n\sigma(n)}$ is insignificant.

Let $C_n$ be the algebra over a field $F$ generated by  $(a_{ij})$,
$i,j=1,\dots,n$, with  relations $a_{ij}a_{kl}=a_{kl}a_{ij}$ if $i\neq k$.
Algebra $C_n$ admits the ring of fractions.

\proclaim {Theorem 3.9.1} In the ring of fractions of algebra $C_n$, let
$A=(a_{ij})$, $i,j=1,\dots , n$ be a matrix such that the entries $a_{ij}$ and
$a_{kl}$ commute when $i\neq k$.
$$
|A|_{pq}=(-1)^{p+q}\det {}_{CF}\ (A^{pq})^{-1}\det {}_{CF} A
\tag 3.9.1
$$
and all factors in (3.9.1) commute.
\endproclaim

\proclaim {Corollary 3.9.2} In the ring of fractions of algebra $C_n$ we have
$$
det _{CF}=|A|_{11}|A^{11}|_{22}\dots a_{nn}
$$
and all factors commute.
\endproclaim

\head 4. Noncommutative Pl\"ucker and Flag Coordinates\endhead
Most of the results described in this section were obtained in \cite {GR4}.

\subhead 4.1. Commutative Pl\"ucker coordinates \endsubhead
Let $k\leq n$ and $A$ be a $k\times n$-matrix over a commutative ring $R$.
Denote by $A(i_1,\dots , i_k)$ the $k\times k$-submatrix of $A$ consisting of columns
labeled by the indices $i_1,\dots , i_k$. Define 
$p_{i_1\dots i_k}(A):=\det A(i_1,\dots , i_k)$. The elements 
$p_{i_1\dots i_k}(A)\in R$ are called Pl\"ucker coordinates of the matrix $A$.
The Pl\"ucker coordinates $p_{i_1\dots  i_k}(A)$
satisfy the following properties:

(i) (invariance) $p_{i_1\dots i_k}(XA)=\det X \cdot p_{i_1\dots  i_k}(A)$ for
any $k\times k$-matrix $X$ over $R$;

(ii) 	(skew-symmetry) $p_{i_1\dots i_k}(A)$ are skew-symmetric in indices
$i_1,\dots,i_k$; 
in particular, $p_{i_1\dots i_k}(A)=0$ if a pair of indices coincides;

(iii) (Pl\"ucker relations) Let $i_1,\dots , i_{k-1}$ be $k-1$ distinct numbers
which are chosen from the set $1,\dots , n$, and $j_1,\dots , j_{k+1}$ be $k+1$
distinct numbers chosen from the same set. Then
$$
\sum _{t=1}^k(-1)^t p_{i_1\dots i_{k-1}j_t}(A) 
p_{j_1\dots j_{t-1}j_{t+1}\dots j_{k+1}}(A)=0.
$$

\example{Example} For $k=2$ and $n=4$ the Pl\"ucker relations in (iii)
imply the famous identity
$$
p_{12}(A)p_{34}(A)-p_{13}(A)p_{24}(A)+p_{23}(A)p_{14}(A)=0. \tag 4.1.1
$$
\endexample

Historically, Pl\"ucker coordinates were introduced as coordinates on Grassmann
manifolds. Namely, let $R=F$ be  a field and $G_{k,n}$ the Grassmannian of
$k$-dimensio\-nal subspaces in the $n$-dimensio\-nal vector  space $F^n$. To each
$k\times n$-matrix $A$ of rank $k$ we associate the subspace of $F^n$ 
generated by the rows of $A$. By the invariance property (i), we can view each 
Pl\"ucker coordinate $p_{i_1\dots   i_k}$ as a section of a certain ample line
bundle on $G_{k,n}$, and all these sections together  define an embedding of
$G_{k,n}$ into the projective space $\Bbb P^N$ of dimension $N=\binom  kn -1$. 
In this sense, Pl\"ucker coordinates are projective coordinates on $G_{k,n}$.

\subhead 4.2. Quasi-Pl\"ucker coordinates for $n\times (n+1)$- and
$(n+1)\times n$-matrices \endsubhead
Let $A=(a_{ij})$, $i=1,\dots , n$, $j=0,1,\dots , n$, be a matrix over a
division ring $R$.  Denote by $A^{(k)}$ the $n\times n$-submatrix 
of $A$ obtained from
$A$ by removing the $k$-th column and suppose that all $A^{(k)}$ are invertible.
Choose an arbitrary $s\in \{1,\dots , n\}$, and denote
$$
q^{(s)}_{ij}(A)=|A^{(j)}|_{si}^{-1}|A^{(i)}|_{sj}. 
$$
\proclaim{Proposition 4.2.1} The element $q^{(s)}_{ij}(A)\in R$ does not depend on
$s$. 
\endproclaim
We denote the common value of $ q^{(s)}_{ij}(A)$ by $ q_{ij}(A)$
and call $ q_{ij}(A)$ the {\it left quasi-Pl\"ucker coordinates\/} of the matrix $A$.

\demo{Proof of Proposition 4.2.1} 
Considering the columns of the matrix $A$ as $n+1$ vectors
in the right $n$-dimensional space $R^n$ over $R$, we see 
that there exists a nonzero $(n+1)$-vector $(x_1,\dots,x_{n+1})\in R^{n+1}$
such that
$$
A\left (\matrix x_0\\ \dots \\ x_n \endmatrix \right ) =0.
$$ 
This means that
$$
A^{(j)}\left (\matrix x_0\\ \dots \\ \widehat{x}_j \\ \dots \\ x_n \endmatrix \right ) =
-\left (\matrix a_{1j}\\ \ \\ \dots  \\  \\a_{nj} \endmatrix \right )x_j.
$$ 
Since all submatrices $A^{(k)}$ are invertible, each $x_I$ is a nonzero element
of $R$. Cramer's rule and transformations properties for quasideterminanats imply that
$|A^{(j)}|_{si}x_i=-|A^{(i)}|_{sj}x_j$. Therefore, 
$$
q_{ij}^{(s)}(A)=|A^{(j)}|_{si}^{-1}|A^{(i)}|_{sj}= -x_ix_j^{-1}\tag 4.2.1 
$$
does not depend on $s$.\qed
\enddemo
\proclaim{Proposition 4.2.2} If $g$ is an invertible $n\times n$-matrix over $R$, then
$q_{ij}(gA)=q_{ij}(A)$. 
\endproclaim
\demo{Proof} We have 
$gA\left (\matrix x_0\\ \dots \\ x_n \endmatrix \right ) =0$.
Therefore, $q_{ij}(gA)=-x_ix_j^{-1}=q_{ij}(A)$. \qed\enddemo

In the commutative case, $q_{ij}(A)$ is a ratio of
two Pl\"ucker coordinates: $q_{ij}(A)= p_{1,\dots,\widehat
j,\dots,n}/p_{1,\dots,\widehat i,\dots,n} = {\det A^{(j)}}/{\det A^{(i)}}$.

Similarly, we define the right quasi-Pl\"ucker coordinates $r_{ij}(B)$ for    
$(n+1)\times n$-matrix $B=(b_{ji})$. 
Denote by $B^{(k)}$ the submatrix of $B$ obtained from $B$ by removing the
$k$-th row. Suppose that all $B^{(k)}$ are invertible, choose $s\in \{1,\dots ,
n\}$,  and set $r^{(s)}_{ij}(B)=|B^{(j)}|_{is}|B^{(i)}|_{js}^{-1}$. 

\proclaim{Proposition 4.2.3} {\rm (i)} The element $r^{(s)}_{ij}(B)$ does not depend of $s$.

Denote the common value of elements $r^{(s)}_{ij}(B)$ by $r_{ij}(B)$.

{\rm(ii)} If $g$ is an invertible $n\times n$-matrix over $R$, then
$r_{ij}(Bg)=r_{ij}(B)$. 
\endproclaim

In the commutative case, $r_{ij}(A)= {\det B^{(j)}}/{\det B^{(i)}}$.

\subhead 4.3.  Definition of left quasi-Pl\"ucker coordinates. General case\endsubhead Let $A=(a_{pq})$, 
$p=1, 
\dots, k$,  $q=1, \dots, n$, $ k< n$, be a matrix over
a division ring $R$.
Choose $1\leq i, j, i_1,\dots, i_{k-1}\leq n$ such
that $i\notin I=\{ i_1,\dots, i_{k-1}\}$.  
Let $A(i,j,i_1,\dots , i_{k-1})$ be
the $k\times (k+1)$-submatrix of $A$ with columns labeled by $i,j,i_1,\dots , i_{k-1}$.

\definition{Definition 4.3.1} Define left
quasi-Pl\"ucker coordinates $ q^{I}_{ij}(A)$ of the matrix $A$ by the formula
$$
q^I_{ij}(A)=q_{ij}(A(i,j,i_1,\dots , i_{k-1})).
$$
\enddefinition
By Proposition 4.2.1, left quasi-Pl\"ucker coordinates are given by the formula
$$
q^{I}_{ij}(A) =\left|\matrix a_{1i}a_{1i_1}&\dots & a_{1,i_{k-1}}\\
{}  &\dots & {}\\
a_{ki}a_{ki_1} & \dots &a_{ki_{k-1}}\endmatrix\right|^{-1}_{si} \cdot
\left|\matrix a_{1j} a_{1,i_1} &\dots &a_{1,i_{k-1}}\\
 {}  &\dots  &{}\\
a_{kj} a_{ki_1} &\dots &a_{ki_{k-1}}\endmatrix\right|_{sj}
$$
for an arbitrary $s$, $1\le s\le k$.

\proclaim{Proposition 4.3.2} If $g$ is an invertible
$k\times k$-matrix over $R$, then $q^I_{ij}(gA)=q^I_{ij}(A)$. 
\endproclaim

\demo{Proof} Use Proposition 4.2.2.
\qed\enddemo

In the commutative case $q^I_{ij} = {p_{jI}/ p_{iI}}$, where
$p_{\alpha_1\dots \alpha_k}$ are the standard Pl\"ucker
coordinates.

\subhead 4.4. Identities for the left quasi-Pl\"ucker coordinates \endsubhead  
The following
properties of $q^I_{ij}$ immediately follow from the definition.
\roster
\item"{(i)}" $q^I_{ij}$ does not depend on the ordering on elements in $I$;
\item"{(ii)}" $q^I_{ij} = 0$ for $j\in I$;
\item"{(iii)}" $q^I_{ii} = 1$ and  $q^I_{ij} \cdot q^I_{jk} = q^I_{ik}$.
\endroster

\proclaim{Theorem 4.4.1 \rm (Skew-symmetry)}  Let $N$, $|N|=k+1$, be a set
of indices, $i,j,m\in N$.  Then
$$
q^{N\smallsetminus\{i,j\}}_{ij} \cdot q^{N\smallsetminus\{j,m\}}_{jm}\cdot
q^{N\setminus\{m,i\}}_{mi} =- 1.
$$
\endproclaim

\proclaim{Theorem 4.4.2 \rm (Pl\"ucker relations)}
Fix $M=(m_1,\dots, m_{k-1})$, $L=(\ell_1,\dots, \ell_k)$. Let $i\notin M$. Then
$$
\sum_{j\in L} q^M_{ij} \cdot q_{ji}^{L\setminus\{j\}}=1.
$$
\endproclaim

\example{Examples} Suppose that $k=2$. 

1) From Theorem 4.4.1 it follows that
$$
q^{\{\ell\}}_{ij}\cdot q^{\{i\}}_{j\ell}\cdot q^{\{j\}}_{\ell i}=-1.
$$
In the commutative case, $q^{\{\ell\}}_{ij} = {p_{j\ell}\over
p_{i\ell}}$ so this identity follows from the skew-symmetry $p_{ij} =- p_{ji}$.

2)  From Theorem 4.4.2 it follows that for any $i,j,\ell, m$
$$
q^{\{\ell\}}_{ij} \cdot q_{ji}^{\{m\}} + q_{im}^{\{\ell\}}\cdot
q^{\{j\}}_{mi} =1.
$$
In the commutative case this identity implies the standard identity
(cf. (4.1.1))
$$
p_{ij}\cdot p_{\ell m} - p_{i\ell}\cdot p_{jm} + p_{im}\cdot p_{\ell j} =0.
$$
\endexample

\remark{Remark} The products $p^{\{\ell \}}_{ij}p^{\{m \}}_{ji}$
(which in the commutative case are equal to
${p_{j\ell}\over p_{i\ell}}\cdot {p_{im}\over p_{jm}}$)
can be viewed as noncommutative cross-ratios.
\endremark

To prove Theorems 4.4.1 and 4.4.2 we need the following lemma.
Let $A=(a_{ij})$, $i=1,\dots ,k$, $j=1,\dots , n$, $k< n$, be a matrix over
a division ring. Denote by $A_{j_1,\dots ,j_{\ell }}$, $\ell \leq n$,
the $k\times \ell$-submatrix $(a_{ij})$,
$i=1,\dots ,k$, $j=j_1,\dots ,j_{\ell}$. Consider the $n\times n$-matrix
$$
X=\pmatrix A_{1\dots k}&A_{k+1\dots n}\\0&E_{n-k}\endpmatrix ,
$$
where $E_m$ is the identity matrix of order $m$.

\proclaim{Lemma 4.4.3} Let $j< k< i$. If
$q^{1\dots \hat {j}\dots k}_{ij}(A)$ is defined, then $|X|_{ij}$ is
defined and
$$
|X|_{ij}=-q^{1\dots \hat {j}\dots k}_{ij}(A). \tag 4.4.1 
$$
\endproclaim

\demo{Proof} We must prove that
$$
|X|_{ij}=-|A_{1\dots \hat {j}\dots k i}|^{-1}_{si}
\cdot |A_{1\dots k}|_{sj} \tag 4.4.2 
$$
provided the right-hand side is defined. We will prove this by induction on
$\ell =n-k$. Let us assume that formula (2.2) holds for $l=m$ and prove it
for $\ell =m+1$. Without loss of generality we can take $j=1, i=k+1$.
By homological relations (Theorem 1.4.3)
$$
|X|_{k+1,1}=-|X^{k+1,1}|^{-1}_{s,k+1}\cdot |X^{k+1,k+1}|_{s1}$$
for an appropriate $1\leq s\leq k$. Here
$$
\align
X^{k+1,1}&=
\pmatrix A_{2\dots k+1}&A_{k+2\dots n}\\0&E_{n-k-1}\endpmatrix ,\\
X^{k+1,k+1}&=
\pmatrix A_{1\dots k}&A_{k+2\dots n}\\0&E_{n-k-1}\endpmatrix .
\endalign
$$
By the induction assumption
$$
\align
|X^{k+1,1}|_{s,k+1}&=
-|A_{23\dots k k+2}|^{-1}_{s,k+2}\cdot |A_{23\dots k+1}|_{s,k+1}, \\
|X^{k+1,k+1}|_{s1}&=
-|A_{23\dots k k+2}|^{-1}_{s,k+2}\cdot |A_{1\dots k}|_{s1} 
\endalign
$$
and $|X|_{k+1,1}=-p^{23\dots k}_{k+1,1}$. \qed\enddemo

To prove Theorem 4.4.2 we apply the second formula in Theorem 1.4.4 
to the matrix
$$
X=\pmatrix A_{1\dots k}&A_{k+1\dots n}\\0&E_{n-k}\endpmatrix 
$$
for $M=(k+1,\dots ,n)$ and any $L$ such that $|L|=n-k-1$. 
By Lemma 4.4.3,
$|X|_{m_i\ell }=-q^{1\dots \hat \ell \dots k}(A)$,
$|X^{M_i,L}|_{m_i q}=-p^{{1\dots n}\setminus L}_{m_i q}(A)$, and
Theorem 4.4.2 follows from Theorem 1.4.4. \qed
\medskip

To prove Theorem 4.4.1 it is sufficient to take the matrix $X$ for $n=k+1$ and
use homological relations.
\qed

\proclaim{Theorem 4.4.4}  Let $A =(a_{ij})$, $i=1, \dots, k$, $j=1,\dots, n$,
be a matrix with formal entries and $f (a_{ij})$  an element of
a free skew-field $F$ generated by $a_{ij}$.  Let $f$ be invariant
under the transformations
$$
A\to gA
$$
for all invertible $k\times k$-matrices $g$ over $F$.  Then $f$
is a rational function of the quasi-Pl\"ucker coordinates.
\endproclaim

\demo{Proof} Let $b_{ij}=a_{ij}$ for $i,j=1,\dots, k$. Consider the matrix
$B=(b_{ij})$. Then $B^{-1}=(|B|^{-1}_{ji})$.
Set $C=(c_{ij})=B^{-1}A$. Then
$$
c_{ij}=\cases \delta _{ij}&\quad \text{for $j\leq k$}, \\
q^{1\dots \hat {i}\dots k}_{ij}(A) & \quad\text{for $j>k$}.
\endcases $$
By invariance, $f$ is a rational expression of $c_{ij}$ with $j>k$.
\enddemo

\subhead 4.5. Right quasi-Pl\"ucker coordinates\endsubhead 
Consider a matrix $B=(b_{pq})$, $p=1,\dots, n$;
$q=1,\dots, k$, $k< n$ over a division ring $F$. Choose $1\leq i, j$,
$i_1,\dots, i_{k-1} \leq n$ such that $j\notin I = (i_1,\dots,
i_{k-1})$. Let $B(i,j,i_1,\dots,i_{k-1})$ be the $(k+1)\times k$-submatrix 
of $B$ with rows labeled by $i,j,i_1,\dots,i_{k-1}$.

\definition{Definition 4.5.1} Define right
quasi-Pl\"ucker coordinates $ r^{I}_{ij}(B)$ of the matrix $B$ by the formula
$$
r^I_{ij}(B)=r_{ij}(B(i,j,i_1,\dots , i_{k-1})).
$$
\enddefinition
By Proposition 4.2.3, right quasi-Pl\"ucker coordinates are given by the formula
$$
r^I_{ij}(B) = \vmatrix b_{i1} &\dots & b_{ik}\\
b_{i_11} & \dots & b_{i_1k}\\
{} &\dots &{}\\
b_{i_{k-1}1} &\dots &b_{i_{k-1}k}\endvmatrix_{it} \cdot
\vmatrix b_{j1} & \dots &b_{jk}\\
b_{i_11} &\dots &b_{i_1 k}\\
{} &\dots &{} \\
b_{i_{k-1} 1} &\dots &b_{i_{k-1}k}\endvmatrix^{-1}_{jt} 
$$
for an arbitrary $t$, $1\leq t\leq k$.

\proclaim{Proposition 4.5.2}
$r^I_{ij}(Bg) = r^I_{ij}(B)$ for each invertible
$k\times k$-matrix $g$ over $F$.
\endproclaim

\subhead 4.6.  Identities for the right quasi-Pl\"ucker coordinates\endsubhead  
Identities for $r^I_{ij}$ are dual to correspoding
identities for the left quasi-Pl\"ucker coordinates $q^I_{ij}$. Namely,

\roster
\item"{(i)}" $r^I_{ij}$ does not depend on the ordering on elements of $I$;
\item"{(ii)}" $r^I_{ij} = 0$ for $i\in I;$
\item"{(iii)}"  $r^I_{ii} = 1$ and $r^I_{ij} \cdot r^I_{jk} = r^I_{ik}$.
\endroster

\proclaim{Theorem 4.6.1 \rm (Skew-symmetry)}  Let $N$, $|N|=k+1$,  be a set
of indices, $i,j,m\in N$.  Then
$$
r^{N\smallsetminus\{i,j\}}_{ij} \cdot r^{N\smallsetminus\{j,m\}}_{jm}\cdot
r^{N\setminus\{m,i\}}_{mi} =- 1.
$$
\endproclaim

\proclaim{Theorem 4.6.2 \rm (Pl\"ucker relations)}
Fix
$M=(m_1,\dots, m_{k-1})$, $L=(\ell_1,\dots, \ell_k)$. Let $i\notin M$. Then
$$
\sum_{j\in L} r_{ij}^{L\setminus\{j\}}r^M_{ij} =1.
$$
\endproclaim

\subhead 4.7. Duality between quasi-Pl\"ucker
coordinates\endsubhead
Let $A=(a_{ij})$, $i=1,\dots, k$, $j=1,\dots, n$; and $B=(b_{rs})$, $r=1,\dots,
n$, $s=1,\dots, n-k$.  Suppose that $AB=0$.  (This is equivalent to the
statement that the subspace generated by the rows of $A$ in the left
linear space $F^n$ is dual to the subspace generated by the columns of $B$
in the dual right linear space.)  Choose indices $1\leq i, j\leq n$ and a subset
$I\subset [1,n]$, $|I|=k-1$, such that
$i\notin I$. Set $J=([1,n]\setminus I)\setminus \{i,j\}$.

\proclaim{Theorem 4.7.1} We have
$$
q^I_{ij}(A) + r^J_{ij} (B) = 0.
$$
\endproclaim

\subhead 4.8. Quasi-Pl\"ucker coordinates for $k\times n$-matrices for
different $k$\endsubhead
Let  $A=(a_{\alpha\beta})$, $\alpha=1,\dots, k$, $\beta = 1, \dots, n$, be a
$k\times n$-matrix over a noncommutative division ring $R$ and $A'$  a
$(k-1)\times n$-submatrix of $A$. Choose $1\leq i, j, m, j_1,\dots,
j_{k-2} \leq n$  such that $i\neq m$ and $i,m\notin J=\{j_1, \dots, j_{k-2}\}$.

\proclaim{Proposition 4.8.1} We have
$$
q^J_{ij}(A')= q_{ij}^{J\cup\{m\}}(A) + q^J_{im}(A')\cdot
q_{mj}^{J\cup\{i\}}(A).
$$
\endproclaim

\subhead 4.9. Applications of quasi-Pl\"ucker coordinates
\endsubhead

\subhead Row and column expansion of a
quasideterminant\endsubhead
Some of the results obtained in [GR], [GR1], [GR2] and partially
presented in Section I can be rewritten in terms of quasi-Pl\"ucker
coordinates.

Let $A=(a_{ij}), i, j=1, \dots, n$,  be a matrix over a division ring $R$.
Choose $1\leq \alpha,\beta\leq n $. Using the notation of section I
let $B=A^{\{\alpha \}, \emptyset},
C=A^{\emptyset , \{\beta \}}$ be the $(n-1)\times n$ and $n\times (n-1)$
submatrices of $A$ obtained by deleting the $\alpha $-th row and
the $\beta $-th column respectively.
For $j\neq \beta$ and $i\neq\alpha$ set
$$
\align
q_{j\beta} &= q^{1\dots\hat j \dots \hat \beta\dots n}_{j\beta}(B),\\
r_{\alpha i}&= r^{1\dots\hat\alpha\dots\hat i \dots n}_{\alpha i}(C).
\endalign
$$

\proclaim{Proposition 4.9.1}
{\rm(i)}  $|A|_{\alpha \beta} = a_{\alpha\beta} - \sum_{j\neq \beta}
a_{\alpha j} q_{j\beta}$,

{\rm(ii)} $|A|_{\alpha\beta} = a_{\alpha\beta} - \sum_{i\neq \alpha}
r_{\alpha i} a_{i\beta}$

\noindent provided the terms in the right-hand side of these formulas are defined.
\endproclaim

\subhead Homological relations\endsubhead
\proclaim{Proposition 4.9.2}  In the previous notation, 

{\rm(i)}  $|A|^{-1}_{ij} \cdot |A|_{i\ell} =- q_{j\ell}$  \quad (row
relations)

{\rm(ii)}  $|A|_{ij} \cdot |A|^{-1}_{kj} =- r_{ik}$      \quad (column relations).
\endproclaim

\proclaim{Corollary 4.9.3}  In the previous notation, let $(i_1,\dots, i_s),
(j_1,\dots, j_t)$ be sequences of indices such that $i\neq i_1$, $i_1\neq
i_2$,\dots, $i_{s-1}\neq i_s$; $j \neq j_1$, $j_1\neq j_2$,\dots, $j_{t-1}\neq j_t$.
Then
$$
|A|_{i_sj_t} = q_{i_si_{s-1}}\dots q_{i_2i_1} q_{i_1i}\cdot |A|_{ij}
\cdot r_{jj_1}r_{j_1j_2}\dots r_{j_{t-1}j_t}.
$$
\endproclaim

\example{Example} For a matrix $A=\pmatrix a_{11}&
a_{12}\\ a_{21} & a_{22}\endpmatrix$ we have
$$
\align
|A|_{22}&=a_{21} \cdot a^{-1}_{11}\cdot |A|_{11} \cdot a^{-1}_{22}\cdot a_{22},
\\
|A|_{11}&=a_{12} \cdot a_{22}^{-1} \cdot a_{21} \cdot a^{-1}_{11}  \cdot |A|_{11} \cdot a^{-1}_{21}  
\cdot a_{22} \cdot a^{-1}_{12} \cdot a_{11}. \endalign
$$
\endexample
\subhead Matrix multiplication\endsubhead The following formula
was already used in the proof of Theorem 4.4.4.
Let $A= (a_{ij})$, $i=1, \dots, n$, $j=1,\dots, m$, $n< m$, $B=(a_{ij})$, $i =
1,\dots, n$, $j=1, \dots, n$, $C=(a_{ik})$, $i=1,\dots, n$, $k=n+1, \dots,m$.

\proclaim{Proposition 4.9.4} Let the matrix $B$ be invertible. Then
$q^{1\dots \hat i\dots n}_{ik}(A)$ are defined for $i=1,\dots n,
k=n+1,\dots m$, and
$$
B^{-1}C=(q^{1\dots\hat i\dots n}_{ik}(A)), \qquad i=1,\dots, n, \quad k= n+1,
\dots, m.
$$
\endproclaim

\subhead Quasideterminant of the product\endsubhead
Let $A=(a_{ij})$, $B=(b_{ij})$,
$i,j=1,\dots n$ be matrices over a division ring $R$. Choose $1\leq k\leq 
n$.
Consider the $(n-1)\times n$-matrix $A'=(a_{ij})$, $i\neq k $,
and the $n\times (n-1)$-matrix $B''=(b_{ij})$,  $j\neq k$.

\proclaim{Proposition 4.9.5} We have
$$
|B|_{kk}\cdot |AB|^{-1}_{kk} \cdot |A|_{kk} =
1+\sum_{\alpha \neq k} r_{k\alpha}\cdot
q_{\alpha k},
$$
where
$r_{k\alpha} = r^{1\dots\hat \alpha\dots n}_{k\alpha}(B'')$
are right quasi-Pl\"ucker coordinates and
$q_{\alpha k} = q_{\alpha k}^{1\dots\hat\alpha\dots n}(A')$
are left quasi-Pl\"ucker coordinates,
provided all expressions are defined.
\endproclaim

The proof follows from the multiplicative property of quasideterminants
and Proposition 4.9.2.

\subhead Gauss decomposition \endsubhead Consider a matrix
$A=(a_{ij})$, $i,j=1,\dots,n$, over a division ring $R$. Let 
$A_k=(a_{ij})$,
$i,j =k,\dots n$, $B_k=(a_{ij})$, $i=1,\dots n$, $j=k, \dots n$, and
$C_k=(a_{ij})$, $i=k,\dots n$, $j=1,\dots n$. These are submatrices of sizes
$(n-k+1)\times (n-k+1)$, $n\times (n-k+1)$, and $(n-k+1)\times n$
respectively.
Suppose that the quasideterminants
$$
y_k=|A_k|_{kk},\qquad  k=1,\dots, n,
$$
are defined and invertible in $R$. 

\proclaim{Theorem 4.9.6 \rm(see \cite{GR1, GR2})} 
$$
A=\pmatrix 1 &{}&x_{\alpha \beta}\\
{}&\ddots&{}\\
0&{}&1\endpmatrix
\pmatrix y_1&{}&0\\
{}&\ddots&{}\\
0&{}&y_n\endpmatrix\pmatrix1&{}&0\\
{}&\ddots&{}\\
z_{\beta \alpha}&{}&1\endpmatrix,
$$
where
$$
\align
x_{\alpha \beta} &=r^{\beta+1\dots n}_{\alpha\beta}(B_{\beta}),
\quad 1\leq \alpha <\beta\leq n,\\
z_{\beta\alpha} &=q^{\beta+1\dots n}_{\beta\alpha} ( C_{\beta}),
\quad 1\leq \alpha <\beta\leq n.\\
\endalign
$$
\endproclaim

Similarly, let $A^{(k)}=(a_{ij})$, $i,j=1,\dots , k$,
$B^{(k)}=(a_{ij})$, $i=1,\dots , n$, $j=1,\dots , k$,
$C^{(k)}=(a_{ij})$, $i=1,\dots , k$, $j=1,\dots , n$.
Suppose that the quasideterminants
$$
y'_k=|A^{(k)}|_{kk}, \qquad k=1,\dots , n,
$$
are defined and invertible in $R$.

\proclaim {Theorem 4.9.7} We have
$$
A=\left (\matrix 1& &0\\
&\ddots & \\
x'_{\beta \alpha}& &1\endmatrix \right )
\left (\matrix y_1'& &0\\
&\ddots & \\
0& &y_n'\endmatrix \right )
\left (\matrix 1& &z_{\alpha \beta }'\\
&\ddots & \\
0& &1\endmatrix \right ),
$$
where
$$
\align
x_{\beta \alpha}'&=r^{1\dots \alpha -1}_{\beta \alpha}
(B^{(\alpha )}), \quad 1\leq \alpha <\beta \leq n,\\
z_{\alpha \beta}'&=q^{1\dots \alpha -1}_{\alpha \beta}
(C^{(\alpha )}), \quad 1\leq \alpha <\beta \leq n.
\endalign
$$
\endproclaim

\subhead Bruhat decompositions \endsubhead
A generalization of Theorem 4.9.6 is given by the following noncommutative analog of the 
Bruhat decomposition. 

\definition{Definition} A square matrix $P$ with entries $0$ and $1$ is called a permutation
matrix if in each row of $P$ and in each column of $P$ there is exactly one entry $1$.
\enddefinition

\proclaim {Theorem 4.9.8 \rm(Bruhat decomposition)}
For an invertible matrix $A$ over a division ring
there exist an upper-unipotent matrix $X$, a low-unipotent matrix $Y$, a diagonal
matrix $D$ and a permutation matrix $P$ such that
$$
A=XPDY.
$$
Under the additional condition that $P^{-1}XP$ is an upper-unipotent matrix, the
matrices $X,P, D, Y$ are uniquely determined by $A$.
\endproclaim

Note that one can always find a decomposition $A=XPDY$ that satisfies the
additional condition.

The entries of matrices $X$ and $Y$ can be written in terms of
quasi-Pl\"ucker coordinates of submatrices of $A$. The entries of
$D$ can be expressed as quasiminors of $A$.

\example{Examples} Let 
$A=\left (\matrix a_{11}&a_{12}\\
a_{21}&a_{22}\endmatrix \right )$.
If $a_{22}\neq 0$, then
$$
A=\pmatrix 1&a_{12}a_{22}^{-1}\\0&1\endpmatrix
\pmatrix |A|_{11}&0\\0&a_{22}\endpmatrix 
\pmatrix 1&0\\a_{22}^{-1}a_{21}&1\endpmatrix.
$$
If $a_{22}=0$ and the matrix $A$ is invertible, then $a_{12}\neq 0$. In this case,
$$
\pmatrix a_{11}&a_{12}\\a_{21}&0\endpmatrix=
\pmatrix 0&1\\1&0\endpmatrix
\pmatrix a _{21}&0\\0&a_{12}\endpmatrix
\pmatrix 1&0\\a_{12}^{-1}a_{11}&1\endpmatrix.
$$
\endexample

An important example of quasi-Pl\"ucker coordinates for the Vandermonde
matrix will be considered later.

\subhead 4.10. Flag coordinates \endsubhead  
Noncommutative flag coordinates were introduced in \cite {GR1, GR2}.

Let $A=(a_{ij})$, $i=1,\dots, k$, $j=1,\dots,n$, be a matrix over a
division ring $R$.  Let $F_p$ be the subspace of the left vector
space $R^n$ generated by the first $p$ rows of $A$. Then $\Cal
F=(F_1\subset F_2\subset\dots\subset F_k)$ is a flag in $R^n$. Put
$$
f_{j_1 \dots j_k} (\Cal F) =\vmatrix a_{1j_1}& \dots  &a_{1j_k}\\
{}&\dots &{}\\
a_{kj_1}&\dots &a_{kj_k}\endvmatrix_{kj_1}.
$$
In \cite{GR1, GR2} the functions $f_{j_1 \dots j_k} (\Cal F)$ were called 
the {\it flag coordinates} of $\Cal F$. 
Transformations properties of quasideterminants imply that
$f_{j_1 \dots j_k} (\Cal F)$
does not depend on the order of the indices $j_2,\dots , j_k$.

\proclaim{Proposition 4.10.1 \rm (see \cite{GR1, GR2})} The
functions $f_{j_1\dots j_m}(\Cal F)$ do not change under left
multiplication of $A$ by an upper unipotent matrix.
\endproclaim

\proclaim{Theorem 4.10.2 \rm (see \cite{GR1, GR2})} The functions 
$f_{j_1 \dots j_k} (\Cal F)$ possess the following
relations:
$$
\align
&f_{j_1j_2j_3 \dots j_k } (\Cal F)f_{j_1j_3 \dots j_k} (\Cal F)^{-1}=
-f_{j_2j_1 \dots j_k} (\Cal F)f_{j_2j_3 \dots j_k} (\Cal F)^{-1},
\\
&f_{j_1 \dots j_k } (\Cal F)f_{j_1\dots j_{k-1}} (\Cal F)^{-1}+
f_{j_2\dots j_kj_1 } (\Cal F)f_{j_2\dots j_k} (\Cal F)^{-1}\\
&\qquad\qquad\qquad+\dots +
f_{j_k j_1\dots j_{k-1}}(\Cal F)f_{j_k j_1\dots j_{k-2}} (\Cal F)^{-1}=0
\endalign
$$
\endproclaim

\example{Example} Let 
$A=\pmatrix a_{11}&a_{12}&a_{13}\\
a_{21}&a_{22}&a_{23}\endpmatrix$.
Then
$f_{12}(\Cal F)a_{11}^{-1}=-f_{21}(\Cal F)a_{12}^{-1}$ and
$f_{12}(\Cal F)a_{11}^{-1}+f_{23}(\Cal F)a_{12}^{-1}+
f_{31}(\Cal F)a_{13}^{-1}=0$.
\endexample

It is easy to see that
$$
q_{ij}^{i_1\dots i_{k-1}}(A)
=\left(f_{ii_1\dots i_{k-1}}(\Cal F)\right)^{-1}\cdot f_{ji_1\dots
i_{k-1}}(\Cal F).
$$
Theorems 4.4.1 and 4.4.2
can be deduced from Theorem 4.10.2.

\head 5. Factorization of Vandermonde quasideterminants and the Vi\`ete Theorem
\endhead

In this section we study factorizations of quasideterminants of
Vandermonde matrices. It is well known that factorizations of Vandermonde
determinants over commutative rings play a fundamental role in mathematics.
Factorizations of non-commutative Vandermonde quasideterminants turn 
out to be equally important. This is why we devote a separate section to these results.
We also use these factorizations to prove the noncommutative Viet\'e theorem, 
which was formulated in \cite {GR3, GR4} using our noncommutative form of
the Sylvester identity.  In Section 6, 7 we will give other applications 
of factorizations of quasideterminants of Vandermonde matrices. 
A good exposition of decompositions of Vandermonde quasideterminants
is given in \cite {Os}.

\subhead 5.1. Vandermonde quasideterminants\endsubhead
Let $x_1, x_2, \dots, x_k$ be a set of elements of a division ring $R$. 
For
$k>1$  the quasideterminant
$$
V(x_1,\dots,x_k)=\vmatrix x_1^{k-1} & \dots &x_k^{k-1}\\
{} &\dots &{}\\
x_1&\dots &x_k\\
1&\dots & 1\endvmatrix_{1k}
$$
is called the {\it Vandermonde\/} quasideterminant.

We say that a sequence of elements $x_1,\dots, x_n\in R$ is {\it independent} if
all quasideterminants $V(x_1,\dots, x_k),\ k=2,\dots , n$, are 
defined and invertible.  For independent sequences $x_1, \dots, x_n$
and $x_1,\dots , x_{n-1}, z$ set 
$$
\align
y_1 &=x_1, \quad z_1 = z\\
y_k &=V(x_1,\dots,x_k) x_k V(x_1,\dots, x_k)^{-1},\quad k\geq 2\\
z_k &=V(x_1,\dots,x_{k-1}, z) z V(x_1,\dots, x_{k-1},z)^{-1},\quad k\geq 2.
\endalign
$$
In the commutative case $y_k = x_k$ and $z_k = z$ for $k=1,\dots, n$.

\subhead 5.2. Bezout and Vi\`ete decompositions of 
the Vandermonde quasideterminants \endsubhead

\proclaim{Theorem 5.2.1 \rm(Bezout decomposition of the Vandermonde
quasideterminant)} Sup\-pose that sequences $x_1,\dots, x_n$ and
$x_1,\dots , x_{n-1}, z$ are independent. Then
$$
V(x_1,\dots, x_n, z) = (z_n-y_n)(z_{n-1}-y_{n-1})\dots(z_1-y_1).\tag 5.2.1
$$
\endproclaim

Note that if $z$ commutes with $x_i, i=1,\dots, n$, then
$$
V(x_1,\dots , x_n, z)= (z-y_n)(z-y_{n-1})\dots (z-y_1).
$$

\proclaim{Theorem 5.2.2 \rm(Vi\`ete decomposition of the Vandermonde
quasideterminant)} For an independent sequence $x_1,\dots, x_n, z$ we have
$$
V(x_1,\dots, x_n, z) = z^n + a_1 z^{n-1} + \dots + a_{n-1} z+ a_n,\tag 5.2.2
$$
where
$$
a_k=(-1)^k\sum_{1\leq i_1<i_2<\dots<i_k\leq n} y_{i_k} y_{i_{k-1}}\dots y_1 .\tag 5.2.3
$$
In particular
$$
\align
a_1&= - (y_1+\dots+ y_n),\\
a_2&= \sum_{1\leq i<j\leq n} y_jy_i,\\
 & \dots\\
a_n&=(-1)^ny_n\dots y_1.
\endalign
$$
\endproclaim

\subhead 5.3.  Proof of Theorem 5.2.2\endsubhead
By induction on $n$ we show that Theorem 5.2.2 follows from
Theorem 5.2.1. For $n=1$ one has $V(x_1,z)= z-x_1$ and
formulas (5.2.1) and (5.2.2) hold. Suppose that these formulas hold for
$m=n-1$. By Theorem 5.2.1
$$
\align
V(x_1,\dots, x_n, z) &= (z_n-y_n) V(x_1,\dots, x_{n-1}, z) \\
&=(V(x_1,\dots, x_{n-1}, z)\cdot z)- (y_n\cdot V (x_1, \dots, x_{n-1}, z)).
\endalign
$$
By induction,
$$
V(x_1,\dots, x_{n-1},z)= z^{n-1} + b_1 z^{n-2} + \dots + b_{n-1},
$$
where
$$
\align
b_1 &=-(y_1+\dots+y_{n-1}),\\
{} &\dots\\
 b_{n-1}&=(-1)^n y_{n-1}\cdot\dots\cdot y_1.
\endalign
$$
Therefore,
$$
\align
V(x_1,\dots,x_n, z) &= z^n+(b_1-y_n) z^{n-1}+(b_2-y_n b_1) z^{n-2} +
\dots -y_n b_n
\\
&= z^n+a_1 z^{n-1} + \dots + a_n,
\endalign
$$
where $a_1,\dots, a_n$ are given by (5.2.3).\qed

\subhead 5.4. Division lemma \endsubhead  To prove Theorem
5.2.1 we need the following result.
\proclaim{Lemma 5.4.1} We have
$$
V(x_1,\dots, x_n, z) = V(\hat x_2,\dots, \hat x_n, \hat z)(z-x_1),
$$
where
$$
\align
\widehat x_k &=(x_k - x_1) x_k(x_k-x_1)^{-1}, \quad k=2,\dots, n\\
\widehat z &=(z-x_1) z(z-x_1)^{-1}.
\endalign
$$
\endproclaim

\demo{Proof} By definition,
$$
V(x_1,\dots, x_n, z) =\vmatrix x_1^n &x^n_2 &\dots &z^n\\
x_1^{n-1} &x_2^{n-1} &\dots &z^{n-1}\\
{}&{} &\dots  &{}\\
x_1 &x_2 &\dots &z\\
1 &1 &\dots &1\endvmatrix_{1, n+1}.
$$
Multiply the $k$-th row by $x_1$ from the left and subtract it from
the $(k-1)$-st row for $k=2,\dots, n$.  Since the quasideterminant does not
change, we have
$$
\align
V(x_1,\dots, x_n, z) &= \vmatrix
0 &x^n_2-x_1x_2^{n-1} &\dots &z^n-x_z^{n-1}\\
0 &x_2^{n-1} - x_1x_2^{n-2} &{} &z^{n-1}-x_1 z^{n-2}\\
\vdots &\vdots  &{} &\vdots\\
0 &x_2-x_1 &{} &z-x_1\\
1 & 1 &{} & 1\endvmatrix_{1, n+1}\\
&=\vmatrix &0 &(x_2-x_1)x_2^{n-1} &\dots &(z-x_1)z^{n-1}\\
\vdots&\vdots &{} &\vdots\\
0 & x_2-x_1& {} &z- x_1\\
1  &1 &{} & 1\endvmatrix_{1, n +1}.
\endalign
$$
Applying to the last quasideterminant Sylvester's theorem with the
element of index $(n+1,1)$ as the pivot we obtain
$$
V(x_1,\dots, x_n, z) = \vmatrix
(x_2-x_1)x_2^{n-1}&\dots &(z-x_1)z^{n-1}\\
\dots &{} &\dots\\
x_2-x_1 &{} &z-x_1\endvmatrix_{1n} .
$$
According to elementary properties of quasideterminants,
multiplying the $k$-th column on the right by 
$(x_{k+1} - x_1)^{-1}$ for $k=1,\dots, n-1$ and
the last column by $(z-x_1)^{-1}$ results in the multiplication
of the value of the quasideterminant on the right by
$(z-x_1)^{-1}$. Therefore,
$$
\align
&V(x_1,\dots, x_n, z)]]
\\&=\vmatrix
(x_2-x_1)x_2^{n-1}(x_2-x_1)^{-1} &\dots &(z-x_1)z^{n-1}(z-x_1)^{-1}\\
\vdots  &{} &\vdots\\
1 & {} & 1\endvmatrix _{1n} (z-x_1) \\
&=\vmatrix
\widehat x_2^{n-1}  &\dots &\widehat z^{n-1}\\
\vdots  & {} &\vdots\\
1   &\dots &1\endvmatrix _{1n}\cdot (z-x_1) = V(\widehat x_2,\dots, \widehat
x_n, \widehat z)\cdot (z-x_1).
\endalign
$$
\enddemo
\subhead 5.5. Proof of Theorem 5.2.1\endsubhead
We proceed by induction on $n$.  By Lemma 5.4.1,
Theorem 5.2.1 is valid for $n=2$.  Also by Lemma 5.4.1,
$$
V(x_1,\dots, x_n, z) = V(\widehat x_2,\dots, \widehat x_n, \widehat z)(z-x_1).\tag 5.5.1
$$
Suppose that our theorem is valid for $m=n-1$.
Then
$$
V(\widehat x_2,\dots,\widehat x_n,\widehat z) = (z'_n - y'_n)\dots (z'_2- y'_2),
$$
where
$$
\align
z'_2 &=\widehat z, \\
y'_2 &= \widehat x_2, \\
z'_k &=V(\widehat x_2,\dots, \widehat x_{k-1},\widehat z) \widehat z V^{-1}(\widehat
x_2,\dots,\widehat x_{k-1},\widehat z),\\
y'_k &= V(\widehat x_2,\dots, \widehat x_k) \widehat x_k V^{-1}(\widehat x_2,\dots,
\widehat x_k) \text{ for } k=3,\dots, n.
\endalign
$$
It suffices to show that
$z'_k = z_k$ and $y'_k=y_k$ for $k=2,\dots, n$.
For $k=2$ this is obvious.  By Lemma 5.4.1,
$$
V(\widehat x_2,\dots, \widehat x_{k-1}, \widehat z)= V(x_1, \dots, x_{k-1}, z)(z-x_1)^{-1},
$$
and by definition $\widehat z=(z-x_1) z(z-x_1)^{-1}$.  So,
$$
\align
z'_k&=\{V(x_1,\dots, x_{k-1},z)(z-x_1)^{-1}\}(z-x_1) z(z-x_1)^{-1}\\
&\quad\times\{(z-x_1) V^{-1}(x_1,\dots, x_{k-1}, z)\}= z_k\text{ for }
k=3,\dots, n.
\endalign
$$
Similarly, $y'_k=y_k$ for $k=3,\dots, n$ and from (5.5.1)
we have
$$
V(x_1,\dots, x_n, z) = (z_n- y_n)\dots (z_2-y_2)(z_1- y_1).
\qquad \qed
$$

\subhead 5.6. Another expression for the coefficients in Vi\`ete decomposition\endsubhead
Ano\-ther expression for the coefficients $a_1,\dots ,a_n$ in Vi\`ete
decomposition of $V(x_1,\dots, x_n, z)$ can be obtained from Proposition 1.5.1.
\proclaim{Theorem 5.6.1 \rm \cite{GKLLRT}} We have
$$
V(x_1,\dots, x_n,z) = z^n+a_1 z^{n-1}+ \dots + a_n,
$$
where for $k=1,\dots, n$
$$
a_k=-\vmatrix
x_1^n&\dots &x_n^n\\
{}&\dots&{}\\
x_1^{n-k+1}&\dots& x_n^{n-k+1}\\
x_1^{n-k-1}&\dots &x_n^{n-k-1}\\
{}&\dots &{}\\
1&\dots &1\endvmatrix _{1n}\cdot
\vmatrix
x_1^{n-1}&\dots&x_n^{n-1}\\
{}&\dots&{}\\
x_1^{n-k}&\dots &x_n^{n-k}\\
{}&\dots&{}\\
1&\dots&1\endvmatrix^{-1} _{kn}. \tag 5.6.1
$$
\endproclaim

>From Theorem 5.6.1 we will get the Bezout and Vi\`ete formulas expressing
the coefficients of the equation
$$
z^n + a_1 z^{n-1} + \dots + a_n = 0.\tag 5.6.2
$$
as polynomials in $x_1,\dots,x_n$ conjugated by Vandermonde determinants.

\subhead 5.7. The Bezout and Vi\`ete Theorems\endsubhead
Recall that the set of elements $x_1,\dots,x_n$ of a ring with unit is
independent if all Vandermonde quasideterminants $V(x_{i_1},\dots , x_{i_k})$ 
for $k\geq 2$ are defined and invertible.

\proclaim{Lemma 5.7.1} Suppose that $x_1,\dots, x_n$ is an independent
set of roots of equation (5.6.2).  Then the coefficients $a_1,\dots, a_n$
can be written in the form (5.6.1).
\endproclaim

\demo{Proof}Consider the system of right linear equations
$$
x^n_i + a_1x^{n-1}_i+\dots+ a_{n-1}x_i+ a_n = 0, \quad i=1,\dots, n
$$
in variables $a_1,\dots ,a_n$ and use Cramer's rule.
\qed
\enddemo

\proclaim{Theorem 5.7.2 \rm (Noncommutative Bezout Theorem)}  Let
$x_1,\dots,x_n$ be an independent set of roots of equation (5.6.2).  In
notations of Theorem 5.2.1,
$$
z^n+a_1 z^{n-1} + \dots+ a_n=(z_n-y_n)\dots(z_1-y_1).
$$
\endproclaim
\demo{Proof} Use Lemma 5.7.1, Theorem 5.6.1 and Theorem 5.2.1.
\qed
\enddemo

\proclaim{Theorem 5.7.2 \rm(Noncommutative Vi\`ete Theorem,
see \cite{GR3})}
Let $x_1,\dots, x_n$ be an independent set of roots of equation (5.6.2).
Then the coefficients $a_1,\dots, a_n$ of the equation are given by
formulas (5.2.3).
\endproclaim
\demo{Proof} Use Lemma 3.2.1, Theorem 3.1.4 and Theorem 3.1.2. \qed
\enddemo

A different proof of this theorem, using
differential operators, appeared in \cite{EGR}. Another noncommutative
version of the Vi\`ete Theorem, based on notions of traces and determinants,
was given by Connes and Schwarz in \cite{CS}.

\head 6. Noncommutative symmetric functions\endhead
\medskip

General theory of noncommutative symmetric functions was developed in the paper
\cite{GKLLRT}. In fact, \cite{GKLLRT} was devoted to the study of different
systems of multiplicative and linear generators in a free algebra $\bold S\bold
y\bold m$ generated by a system of noncommuting variables $\Lambda _i,
i=1,2,\dots $. In \cite{GKLLRT} these variables were called elementary
symmetric functions, but the theory was developed independently of the origin
of $\Lambda _i$. Thus,  in [GKLLRT] only a formal theory of noncommutative
symmetric functions ``without variables" was introduced. The real theory of
noncommutative symmetric functions got ``the right to exist" only after the
corresponding variables were introduced in \cite {GR3, GR4} following the
Viet\'e theorem and the basic theorem in the theory of noncommutative
symmetric functions has been proved in \cite {Wi}.

In this section we apply the general theory to noncommutative symmetric
functions generated by specific $\Lambda _i$. As in the commutative case, they
depend of a set of roots of a polynomial equation.

\subhead 6.1. Formal noncommutative symmetric functions \endsubhead
This theory was started in \cite{GKLLRT} and developed in
several papers (see, for example, \cite{KLT, LST}). An extensive review
was given in \cite {Thi}. Here we just recall some basic constructions.

The algebra $\bSym$ is a free graded associative algebra over
a field $F$ generated by an infinite sequence of variables
$(\Lambda _k)$, $\deg \Lambda_k=k$, $k\geq 1$. The homogeneous
component of degree $n$ is denoted by $\bold S\bold y\bold m_n$.The direct
sum $\oplus _{n\geq 1}\bold S\bold y\bold m_n$ is denoted by
$\bold S\bold y\bold m_+$.
Initially the $\Lambda _k$'s were regarded as as the elementary symmetric
functions of some virtual set of arguments. A natural set
of arguments was found later, see \cite {GR3, GR4}.

Recall some properties of the algebra $\Sym _N$  of symmetric commutative 
polynomials in variables $t_1,\dots , t_N$. The algebra $\Sym_N$ has a natural
grading, $\deg t_i=1$, $i=1,\dots , N$, and is freely generated by the
elementary symmetric functions $e_1(N)=\sum _it_i$,
$e_2(N)=\sum _{i<j}t_i t_j$,\dots, $e_N(N)=t_1t_2\dots t_N$. (There are 
other natural sets of generators in $\Sym _N$). 
Setting $t_N=0$ one gets a canonical epimorphism of graded algebras
$p_N:\Sym_N\to \Sym_{N-1}$.

The projective limit of graded algebras $\Sym _N$ with respect to the system $\{p_N\}$ 
is called the algebra
$\Sym$ of symmetric functions in infinite set of variables $t_1, t_2,\dots $ 
(see \cite{Mac}). 
One can view the algebra $\Sym $ as a free commutative algebra generated
by formal series $e_1=\sum _it_i$, $e_2=\sum _{i<j}t_it_j$, \dots,
$e_k=\sum {i_1<\dots <i_k}t_{i_1}\dots t_{i_k},\ \dots$. 
The series $e_k$ is called the $k$-th elementary symmetric function in
$t_1, t_2, \dots $. 
 
In $\text{Sym}$, there are also other standard sets of generators 
(see, e.g., \cite{Mac}). The most common among them are the complete symmetric 
functions 
$(h_k)_{k\geq 1}$ and the power symmetric functions $(p_k)_{k\geq 1}$. To express them
in terms of $(e_k)$ one can use generating functions. Namely, set
$e_0=h_0=1$. Let $\tau $ be a formal variable. Set $\lambda (\tau
)=\sum _{k\geq 0}h_k\tau _k$, $\sigma (\tau )=\sum _{k\geq 0}$,
$\psi (\tau )=\sum _{k\geq 1}p_kt^{k-1}$. Then
$$
\align
\lambda (\tau )&=\sigma (-\tau)=1,\\
\psi (\tau )&=\frac {d}{d\tau }\log \sigma (\tau ).
\endalign
$$

Define the canonical epimorphism $\pi :\bold S\bold y\bold m
\to\text {Sym}$ by setting $\pi (\Lambda _k)=e_k$,
$k\geq 1$. Let $\Cal I_N$ be an ideal $\bold S\bold y\bold m$
generated by all $\Lambda _k$, $k>N$. The epimorphism $\pi $
induces the canonical epimorphism $\pi _N$ of $\bold S\bold y\bold
m $ onto the algebra ${\text Sym}_N$ of symmetric polynomials in
commuting variables $t_1,\dots , t_N$, $N\geq 1$. Note that 
$\pi_N (\Cal I_N)=0$.

Noncommutative analogs of functions $(h_k)$ and $(p_k)$ can be
constructed in the following way.
Let $\tau $ be a formal variable commuting with all $\Lambda_k$. 
Set $\Lambda _0=1$ and define the generating series
$$
\lambda (\tau):=\sum _{k\geq 0}\Lambda _k\tau ^k.
$$

\definition{Definition 6.1.1} The complete homogeneous symmetric functions
are the coefficients $S_k$ in the generating series
$$
\sigma (\tau ):=\sum _{k\geq 0}S_k\tau ^k=\lambda (-t)^{-1}.\tag
6.1.1a
$$
The power sums symmetric functions of the first kind $\Psi _k$
are the coefficients $\Psi _k$ in the generating series
$$
\sum _{k\geq 1}\Psi_k \tau ^{k-1}:=\sigma (\tau)^{-1}\frac
{d}{dt}\sigma (\tau). \tag 6.1.1b
$$
The power sums symmetric functions of the second kind $\Phi _k$
are defined by
$$
\sum _{k\geq 1}\Phi_k \tau ^{k-1}:=\frac {d}{d\tau }\log \ \sigma
(\tau ).
\tag 6.6.1c
$$
\enddefinition

By using formulas (6.1.1a)--(6.6.1c) one can prove that $\pi _N (S_k)$ is the
$k$-th complete symmetric function and $\pi _N(\Psi _k)=\pi
_N(\Phi _k)$ is the $k$-th power symmetric function in $N$
commuting variables. Note that in the
right-hand sides of (6.1.1b) and (6.1.1c) different noncommutative
analogs of the logarithmic derivative of $\sigma (t)$ are used.

Definition 6.1.1 leads to the following quasideterminantal
formulas.

\proclaim{Proposition 6.1.2} For every $k\geq 1$, one has
$$
\allowdisplaybreaks
\align
S_k&=(-1)^{k-1}\left |\matrix
\Lambda _1&\Lambda _2&\dots &\Lambda _{k-1}&\boxed {\Lambda _k}\\
         1&\Lambda _1&\dots &\Lambda _{k-2}&\Lambda _{k-1}\\
         0&         1&\dots &\Lambda _{k-3}&\Lambda _{k-2}\\
          &          &\dots &              &\\
         0&         0&\dots &             1&\Lambda _1
\endmatrix \right |,
\\
\Lambda _k&=(-1)^{k-1}\left |\matrix
S_1&1&0&\dots &0\\
 S_2&S_1&1&\dots &0\\
 S_3&S_2&S_1&\dots &0\\
    &   &\dots &  &\\
 \boxed {S_k}&S_{k-1}&S_{k-2}&\dots &S_1
\endmatrix \right |,
\\
kS_k&=\left |\matrix
\Psi _1&\Psi _2&\dots &\Psi _{k-1}&\boxed {\Psi _k}\\
     -1&\Psi _1&\dots &\Psi _{k-2}&\Psi _{k-1}\\
      0& -2    &\dots &\Psi _{k-3}&\Psi _{k-2}\\
          &          &\dots &              &\\
         0&         0&\dots & -n+1&\Psi _1
\endmatrix \right |,
\\
k\Lambda _k&=\left |\matrix
\Psi _1&1&0&\dots &0\\
\Psi _2&\Psi_1&2&\dots &0\\
    &   &\dots &  &\\
\boxed {\Psi _k}&\Psi _{k-1}&\Psi _{k-2}&\dots &\Psi _1
\endmatrix \right |,
\\
\Psi _k&=(-1)^{k-1}\left |\matrix
\Lambda _1&2\Lambda _2&\dots &(k-1)\Lambda _{k-1}&\boxed {k\Lambda _k}\\
         1&\Lambda _1&\dots &\Lambda _{k-2}&\Lambda _{k-1}\\
         0&         1&\dots &\Lambda _{k-3}&\Lambda _{k-2}\\
          &          &\dots &              &\\
         0&         0&\dots &             1&\Lambda _1
\endmatrix \right |,
\\
\Psi _k&=\left |\matrix
S_1&1&0&\dots &0\\
2 S_2&S_1&1&\dots &0\\
    &   &\dots &  &\\
\boxed {kS_k}&S_{k-1}&S_{k-2}&\dots &S_1
\endmatrix \right |.
\endalign
$$
\endproclaim

Each of the four sequences $(\Lambda _k)$, $(S_k)$, $(\Psi _k)$,
and $(\Phi _k)$ is a set of generators in $\bold S\bold y\bold m$.
Therefore, each of the four sets of products $F_{i_1}\dots
F_{i_N}$, $i_1,\dots , i_N\geq 1$, where $F_{i_k}$ equals to
$\Lambda _{i_k}$, $S_{i_k}$, $\Psi_{i_k}$, or $\Phi {i_k}$,
is a linear basis in $\bold S\bold y\bold m_+$.
Linear relations between these bases were given in 
\cite{GKLLRT}.

 Another important example of a linear basis in $\bold
S\bold y\bold m_+$ is given by {\it ribbon Schur functions}.

\subhead 6.2. Ribbon Schur functions \endsubhead
Commutative ribbon Schur functions were defined by MacMahon 
\cite{M}. Here we follow his ideas.

Let $I=(i_1,\dots , i_k)$, $i_1,\dots , i_k\geq 1$, be an ordered
set. 

\definition{Definition 6.2.1 \rm \cite {GKLLRT}}
The ribbon Schur function $R_I$ is defined by the formula
$$
R_I=(-1)^{k-1}\left |\matrix
S_{i_1}&S_{i_1+i_2}&S_{i_1+i_2+i_3}&\dots &\boxed {S_{i_1+\dots +i_k}}\\
1      &S_{i_2}&S_{i_2+i_3}&\dots &S_{i_2+\dots +i_k}\\
0      &1      &S_{i_3}&\dots &S_{i_3+\dots +i_k}\\
       &       &       &\dots &                  \\
0      &0      &0      &\dots &S_{i_k}\endmatrix \right |.
$$
\enddefinition

Definition 6.2.1 allows us to express $R_I$'s as polynomials
in $S_k$'s. To do this we need the following ordering of
sets of integers.

Let $I=(i_1,\dots , i_r)$ and $J=(j_1,\dots , j_s)$. We say that
$I\leq J$ if $i_1=j_1+j_2+\dots +j_{t_1}$, $i_2=j_{t_1+1}+\dots
+j_{t_2}$, $\dots $, $i_s=j_{t_{s-1}+1}+\dots +j_s$. For example,
if $I\leq (12)$, then $I=(12)$ or $I=(3)$. If $I\leq (321)$, then
$I$ is equal to one of the sets $(321)$, $(51)$, $(33)$, or $(6)$.

For $I=(i_1,\dots , i_r)$ set $l(I)=r$ and $S^I=
S_{i_1}S_{i_2}\dots S_{i_r}$.

\proclaim {Proposition 6.2.2} \rm \cite {GKLLRT} (p. 254)
$$
R_J=\sum _{I\leq J}(-1)^{l(J)-l(I)}S^I.
$$
\endproclaim

\example{Example} $R_{123}=S_6-S^2_3-S_1S_5+S_1S_2S_3$.
\endexample

Definition 6.1.2 implies that $R_I=S_m$ for $I=\{m\}$ and
$R_I=\Lambda _k$ for $i_1=\dots =i_k=1$. For each $N$ the
homomorphism $\pi _N$ maps $R_I$ to the corresponding
MacMahon ribbon Schur function.
 
In \cite {GKLLRT} similar formulas expressing $R_I$ as
quasideterminants of matrices with entries $\Lambda _k$, as well
as linear relations with different bases in $\bold S\bold y\bold
m_+$ defined in Section 6.1, are given.

Natural bases in algebra $\Sym$ of commutative symmetric functions
are indexed by weakly decreasing (or, weakly increasing) finite
sequences of integers. Examples are products of elementary
symmetric functions $e_{i_1}\dots e_{i_k}$ where $i_1\geq i_2
\dots \geq i_k$ and Schur functions $s_{\lambda }$ where
$\lambda =(i_1,\dots , i_k)$.
The following theorem gives a natural basis in the algebra
of noncommutative symmetric functions. Elements of this basis
are indexed by all finite sequences of integers.
 
\proclaim{Theorem 6.2.3 \rm \cite {GKLLRT}} The ribbon Schur
functions $R_I$ form a linear basis in $\bold S\bold
y\bold m$.
\endproclaim

Let $\pi:\bSym\to\Sym$
be the canonical morphism. Then it is known (see \cite{M})
that the commutative ribbon Schur functions $\pi(R_I)$ are not linearly
independent. For example, commutative ribbon Schur functions
defined by sets $(ij)$ and $(ji)$ coincide. This means that the kernel
$\Ker \pi$ is nontrivial.

\remark {Remark}  In the commutative case, ribbon Schur functions 
$\pi (R_I)$ with with weakly decreasing $I$ constitute a basis in 
the space of symmetric functions. However,
this basis is not frequently used.
\endremark

The description of the kernel $\Ker \pi$ in terms of ribbon Schur functions 
is given by the following theorem.

For an ordered set $I$ is denote by $u(I)$ the
corresponding unordered set.

\proclaim{Theorem 6.2.4} The kernel of $\pi $ is linearly generated
by the elements
$$
\Delta _{J,J'}=\sum _{I\leq J}R_I - \sum _{I'\leq J'}R_{I'}
$$
for all $J$, $J'$ such that $u(J)=u(J')$.
\endproclaim

\example{Example} 1. Let $J=(12)$, $J'=(21)$. Then
$\Delta _{J,J'}=(R_{12}+R_3) - (R_{21}+R_3)=R_{12}-R_{21}$
and $\pi (R_{12})=\pi(R_{21})$.

2. Let $J=(123)$, $J'=(213)$. Then $\Delta _{J,J'}=(R_{123}+R_{33}
+R_{15}+R_{6}) - (R_{213}+R_{33}+R_{24}+R_{6})=R_{123}+R_{15}-R_{213}-R_{24}$. 
This shows, in particular, that $\pi (R_{123})-\pi(R_{213})=\pi(R_{24}) - 
\pi(R_{15})\ne 0$.
\endexample

The homological relations for quasideterminants imply the
multiplication rule for the ribbon Schur functions.
Let $I=(i_1,\dots , i_r)$, $J=(j_1,\dots , j_s)$, $i_p\geq 1$, $j_q\geq 1$
for all $p,q$. Set $I+J=(i_1,\dots , i_{r-1}, i_r+j_1, j_2,\dots , j_s)$
and $I\cdot J=(i_1,\dots ,i_r,j_1,\dots , j_s)$.

The following picture illustrates this definition (and explains the origin of the name "ribbon Schur functions").
To each ordered set $I=\{i_1,i_2,\dots,i_k\}$ we can associate a ribbon, 
i.e., a sequence of square cells on the square rules paper starting at 
the square $(0,0)$ and
going right and down, with $i_1$ squares in the first column, $i_2$ 
squares in the second column, and so on, see Figure 1 for the ribbons
corresponding to $I=(2,1,3)$ and $J=(3,1,2)$. 
Then the 
construction of ribbons $I+J$ and $I\cdot J$ has a 
simple geometric meaning as shown in Figure 2 for $I=(2,1,3)$ and
$J=(3,1,2)$. 

\midinsert
\centerline{
\epsfbox{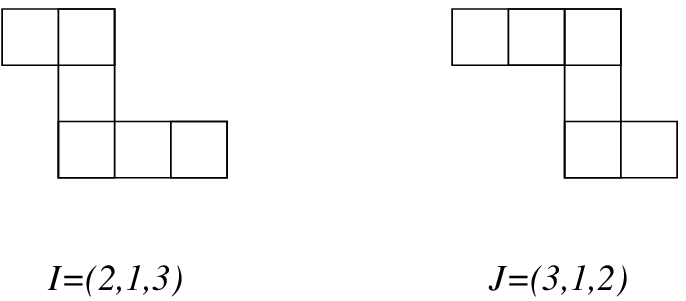}}
\botcaption{Figure 1}
\endcaption
\endinsert

\midinsert
\centerline{
\epsfbox{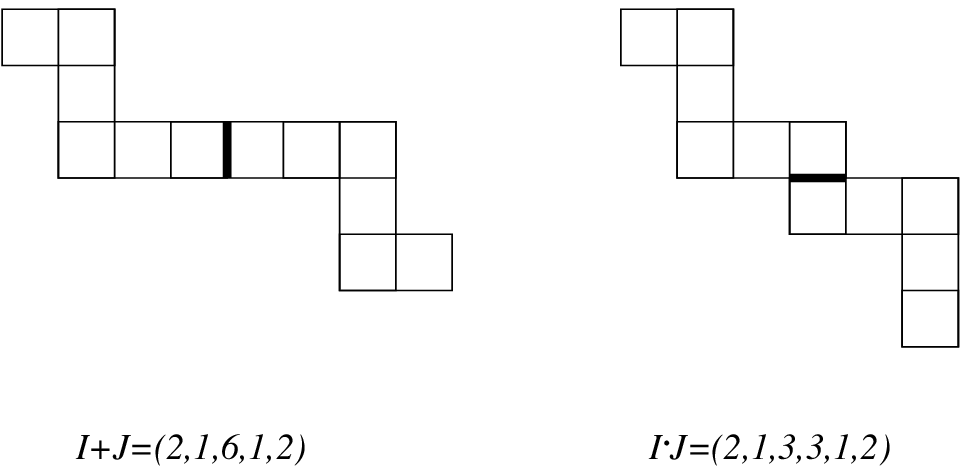}}
\botcaption{Figure 2}
\endcaption
\endinsert

\proclaim{Theorem 6.2.5 \rm \cite {GKLLRT}} We have
$$
R_IR_J=R_{I+J}+R_{I\cdot J}.
$$
\endproclaim

The commutative version of this multiplication rule is due to MacMahon.

Naturality of ribbon Schur functions $R_I$ can be explained
in terms of the following construction.

\subhead
6.3. Algebras with two multiplications
\endsubhead
The relations between the functions $\Lambda _k$ and
the functions $S_k$ can be illuminated by noting that the ideal $\bold S\bold y\bold m_+$
has two natural associative multiplications $*_1$ and $*_2$.  In terms of ribbon Schur functions
it can be given as $R_I*_1R_J=R_{I.J}$ and $R_I*_2R_J=R_{I+J}$.
We formalize this notion as follows.

\definition{Definition 6.3.1} A linear space $A$ with two bilinear products $\circ_1$ and $\circ_2$
is called a biassociative algebra with products $\circ_1$ and $\circ_2$
if 
$$
(a \ \circ_i \ b) \ \circ_j \ c = a \ \circ_i \ (b \ \circ_j \ c)
$$ 
for all $a,b,c \in A$ and all $i,j \in \{1,2\}.$
\enddefinition

Note that if the products $\circ_1$ and $\circ_2$ in a biassociative algebra
$A$ have a common identity element $1$  (i.e., if $1 \ \circ_i \ a = a \
\circ_i 1 = a$ for all $a \in A$ and $i=1,2$, then
$$
a \ \circ_1 \ b = (a \ \circ_2 1 )\ \circ_1 b = a \ \circ_2 (1 \ \circ_1 b) = a
\ \circ_2 b 
$$ 
for all $a,b \in A$ and so $\circ_1 = \circ_2$.

Note also that if $A$ is a biassociative algebra with two products $\circ_1$
and $\circ_2$, then for $r,s \in F$ one can define the linear combination
$\circ_{r,s} 
= r\circ_1 + s\circ_2$ by the formula
$$
a \ \circ_{r,s} \ b = r(a \ \circ_1 \ b) + s(a \ \circ_2 \ b), \qquad a,b\in A. 
$$
Then $A$ is a biassociative algebra with the
products $\circ_{r,s}$ and  $\circ_{t,u}$ for each
$r,s,t,u \in F$.

Jacobson's discussion of isotopy and homotopy of Jordan algebras
(see \cite{Ja2}, p.~56,ff) shows that if $A$ is an associative algebra
with the product $\circ$ and $\circ_a$ for $a\in A$ is defined by
the formula
$$
b \ \circ_a c = b \ \circ a \ \circ c,
$$ 
then $A$ is a biassociative algebra with
the products $\circ$ and $\circ_a$.

We now endow the ideal $\bSym_+\subset \bSym$ with the structure of a
biassociative algebra in two different ways. Recall that the nontrivial
monomials $(\Lambda _{i_1}\dots \Lambda_{i_r})$
as well as the nontrivial monomials $(S _{i_1}\dots S _{i_r})$ form linear
bases in $\bSym+$.

\definition{Definition 6.3.2} Define the linear map
$*_1: \bold S\bold y\bold m _+\otimes \bold S\bold y\bold m _+
\rightarrow \bold S\bold y\bold m_+$
by
$$
(\Lambda _{i_1}\dots \Lambda _{i_r})*_1
(\Lambda _{j_1}\dots \Lambda _{j_s})=
\Lambda _{i_1}\dots \Lambda _{i_{r-1}}
\Lambda _{i_r+j_1}\Lambda _{j_2}\dots \Lambda _{j_s}
$$
and the linear map
$*_2: \bold S\bold y\bold m_+\otimes \bold S\bold y\bold m_
+\rightarrow \bold S\bold y\bold m_+$
by
$$
(S _{i_1}\dots S _{i_r})*_2
(S _{j_1}\dots S _{j_s})=
S _{i_1}\dots S _{i_{r-1}}
S _{i_r+j_1}S _{j_2}\dots S _{j_s}.
$$
\enddefinition

Write $ab=a*_0b$ for $a,b\in \bold S\bold y\bold m _+$. Then it is clear
that
$a*_i(b*_jc)=(a*_ib)*_jc$ for all $a,b,c\in \bSym_+$ and
$i,j=0,1$ or $i,j=0,2$. Thus we have the following result.

\proclaim{Lemma 6.3.3} $\bSym_+$ is a biassociative algebra with products
$*_0$ and $*_1$ and also a biassociative algebra with products $*_0$ and $*_2$.
\endproclaim

In fact, $*_0$, $*_1$ and $*_2$ are closely related.

The following Lemma is just a restatement of Theorem 6.2.5.

\proclaim{Lemma 6.3.4} $*_0=*_1 + *_2$.
\endproclaim

\demo{Proof} We have
$$
\lambda (-t)^{-1}=\Big(1+\sum _{i>0}(-1)^i\Lambda _it^i\Big)^{-1}=
1+\sum _{j>0}\sum _{i_1+\dots +i_l=j}(-1)^{l+j}
\Lambda _{i_1}\dots \Lambda _{i_l}t^j.
$$
Since $\Lambda _i=\Lambda _1*_1\Lambda _1*_1\dots *_1\Lambda _1$,
where there are $i-1$ occurences of $*_1$, the
coefficient at $t^j$ in $\lambda (-t)^{-1}$ is
$$
\sum _{u_1, \dots, u_{j-1}\in \{0,1\}}(-1)^k
\Lambda _1*_{u_1}\Lambda _1*_{u_2}\dots *_{u_{j-1}}\Lambda _1,
$$
where $k$ is the number of $u_t$ equal to $1$.

Since $1+\sum S_it^i=\lambda (-t)^{-1}$ we have 
$$
S_j=\sum _{u_1, \dots u_{j-1}\in \{0,1\}}(-1)^k
\Lambda _1*_{u_1}\Lambda _1*_{u_2}\dots *_{u_{j-1}}\Lambda _1.
$$
Therefore $S_i*_0S_j- S_i*_1S_j= S_i(*_0 - *_1)S_j = S_{i+j} = S_i*_2S_j$ and 
$*_2=*_0-*_1$, as required.
\qed\enddemo

Now let $U$ be the two-dimensional vector space with basis $\{u_0,u_1\}$ and 
$$
F\langle U\rangle = \sum_{k \ge 0} F\langle U\rangle_k
$$
the (graded) free associative algebra on $U$, with the homogeneous components
$$
F\langle U\rangle_k = U^{\otimes k}.
$$ 
We use the products $*_1$ and $*_2$ to define two isomorphisms,
$\phi_1$ and $\phi_2$,
of $F\langle U\rangle$ to $\bSym _+$. Namely, for a basis element 
$ u_{i_1}\dots u_{i_l}\in F\langle U\rangle_k$ set
$$
\phi_1 : u_{i_1}\dots u_{i_l}\mapsto
\Lambda _1*_{i_1}\Lambda _1*_{i_2}\dots *_{i_l}\Lambda _1\in
(\bSym_+)_{l+1}
$$
and
$$
\phi_2 : u_{i_1}\dots u_{i_l}\mapsto
\Lambda _1*_{j_1}\Lambda _1*_{j_2}\dots *_{j_l}\Lambda _1\in
(\bSym_+)_{l+1}
$$
where $j_t = 0$ if $i_t = 0$ and $j_t = 2 $ if $i_t = 1.$
Note that $\phi_1$ and $\phi_2$ shift degree.

Define the involution $\theta$ of $U$ by $\theta (u_0) = u_0$ and $\theta (u_1) = u_0 - u_1$.
Then $\theta$ extends to an automorphism $\Theta$ of $F\langle U\rangle$
and the restriction $\Theta_k$ of this automorphism to $F\langle U\rangle_k$ is the $k$-th tensor power of $\theta$.
Clearly $\phi_1 \Theta = \phi_2$ and so we recover
Proposition 4.13 in \cite{GKLLRT}, which describes, in terms of tensor powers, the relation between the bases of $\bSym_+$
consisting of nontrivial monomials in $\Lambda_i$ and of nontrivial monomials in $S_i$.

Similarly, taking the identity
$$
a^{n-1}+(-1)^n b^{n-1}+\sum _{k=1}^{n-1}(-1)^{n-k}
a^{k-1}(a+b)b^{n-k-1}=0,
$$
valid in any associative algebra, setting
$a=u_0-u_1$, $b=u_1$, and applying $\phi_1 $, we obtain
the identity 
$$
0=\sum _{k=0}^n (-1)^{n-k}\Lambda _k S_{n-k}
$$
between the elementary and complete
symmetric functions (Proposition 3.3 in \cite {GKLLRT}).
Using Proposition 6.1.2, one can express these identities in terms of quasideterminants.

%We introduce a natural bijective correspondence between monomials
%in $F(U)$ and ribbon Schur functions.

%Write each finite sequence $\omega $ of symbols $1$ and $2$ as
%$(2\dots 2)(1)(2\dots 2)(1)\dots $ bracketing maximal
%number of $2$'s and each $1$. Label the sequence as
%$I_{\omega }=(i_1,\dots , i_k)$ where $i_m$ is the number of elements in
%the $m$-th bracket. For example, sequences
%$2212$ and $12212$ will be labelled as $2,1,1$ and
%and $1,2,1,1$ correspondingly. The sequence $222$ is labelled
%as $I=(3)$ and the sequence $111$ is labelled as $I=(111)$.

%To each monomial $u_{j_1}u_{j_2}\dots u_{j_l}\in F(U)$
%it corresponds a sequence $j_1,\dots , j_l$ of $1$'s and $2$'s.
%To each such sequence it corresponds and ordered set
%$I=(i_1,\dots , i_k)$.

%We constructed
%bijective correspondences between monomials
%in $F\langle U\rangle $ and finite sequences in $1,2$'s and
%between the finite sequences and ordered sets $I=(i_1,\dots , i_k)$.
%Thus, to any monomial $u_{j_1}\dots u_{j_l}\in F\langle U\rangle $
%corresponds a ribbon Schur function $R_I$ establishing a bijective
%correspondence between two natural linear basises.

\subhead 6.4. Quasi-Schur functions \endsubhead
Quasi-Schur functions were defined in \cite{GKLLRT}. They are 
elements not of $\bSym$ but of the free-skew field
generated by $S_1,S_2\dots $. Let $I=(i_1,\dots , i_k)$.

\definition{Definition 6.4.1} Define $\check S_I$
by the formula
$$
\check S_I=(-1)^{k-1}\left |\matrix
S_{i_1}&S_{i_2+1}&\dots &\boxed {S_{i_k+k-1}}\\
S_{i_1-1}&S_{i_2}&\dots &S_{i_k+k-2}\\
         &       &\dots & \\
S_{i_1-k+1}&S_{i_2-k+2}&\dots &S_{i_k}\endmatrix \right | .\tag 6.4.1
$$
If $I=(i_1,\dots , i_k)$ is a partition,
i.e., a weakly increasing sequence of nonnegative sequences,
the element $\check S_I$ is called a quasi-Schur function.
For an arbitrary set $I$, $\check S_I$ is called a generalized quasi-Schur function.
\enddefinition

In particular, $\check S_{\{k\}}=S_k$ and $\check S_{\{1^k\}}=\Lambda _k$,
where $1^k=(1\dots 1)$ with $k$ occurences of $1$.

Definition 6.4.1 is a noncommutative analog of Jacobi--Trudi formula.
In the commutative case,
$\check S_I$ for a partition $I$ is the {\it ratio} of two Schur functions
$S_I/S_J$, where $J=(i_1-1,\dots , i_{k-1}-1)$. It shows that in general 
$\check S$ cannot be represented as a polynomial in $S_k$.

\remark{Remark} The homological relations and the trasformation properties
of quasideterminants imply that any generalized quasi-Schur function
$\check S_I$ can
be expressed as a rational function in the quasi-Schur functions.
For example,
$$
\check S_{42}=-\left |\matrix S_4&\boxed {S_3}\\S_3&S_2\endmatrix \right |
=-\left |\matrix \boxed {S_3}&S_4\\S_2&S_3\endmatrix \right |
=\left |\matrix S_3&\boxed {S_4}\\S_2&S_3\endmatrix \right |
S_3^{-1}S_2=\check S_{33}S_3^{-1}S_2.
$$
\endremark

\subhead 6.5. Symmetric functions in noncommutative variables \endsubhead
We fix $n$ independent indeterminants $x_1,
x_2, \dots, x_n$ and construct new variables $y_1,\dots, y_n$ which are
{\it rational} functions in $x_1, \dots, x_n$ as follows. Recall that
in 5.1 we defined the Vandermonde quasideterminant
$$
V(x_1,\dots , x_k)=\vmatrix 
x_1^{k-1} &\dots&\boxed {x_k^{k-1}}\\
x_1^{k-2}&\dots &x_k^{k-2}\\
{}&\dots&{}\\
1&\dots&1\endvmatrix .
$$
Set
$$
\align
y_1 &= x_1,\\
y_2 &=V(x_1, x_2)x_2V(x_1, x_2)^{-1}
=(x_2-x_1)x_2(x_2-x_1)^{-1},\\
&\dots \\
y_n&=V(x_1,\dots , x_n)x_nV(x_1,\dots , x_n)^{-1}
\endalign
$$

In the commutative case $x_i=y_i$, $i=1,\dots n$. In the noncommutative case
$x_i$ and $y_i$ are obviously different.

\remark{Remark} Consider the free skew-field $R$ generated by $x_1,\dots , x_n$.
Define on $R$ differential operators $\partial _i$ by formula
$\partial _ix_j=\delta _{ij}$ and
the Leibniz rule $\partial _i (fg)=\partial _i(f)g +
f\partial _i(g)$, for $i=1,\dots n$.
It easy to see that $\partial _iy_j\neq \delta_{ij}$. However, denote
$\partial =\partial _1 + \dots +\partial _n$. Then
$$ 
\partial (V(x_1,\dots , x_k))=0, \quad k=2,\dots n
$$
and
$$ \partial (y_i)=\partial (x_i)=1, \quad i=1,\dots ,n.$$
\endremark

\subhead Elementary symmetric functions\endsubhead
\definition{Defintion 6.5.1} The functions
$$
\align
\Lambda_1(x_1,\dots, x_n) &= y_1 + y_2+\dots + y_n,\\
\Lambda_2(x_1,\dots, x_n) &= \sum_{i<j} y_jy_i,\\
      {}&\dots    {}\\
\Lambda_n(x_1,\dots, x_n) &= y_n\dots y_1
\endalign
$$
are called elementary symmetric functions in $x_1, \dots, x_n.$
\enddefinition

In the commutative case these functions are the standard elementary
symmetric functions of $x_1,\dots, x_n$.
By the Noncommutative Vi\`ete Theorem (Theorem 5.7.1),
$\Lambda_i(x_1,\dots, x_n) = (-1)^i a_i$, $i=1,\dots,n$, where
$x_1, \dots, x_n$ are the roots of the equation
$$
x^n+a_1 x^{n-1} + \dots + a_{n-1} x+ a_n = 0.
$$
This implies
\proclaim{Proposition 6.5.2} The functions $\Lambda_i(x_1,\dots, x_n)$ are
symmetric in $x_1, \dots, x_n$.
\endproclaim

Denote by $\overline {\bSym}_n$ the subalgebra of
the algebra of rational functions in $x_1,\dots , x_n$ generated
by $\Lambda _k(x_1,\dots , x_n)$, $k=1,\dots , n$.
Define the surjective homomorphism
$$
\phi : \bold S\bold y\bold m\to\overline {\bSym}_n 
\tag 6.5.1
$$
by setting $\phi (\Lambda _k)=\Lambda _k(x_1,\dots ,x_n)$.

\proclaim{Theorem 6.5.3} The kernel of $\phi $ is generated by
$\Lambda _k$ for $k>n$.
\endproclaim

\remark{Remark} The order of $y_1,\dots, y_n$ is essential in the
definition of $\Lambda_i(x_1,\dots , x_n)$, $i=1,\dots, n$.
For example, $\Lambda_2 (x_1,x_2) = y_2y_1$ is symmetric in $x_1, x_2$, 
whereas the product $y_1 y_2$ is not symmetric.
To see this set $\Delta =x_2-x_1$. The symmetricity in $x_1, x_2$ of the product
$y_1 y_2$ would imply that 
$x_1\Delta^2=\Delta ^2x_1$.
\endremark

\subhead Complete symmetric functions\endsubhead

\definition{Definition 6.5.4}
The functions
$$
S_k(x_1,\dots, x_n) = \sum_{i_1\leq i_2\leq\dots\leq i_k}y_{i_1}\dots
y_{i_k},\quad k=1,2,3,\dots
$$
are called complete symmetric functions in $x_1,\dots , x_n$.
\enddefinition

In the commutative case these functions are the standard complete
symmetric functions in $x_1,\dots, x_n$.

Let $t$ be a formal variable commuting with $x_i$, $i=1,\dots n$.
Define the generating functions
$$
\align
\lambda (t)&=1+\Lambda (x_1,\dots ,x_n)t+\dots
+\Lambda _n (x_1,\dots,x_n)t^n\\
\sigma (y)&=1+\sum _i S_i(x_1,\dots , x_n)t^i =\lambda (-t)^{-1}.
\endalign
$$

\proclaim{Proposition 6.5.5} We have 
$$
\sigma (t)\lambda (-t)=1.
$$
\endproclaim

In the commutative case, $S_k(x_1,\dots , x_n)$ are the standard
complete symmetric functions.

\proclaim{Proposition 6.5.6} The functions $S_k(x_1,\dots, x_n)$ are
symmetric in
$x_1,\dots, x_n$.
\endproclaim

\demo{Proof} Use Proposition 6.5.5, Theorem 6.5.3 and Proposition 6.1.2.
\qed
\enddemo

\remark{Remark} The order of elements $y_s$ in the definition of $S_k$ is essential:
$S_2(x_1,x_2) = y^2_1+ y_1y_2 + y_2^2$ is symmetric in $x_1, \dots ,x_n$ whereas
$y_1^2 + y_2y_1 + y_2^2$ is not symmetric (cf. the remark after Theorem 6.5.3).
\endremark

\subhead Ribbon Schur functions\endsubhead
We define "ribbon Schur functions with arguments" $R_I(x_1,\dots , x_n)$
similarly to Definition 6.2.1, replacing $R_I$ with $R_I(x_1,\dots , x_n)$
and $S_k$ with $S_k(x_1,\dots , x_n)$. Evidently,
ribbon Schur functions $R_I(x_1,\dots, x_n)$ are
symmetric in $x_1,\dots , x_n$ and form a linear basis in
$\overline {\bSym}_n$.

\proclaim{Proposition 6.5.7} We have
$R_I(x_1,\dots , x_n)=\phi (R_I)$, where
$\phi $ is defined by formula (6.5.1).
\endproclaim

To express $R_I(x_1,\dots , x_n)$ as a sum of monomials in
$y_1,\dots , y_n$ we need 
some notation. Let $w=a_{i_1}\dots a_{i_k}$ be a word in ordered
letters $a_1<\dots <a_n$.  An integer $m$ is
called a {\it descent\/} of $w$ if $1\leq m\leq k-1$ and $i_m > i_{m+1}$.
Let $M(w)$ be the set of all descents of $w$.

Let $J=(j_i,\dots, j_k)$ be a set of positive integers.

\proclaim{Theorem 6.5.8}
$$
R_J(x_1,\dots, x_n) = \sum y_{i_1}\dots y_{i_m}, \tag 6.5.2
$$
where the sum is taken over all words $w=y_{i_1} \dots y_{i_m}$ such that
$M(w)=\{j_1, j_1+j_2,\dots, j_1 + j_2+\dots+ j_{k-1}\}$.
\endproclaim

The proof of the theorem was essentially given in \cite{GKLLRT}, Section VII.

\subhead 6.6. Main Theorem for noncommutative symmetric functions\endsubhead
In the commutative case the classical main theorem of the theory of symmetric
functions says that every symmetric polynomial of $n$ variables
is a polynomial of (elementary) symmetric functions of these
variables. Its analogue for a noncommutative case is given by
the following theorem. Denote $\Lambda _k(x_1,\dots , x_n)$ as
$\Lambda _k(X)$, $k=1,\dots , n$.

Recall that in the previous section we defined the elements
$y_k$ by the formulas $y_1=x_1$, $y_k=V(x_1,\dots , x_k)x_k
V(x_1,\dots , x_k)^{-1}$ for $k=2,\dots , n$.    

\proclaim{Theorem 6.6.1 \rm(\cite {Wi})}
Let a polynomial $P(y_1, \dots, y_n)$ over $\Bbb Q$ be
symmetric in $x_1, \dots, x_n$. Then $P(y_1,\dots,
y_n)=Q(\Lambda_1(X), \dots, \Lambda_n(X))$, where $Q$ is a noncommutative
polynomial over $\Bbb Q$.
\endproclaim

\remark{Remark} Recall that $P(y_1,\dots , y_n)$ is a polynomial
in $y_i$ and not in $x_i$. We can express this by saying that in the natural variables $x_i$, noncommutative symmetric
polynomials are not polynomials but rational functions.
\endremark

\proclaim{Corollary 6.6.2} A polynomial $P(y_1, \dots, y_n)$ with 
coefficients 
in $\Bbb Q$ is symmetric in $x_1, \dots, x_n$ if and
only if $P(y_1,\dots, y_n)$, viewed as a rational function of $x_1,\dots, x_n$, is a
linear combination of $R_J(x_1,\dots,x_n)$.
\endproclaim

\subhead 6.7. Quasi-Pl\"ucker coordinates of Vandermonde
matrices and symmetric functions\endsubhead
Here we study right quasi-Pl\"ucker coordinates
$r_{ij}^{i_1,\dots ,i_{n-1}}(V_n)$, where  $V_n=(x_j^i)$, $i\in
\Bbb Z$, $j=1,\dots , n$, is the Vandermonde matrix,  $x_j$, $j=1,\dots , n$
are noncommuting variables. The matrix $V_n$ has $n$ columns and infinitely many rows.
The following proposition shows the importance of such coordinates.

\proclaim{Proposition 6.7.1} We have
$$
\gather
r_{n+k-1, n-1}^{0,1,\dots , n-2}(V_n)=S_k(x_1,\dots , x_n),
\quad k=0, 1, 2,\dots ,\\
 r_{n, n-k}^{0,1,\dots , n-k-1, n-k+1,\dots , n-1}(V_n)=
(-1)^{k-1}\Lambda _k(x_1,\dots , x_n),
\quad k=0, 1, \dots, n.
\endgather
$$
\endproclaim

\example {Examples}
$$
\align
S_k(x_1, x_2)&=
\left |  \matrix
x_1^{k+1}&\boxed {x_2^{k+1}}\\
1                &                    1
\endmatrix \right |
\left |  \matrix
x_1&\boxed {x_2}\\
1                &                    1
\endmatrix \right |^{-1}=(x_2^{k+1}-x_1^{k+1})(x_2-x_1)^{-1}
\\
&=y_2^k+y_2^{k-1}y_1+\dots +y_2y_1^{k-1}+y_1^k,
\endalign
$$
where $y_1=x_1$, $y_2=(x_2-x_1)x_2(x_2-x_1)^{-1}$;
$$
\align
\Lambda_1(x_1, x_2)&=
\left |  \matrix
x_1^2&\boxed {x_2^2}\\
1        &                    1
\endmatrix \right |
\left |  \matrix
x_1&\boxed {x_2}\\
1    &         1
\endmatrix \right |^{-1}=(x_2^2-x_1^2)(x_2-x_1)^{-1}=y_1+y_2,
\\
\Lambda_2(x_1, x_2)&=
-\left |  \matrix
x_1^2&\boxed {x_2^2}\\
x_1        &          x_2
\endmatrix \right |
\left |  \matrix
x_1&x_2\\
1    & \boxed  {1}
\endmatrix \right |^{-1}=(x_2^2-x_1x_2)(1-x_1^{-1}x_2)^{-1}=y_2y_1.
\endalign
$$
\endexample
\remark {Remark} Formulas for $S_k$ in Proposition 6.7.1 are valid for all
$k\in \Bbb Z$.
\endremark

An important "periodicity" property of quasi-Pl\"ucker coordinates of
Vandermonde matrices
is given by the following proposition.

\proclaim {Proposition 6.7.2} For any $k\in \Bbb Z$ we have
$$
r_{ij}^{i_1,\dots ,i_{n-1}}(V_n)=r_{i+k,j+k}^{i_1+k,\dots,i_{n-1}+k}(V_n).
$$
\endproclaim

Proposition 6.7.1 can be generalized as follows. Recall that
in 6.4 we defined generalized quasi-Schur functions. Let
$I=\{i_1,i_2,\dots , i_m\}$.

\proclaim{Proposition 6.7.3} Let   $I=\{i_1,i_2,\dots , i_m\}$ and
$i_1\geq i_2\geq \dots \geq  i_m$. Set $J=\{0,1,\dots ,
n-m-1,n-m+i_1,\dots , n-2+i_{m-1}\}$. Then
$$
r_{n-1+i_1, n-m}^J(V_n)=\check S_I.
$$
\endproclaim

%\example{Examples} Special cases of Proposition 6.7.3 are
%$$
%\gather
%r_{n+k-1, n-1}^{01\dots n-2}=S_k,
%\\
%r_{n-1, n-k}^{01\dots n-k-1, n-k+1,\dots , n-2}=\Lambda _k.
%\endgather
%$$
%\endexample

\head 7. Universal quadratic algebras associated with
pseudo-roots of noncommutative polynomials and noncommutative
differential polynomials
\endhead

\subhead {7.1. Pseudo-roots of noncommutative polynomials}\endsubhead
During the last years the
authors introduced and studied
universal algebras associated with pseudo-roots of noncommutative
polynomials. The results appeared in \cite{GRW, GGR,
GGRWS}.

Let $C$ be an algebra with unit and $P(t)\in C[t]$ be a
polynomial (where $t$ is a formal variable commuting with
elements of $C$). We say that an element $c\in C$ is a {\it
pseudo-root} of $P(t)$ if there exist polynomials $L_c(t),
R_c(t)\in C[t]$ such that $P(t)=L_c(t)(t-c)R_c(t)$. If
$P(t)=a_0t^n+a_1t^{n-1}+\dots +a_{n-1}t+a_n$ and $c$ is a
pseudo-root of $P(t)$ with $R_c(t)=1$, then
$$
a_0c^n+a_1c^{n-1}+\dots +a_{n-1}c+a_n=0,
$$
i.e., $c$ is a root of the polynomial
$P(x)=a_0x^n+a_1x^{n-1}+\dots +a_{n-1}x+a_n$ (where $x$ is a
noncommuting variable). Our theory shows that an
the analysis of noncommutative polynomials is impossible without
studying pseudo-roots.

Let $x_1,\dots , x_n$ be roots of a {\it generic} monic
polynomial $P(x)=x^n+a_1x^{n-1}+\dots +a_n$ over an algebra $C$.
There are two important classical problems: a) to express the
coefficients $a_1,\dots , a_n$ via the roots, b) to determine
all factorizations of $P(x)$, or $P(t)$.

When $C$ is a division ring, the first problem was solved in
\cite {GR3, GR4} using the theory of quasideterminants; the solution is 
presented in Section 5. Let
$V(x_{i_1},\dots , x_{i_k})$ be the Vandermonde quasideterminant
corresponding to the sequence $x_{i_1},\dots , x_{i_k}$. 
For an ordering $\{i_1,\dots , i_n\}$ of $\{1,\dots ,n\}$, 
we constructed the elements 
$$
\align
x_{\emptyset , i_1}&=x_{i_1},\\
x_{\{i_1, i_2, \dots ,i_{k-1}\}, i_k}&= V(x_{i_1},\dots ,
x_{i_k})x_{i_k}V(x_{i_1},\dots , x_{i_k})^{-1},\quad  k=2,\dots , n,
\endalign
$$
in $C$ such that for every $m=1, \dots, n$,
$$
(-1)^ma_m=\sum _{j_1>j_2>\dots j_m}y_{j_1}y_{j_2} \dots y_{j_m},
\tag 7.1.1
$$
where $y_1=x_{i_1}$, $y_k=x_{\{i_1,\dots , i_{k-1}\}, i_k}$,
$k=2,\dots , n$.

It is surprising that the left-hand side in formula (7.1.1) does
not depend on the ordering of $\{1,\dots , n\}$ whereas the
right-hand side {\it a priori\/} does depend on it. The
independence of the right-hand side in (7.1.1) on the ordering of
$\{1,\dots , n\}$ was a key point in the theory of noncommutative
symmetric functions developed in \cite {GR3, GR4}, see
Section 6 of the present article.

The element  $x_{\{i_1,\dots , i_{k-1}\}, i_k}$ has an
interesting structure. As we have already mentioned, it is symmetric in $x_{i_1}, \dots ,x_{i_{k-1}}$. Next, it is a rational function in $x_{i_1}, \dots ,x_{i_k}$ containing, in the general case, $k-1$ inversions. In
other words, $x_{\{i_1,\dots , i_{k-1}\}, i_k}$ is a rational
expression of height $k-1$.

Set $A_k=\{i_1, \dots , i_{k-1}\}$ for $k=2,\dots , n$. Recall that
formulas (7.1.1) are equivalent to the decomposition
$$
P(t)=(t-x_{A_n, i_n})(t-x_{A_{n-1}, i_{n-1}})\dots (t-x_{i_1}).
\tag 7.1.2
$$
This formula shows that the elements $x_{A,i}$ are pseudo-roots of
$P(t)$, and in the general case the polynomial $P(t)$ has at
least $n2^{n-1}$ different pseudo-roots. We study all these pseudo-roots
together by constructing the universal algebra of pseudo-roots
$Q_n$.

\subhead {7.2. Universal algebra of pseudo-roots}\endsubhead
 It is easy to see that the elements $x_{A,i}$, $i\notin
A$, satisfy the following simple relations:
$$
\align
x_{A\cup \{i\},j}+ x_{A, i}&=x_{A\cup \{j\},i}+ x_{A, j}, \tag 7.2.1a
\\
x_{A\cup \{i\},j}\cdot x_{A, i}&=x_{A\cup \{j\},i}\cdot x_{A, j},
\tag 7.2.1b
\endalign
$$
for all $A\subseteq \{1,\dots , n\}$, $i,j\notin A$.

In order to avoid inversions and to make our construction
independent of the algebra $C$, we define the universal algebra
of pseudo-roots $Q_n$ over a field $F$ to be the algebra with generators $z_{A,
i}$, $A\subseteq \{1,\dots , n\}$, $i\notin A$, and relations
corresponding to (7.2.1) (with $x$ replaced by $z$).

Each algebra $Q_n$ has a
natural derivation $\partial : Q_n\rightarrow Q_n$,
$\partial (z_{A,i})=1$ and a natural 
anti-involution $\theta : Q_n\rightarrow Q_n$,
$\theta (z_{A,i})=z_{C,i}$ where $C=\{1,\dots , n\}\setminus A\setminus
\{i\}$.

The algebra contains, as a subalgebra, the free associative algebra 
generated by the $z_i=z_{\emptyset, i}$, $i=1,\dots , n$. The
algebra $Q_n$ also admits a natural homomorphism $\alpha _n$ to
the skew-field generated by elements $z_1,\dots , z_n$. Namely, let
$A=\{i_1,\dots , i_k\}$. Set
$$
\alpha _n(z_{A,i})=\left |\matrix
z_{i_1}^k&\dots &z_{i_k}^k&\boxed {z_i^k}\\
z_{i_1}^{k-1}&\dots &z_{i_k}^{k-1}&z_i^{k-1}\\
             &\dots &             & \\
1            &\dots    & 1        &1 \endmatrix \right | z_i
\left |\matrix
z_{i_1}^k&\dots &z_{i_k}^k&\boxed {z_i^k}\\
z_{i_1}^{k-1}&\dots &z_{i_k}^{k-1}&z_i^{k-1}\\
             &\dots &             & \\
1            &\dots    & 1        &1 \endmatrix \right |^{-1}.
$$
\proclaim {Conjecture 7.2.1} The homomorphism $\alpha _n$ is a
monomorphism.
\endproclaim

We can prove this conjecture for $n=2,3$.

The algebra $Q_n$ admits a natural commutative specialization
$\pi _n :Q_n \to \allowmathbreak F[t_1,\dots , t_n]$ given by
$$ 
\pi_n(z_{A,i})=t_i,\quad  i=1,\dots ,n. \tag 7.2.2
$$
In particular, for each $i\notin A$, the image of $z_{A,i}-z_i$ under $\pi_n$ equals zero. 

To study the kernel of $\pi_n$ further, it is convenient to 
define, for each pair
$A, B\subseteq \{1,\dots , n\}$, with $A\cap B=\emptyset $ the
element $z_{A, B}\in Q_n$ by the recurrence formula 
$$
\gather
z_{\emptyset,\emptyset}=0,\\
z_{A\cup\{i\}, B}-z_{A, B\cup \{i\}}=z_{A,B}\quad\text{for $i\notin A,B$.}
\endgather
$$
One can easily see that when $B$ contains more than one element, the element $z_{A,B}$ is ``invisible" in the commutative case, i.e., $\pi_n(z_{A,B})=0$.

In \cite {GRW} it was proved that
$$
z_{A,i}-z_i=\sum_{\emptyset \neq C\subseteq A}z_{\emptyset , C\cup\{i\}}.
$$
The terms in the right-hand side in (7.2.2) measure the
``noncommutativity" of $z_{A,i}$. Moreover, in a sense, the ``degree of
noncommutativity'' carried by $z_{\emptyset ,B}$ depends on the
size of $B$: the greater $|B|$, the more ``noncommutative" the element
by $z_{\emptyset ,B}$ is.

\subhead 7.3. Bases in the algebra $Q_n$ \endsubhead
The algebra $Q_n$ has a natural graded structure $Q_n = \sum_{l
\ge 0} Q_{n,l}$ where $Q_{n,l}$ is the span of all products of
$l$ generators $z_{A,i}$.

One can see  that elements $z_{A, \emptyset}$ and, similarly, the
elements $z_{\emptyset, A}$, for all $A\subseteq \{1,\dots ,
n\}$, $A\neq \emptyset $, constitute a basis in the subspace of
$Q_{n,1}$. These elements satisfy simple
quadratic relations.

Our study of $Q_n$ relies on the construction of a basis in
$Q_n$, which is a hard combinatorial problem.

For $A \subseteq \{1,\dots ,n\}$ let $\min(A)$ denote the smallest
element of $A$.  Then set $A^{(0)} = A$, $A^{(1)} = A \backslash
\{\min(A)\}$, $A^{(i+1)} = (A^{(i)})^{(1)}$. Set
$r_A=z_{\emptyset, A}$.

\proclaim{Theorem 7.3.1 \rm(see \cite{GRW})}  The set of all
monomials
$$
r_{A_1^{(0)}}r_{A_1^{(1)}}\dots
r_{A_1^{(j_1)}}r_{A_2^{(0)}}r_{A_2^{(1)}}\dots
r_{A_2^{(j_2)}}\dots r_{A_l^{(0)}}r_{A_l^{(1)}}\dots
r_{A_l^{(j_l)}}
$$
where $A_1,\dots ,A_l \subseteq \{1,\dots ,n\}$ and for each $1 \le
i \le l - 1$, either $A_{i+1} \not\subseteq A_i$ or $|A_{i+1}|
\ne |A_i| - j_i - 1$, is a basis in  $Q_n$.
\endproclaim

\remark{Remark} It would be interesting to study in details combinatorial
properties of the basis and to give constructions of similar
bases in $Q_n$ using different techniques including
noncommutative Gr\"obner bases and Bergman's Diamond lemma.
\endremark

\subhead 7.4. Algebra of noncommutative symmetric polynomials as
a subalgebra in $Q_n$\endsubhead
For each ordering $I=(i_1,\dots , i_n)$ of $\{1,\dots ,
n\}$, there is a natural free subalgebra $Q_{n, I }\subset Q_n$
generated by $\{z_{\{i_1, \dots , i_k\}, \emptyset } \ |\ 1\leq
k\leq n\}$. Using the basis theorem in \cite{GRW}, we can describe
arbitrary intersections of subalgebras $Q_{n, I}$. In particular, we can prove
the following theorem. Let $\bold S_n$ be a subalgebra in
$Q_n$ generated by all coefficients $a_m$ (see formula (7.1.1.)).
One can identify $\bold S_n$ with the algebra
$\overline {\bSym}_n$ of noncommutative symmetric
functions in $x_1,\dots , x_n$.

\proclaim{Theorem 7.4.1 \rm(see \cite{GRW})}
The intersection of all subalgebras $Q_{n, I}$ coincides with
algebra $\bold S_n$.
\endproclaim

This is a purely noncommutative phenomenon: under the commutative
specialization $\pi_n$, all algebras $Q_{n,I}$ map to the algebra
of all polynomials and algebra $\bold S_n$ maps to the algebra of
symmetric funtions.

\subhead 7.5. The dual algebra $Q_n^!$ \endsubhead
The definition of the dual quadratic algebra and of Koszul
quadratic algebras can be found, e.g., in \cite {L\"o}.

Recall, that the quadratic algebra $Q_n$ has a natural graded structure $Q_n = \sum_{i
\ge 0} Q_{n,i}$ where $Q_{n,i}$ is the span of all products of
$i$ generators.  As usual, we denote the Hilbert series of $Q_n$
by $H(Q_n,\tau ) = \sum_{i \ge 0} \dim(Q_{n,i})\tau ^i$.

In \cite {GGRWS} we computed the Hilbert series of $Q_n$ and of
its dual quadratic algebra $Q_n^!$. 
In particular, the following result was proved.

\proclaim{Theorem 7.5.1} We have
$$
\align
H(Q_n,\tau ) &= \frac{1-\tau }  {1-\tau (2-\tau )^n},\\
H(Q_n^!,\tau ) &= \frac{1+\tau (2+\tau )^n} {1+\tau }.
\endalign
 $$
\endproclaim

In particular, since $H(Q_n^!,\tau )$ is a polynomial in $\tau $, the dual 
algebra $Q_n^!$ is finite dimensional.
Similarly to $Q_n$, it also has a rich and interesting structure. 

Theorem 7.5.1 shows that 
$$
H(Q_n,\tau )\cdot H(Q_n^!, -\tau )=1.
$$ 
This also follows from the koszulity of $Q_n$, which was recently proved by Serconek and
Wilson \cite {SW}.

\subhead 7.6. Quotient algebras of $Q_n$\endsubhead
 There are at least two reasons to study quotient algebras of
$Q_n$. The noncommutative nature of $Q_n$ can be studied by
looking at quotients of $Q_n$ by ideals generated by some
$z_{\emptyset, A}$. These quotients are ``more commutative'' than
$Q_n$. For example, the quotient of $Q_n$ by the ideal generated
by all $z_{\emptyset , A}$ with $|A|\geq 2$ is isomorphic to the
algebra of commutative polynomials in $n$ variables.

To consider more refined cases, we need to turn to a
``noncommutative combinatorial topology''.  In this approach the
algebra $Q_n$ corresponds to an $n$-simplex $\Delta _n$, and we
consider quotients of $Q_n$ by ideals generated by some
$z_{\emptyset, A}$ corresponding to subcomplexes of $\Delta _n$.
In \cite {GGR} we described generators and relations for those
quotients. Special attention was given to the quotients of $Q_n$
corresponding to 1-dimensional subcomplexes of $\Delta _n$ (they
are close to algebras of commutative polynomials). However, we
need to study other quotient algebras of this type and to be able
to ``glue" together such quotient algebras.
This will lead to a construction of a ``noncommutative
combinatorial topology''.

Another interesting class of quotients of $Q_n$ consists of algebras
corresponding to special types of polynomials (such as
$x^n=0$ or polynomials with multiple roots). Here is an example.

\example{Example} Let $F$ be a field. We consider quotients of the $F$-algebra
$Q_2$. The algebra $Q_2$ itself is generated by $z_1$,
$z_2$, $z_{1,2}$, $z_{2,1}$. It corresponds to a polynomial
$P(t)=t^2-pt+q$ with $p=z_1+z_{1,2}=z_2+z_{2,1}$,
$q=z_{1,2}z_1=z_{2,1}z_2$.

There is an ``invisible element" $z_{\emptyset, 12}=
z_{1,2}-z_2=z_{2,1}-z_1$. This element satisfies the relation
$z_{\emptyset, 12}(z_1-z_2)=z_1z_2-z_2z_1$.

There are three quotient algebras corresponding to special cases
of $P(t)$.

(i) The algebra corresponding to the polynomial $t^2$ is 
$Q_2/(p,q)$. This algebra is isomorphic to $F\langle z_1,z_2\rangle/(z_1^2,z_2^2)$.

(ii) The algebra corresponding to a polynomial with multiple roots
is $Q_2/(z_1-z_2)$. It is isomorphic to the free algebra with
generators $z_1$ and $z_{1,2}$.

(iii) The algebra $Q_2/(z_{\emptyset, 12})$ is isomorphic to the
algebra of commutative polynomials in two variables.
\endexample

Another class of quotient algebras of $Q_n$ was introduced in \cite{GGR}.

\subhead 7.7. Quadratic algebras associated with differential
noncommutative polynomials \endsubhead
In \cite {GRW} we also constructed, similarly to algebras $Q_n$,
the universal algebras of pseudo-roots of  noncommutative
differential polynomials.

Let $C$ be an algebra over a field $F$. Recall that
a derivation $D$ of $C$ is an $F$-linear map
$D: C\to C$ such that $D(ab)=D(a)b+aD(b)$.
Any element $a\in C$ acts on $C$ by the multiplication
from the left.  It is obvious that the commutator $[D, a]=Da-aD$ acts
on $C$ as the left multiplication by $D(a)$.
Also, any polynomial $P=P(D)=a_0D^n+a_1D^{n-1}+\dots +a_n$, $a_i\in C$ 
for $i=0,1,\dots,n$ acts on $C$ by the formula
$$
P(D)(\phi)=a_0D^n(\phi)+a_1D^{n-1}(\phi)+\dots +a_{n-1}D(\phi)+ a_n\phi.
$$

We say that an element $c\in C$ is a pseudo-root of $P(D)$ if
there exist polynomials $L_c(D)$ and $R_c(D)$ with coefficients
in $C$ such that $P(D)=L_c(D)(D-c)R_c(D)$ (taking into account the commutation rule (7.7.1)). We say that $c$ is a
root of $P(D)$ if $R_c(D)=1$. 

Suppose that $a_0=1$ and the differential polynomial $P(D)$ has
$n$ different roots $f_1,\dots , f_n\in C$. Following \cite {EGR}, for
any ordering $(i_1,\dots , i_n)$ of $\{1,\dots,n\}$, in \cite{GRW} we constructed, for 
a generic $P$, pseudo-roots $f_{i_1,i_2}$,
\dots , $f_{i_1,\dots i_{n-1},i_n}$ such that
$$
P(D)=(D-f_{i_1,\dots ,i_{n-1},i_n})\dots (D-f_{i_1,i_2})
(D-f_{i_1}).
$$
For $k=2,\dots , n$ the element $f_{i_1,\dots ,i_{k-1},i_k}$ does
not depend on the order of elements $(i_1,\dots ,i_{k-1})$.

Set $f_{\emptyset, i}=f_i$. It was proved in \cite {GRW} that for
any $A\subset \{1,\dots,n\}$ such that  $|A|<n-1$, and for any
$i,j\notin A$ we have
$$
\align
f_{A\cup i,j}+f_{A,i}&=f_{A\cup j,i}+f_{A,j}, \tag 7.7.1a
\\
f_{A\cup i,j}f_{A,i}-D(f_{A,i}) &= f_{A\cup
j,i}f_{A,j}-D(f_{A,j}) . 
\tag 7.7.1b
\endalign
$$
Based on these formulas one can define universal algebras $DQ_n$
of pseudo-roots of  noncommutative differential polynomials. They
are defined by elements $f_{A,i}$ for $i\notin A$ satisfying
relations (7.7.1). The theory of algebras $DQ_n$ seems to be 
useful in the study of noncommutative integrable systems.

\head 8. Noncommutative traces, determinants and
eigenvalues\endhead

In this section we discuss noncommutative traces, determinants
and eigenvalues. Our approach to noncommutative determinants in
this Section is different from our approach described in Section
3. 

Classical (commutative) determinants play a key role in
representation theory. Frobenius developed his theory of group
characters by studying factorizations of group determinants (see
\cite {L}). Therefore, one cannot start a noncommutative
representation theory without looking at possible definition of
noncommutative determinants and traces. The definition of a
noncommutative determinant given in this Section is different
from the definition given in Section 3. However, for matrices 
over commutative algebras, quantum and Capelli matrices both 
approach give the same results.

\subhead 8.1. Determinants and cyclic vectors\endsubhead
Let $R$ be an algebra with unit and 
$A:R^m\to R^m$ a linear map of right vector spaces,
A vector $v\in R^m$ is an $A$-{\it cyclic vector\/} 
if $v, Av,\dots , A^{m-1}v$ is a basis in $R^m$ regarded as a right
$R$-module. In this case there exist $\Lambda _i(v,A)\in R$,
$i=1,\dots ,m$, such that
$$
(-1)^mv\Lambda _m(v,A) + (-1)^{m-1}(Av)\Lambda _{m-1}(v,A) +
\dots -(A^{m-1}v)\Lambda _1(v,A) +A^mv=0.
$$
\definition{Definition 8.1.1}
We call $\Lambda _m(v,A)$ the {\it determinant\/} of $(v,A)$ and
$\Lambda _1(v,A)$ the {\it trace\/} of $(v,A)$. 
\enddefinition

We may express
$\Lambda _i(v,A)\in R$,  $i=1,\dots ,m$, as quasi-Pl\"ucker
coordinates of the $m\times (m+1)$ matrix with columns
$v,Av,\dots,A^n v$ (following  \cite {GR4}).

In the basis $v,Av,\dots , A^{m-1}$ the map $A$ is 
represented by the Frobenius matrix
$A_v$ with the last column equal to
$((-1)^m\Lambda _m(v, A), \dots , -\Lambda _1(v,A))^T$.
>From Theorem 3.1.3 it follows that if determinants
of $A_v$ are defined, then they coincide up to a sign
with $\Lambda _m(V,A)$. This justifies our 
definition.

Also, when $R$ is a commutative algebra, $\Lambda _m(v,A)$ is the
determinant of $A$ and $\Lambda _1(v,A)$ is the
trace of $A$.

When $R$ is noncommutative,  the expressions
$\Lambda _i(v,A)\in R$, $ i=1,\dots ,m$, depend on vector $v$.
However, they provide some information about $A$. For example,
the following statement is true.

\proclaim{Proposition 8.1.2} If
the determinant $\Lambda _m(v,A)$ equals zero, then the map $A$ is
not invertible.
\endproclaim

Definition 8.1.1 of noncommutative determinants and traces was
essentially used in \cite {GKLLRT} for linear maps given by
matrices $A=(a_{ij})$, $i,j=1,\dots , m$ and unit vectors 
$e_s$, $s=1,\dots , m$. In this
case $\Lambda _i(e_s, A)$ are quasi-Pl\"ucker coordinates of the
corresponding Krylov matrix $K_s(A)$. Here (see \cite {G})
$K_s(A)$ is the matrix $(b_{ij})$, $i=m,m-1,\dots , 1,0$,
$j=1,\dots , m$, where $b_{ij}$ is the $(sj)$-entry of $A^i$.

\example{Example} Let $A=(a_{ij})$ be an $m\times m$-matrix and
$v=e_1=(1,0,\dots ,0)^T$. Denote by $a_{ij}^{(k)}$ the corresponding
entries of $A^k$. Then
$$
\Lambda _m(v, A)=(-1)^{m-1}\left |\matrix
\boxed {a_{11}^{(m)}}&a_{12}^{(m)}&\dots &a_{1m}^{(m)}\\
a_{11}^{(m-1)}&a_{12}^{(m-1)}&\dots &a_{1m}^{(m-1)}\\
              &\dots         &\dots &              \\
a_{11}&a_{12}&\dots &a_{1m}\endmatrix \right |.
$$
For $m=2$ the ``noncommutative trace"
$\Lambda _1$ equals $a_{11}+a_{12}a_{22}a_{12}^{-1}$ and
the ``noncommutative determinant"
$\Lambda _2$ equals $a_{12}a_{22}a_{12}^{-1}a_{11}-a_{12}a_{21}$.
\endexample

It was shown in \cite{GKLLRT} that if $A$ is a quantum matrix, then
$\Lambda _m$ equals $\det_q A$ and $A$ is a Capelli
matrix, then $\Lambda _m$ equals the Capelli determinant.

A construction of a noncommutative determinant and a
noncommutative trace in terms of cyclic vectors in a special case was 
used in \cite {Ki}.

One can view the elements $\Lambda _i(v,A)$ as elementary
symmetric functions of ``eigenvalues" of $A$.

Following Section 6 we introduce complete
symmetric functions $S_i(v,A)$, $i=1,2,\dots $,
of ``eigenvalues" of $A$ as follows. Let $t$ be a formal
commutative variable. Set $\lambda (t)=1+\Lambda_1(v,A)t
+\dots +\Lambda_m(v,A)t^m$ and define the elements $S_i(v,A)$
by the formulas
$$ 
\sigma(t):=1+\sum _{k>0}S_kt^k=\lambda(-t)^{-1}.
$$
Recall that in Section 6 we introduced ribbon Schur
functions and that $R_{1^kl}$ is the ribbon Schur
function corresponding to the hook with $k$
vertical and $l$ horizontal boxes.
In particular, $\Lambda_k=R_{1^k}$,
$S_l=R_l$.

Let $A:R^m\rightarrow R^m$ be a linear map of right
linear spaces.

\proclaim {Proposition 8.1.3}
For $k\geq 0$
$$
A^{m+k}v=(-1)^{m-1}vR_{1^{m-1}(k+1)} + (-1)^{m-2}(Av(R_{1^{m-2}(k+1)}
+\dots + (A^{m-1}v)R_{k+1}.
$$
\endproclaim

Let $A=\text{diag}(x_1,\dots , x_m)$. In the general
case for a cyclic vector one can take $v=(1,\dots ,1)^T$.
In this case, the following two results hold.

\proclaim{Proposition 8.1.4} For $k=1,\dots , m$
$$
\multline
\Lambda _k(v,A)\\
=
\left |\matrix
1&\dots&\boxed{x_m^{m-k}}&\dots&x_m^{m-1}\\
 &     &\dots            &     &         \\
1&\dots&x_1              &\dots&x_1^{m-1}
\endmatrix\right |^{-1}\cdot
\left |\matrix
1&\dots&x_m^{m-k-1}&x_m^{m-k+1}&\dots&\boxed{x_m^m}\\
 &     &\dots      &           &     &         \\
1&\dots&x_1^{m-k-1}&x_1^{m-k+1}&\dots&x_1^m
\endmatrix\right |.
\endmultline
$$
\endproclaim

\proclaim{Proposition 8.1.5} For any $k>0$
$$
S_k(v,A)=
\left |\matrix
1&\dots&\boxed{x_m^{m-1}}\\
 &\dots&             \\
1&\dots&x_1^{m-1}    
\endmatrix\right |^{-1}\cdot
\left |\matrix
1&\dots&x_m^{m-2}&\boxed{x_m^{m+k-1}}\\
 &     &\dots    &                   \\
1&\dots&x_1^{m-2}&x_1^{m+k-1}
\endmatrix\right |.
$$
\endproclaim

Note that formulas for $S_k$ look somewhat simpler
than formulas for $\Lambda_k$.

\subhead 8.2. Noncommutative determinants and noncommutative eigenvalues
\endsubhead
One can also express $\Lambda _i(v,A)\in R$ in terms of left eigenvalues
of $A$.

Let a linear map $A:R^m\to R^m$ of the right vector
spaces is represented by the matrix $(a_{ij})$.

\definition{Definition 8.2.1} A nonzero row-vector $u=(u_1,\dots , u_m)$ is
a left eigenvector of $A$ if there exists $\lambda \in R$ such
that $uA=\lambda u$.
\enddefinition

We call $\lambda $ a {\it left eigenvalue} of $A$ corresponding
to vector $u$. Note, that $\lambda$ is the eigenvalue of $A$ corresponding to a left eigenvector $u$ then, for each $\alpha\in R$, $\alpha\lambda\alpha^{-1}$ is the eigenvalue corresponding to the left eigenvector $\alpha u$. Indeed,
$(\alpha u)A=\alpha\lambda \alpha ^{-1} (\alpha u)$.

For a row vector $u=(u_1,\dots , u_m)$ and a column vector $v=(v_1,\dots , v_m)^T$ denote by $\langle u,v\rangle$ the inner product $\langle u,v\rangle=u_1v_1+\dots u_mv_m$.

\proclaim{Proposition 8.2.2} Suppose that $u=(u_1,\dots , u_m)$ is
a left eigenvector of $A$ with the eigenvalue $\lambda$, $v=(v_1,\dots , v_m)^T$ is a cyclic vector of $A$, and $\langle u,v\rangle=1$. Then
The eigenvalue $\lambda$ satisfies the equation
$$
(-1)^m\Lambda _m(v,A) + (-1)^{m-1} \lambda \Lambda _{m-1}(v,A) +
\dots -\lambda ^{m-1}\Lambda _1(v,A) + \lambda ^m=0. \tag 8.2.1
$$
\endproclaim

Equation (8.2.1) and the corresponding Vi\`ete theorem (see Section 3)
show that if the map $A:R^m\to R^m$ has left eigenvectors
$u^1, \dots ,u^m$ with corresponding eigenvalues $\lambda _1,\dots , \lambda _m$
such that $\langle u^i, v\rangle=1$ for $i=1,\dots , m$ and  any submatrix
of the Vandermonde matrix $(\lambda _i^j)$ is invertible, 
then  all $\Lambda_i(v,A)$ can be expressed in terms of 
 $\lambda _1,\dots , \lambda _m$
as ``noncommutative
elementary symmetric functions" by formulas similar to
those in Definition 6.5.1.

\subhead 8.3. Multiplicativity of determinants\endsubhead
In the commutative case the multiplicativity of determinants and
the additivity of traces are related to computations of
determinants and traces with diagonal block-matrices. In the
noncommutative case we suggest to consider the following
construction.

Let $R$ be an algebra with a unit.
Let $A: R^m\rightarrow R^m$ and $D: R^n\rightarrow R^n$
be linear maps of right vector spaces, $v\in R^m$ 
an $A$-cyclic vector and $w\in R^n$ a $D$-cyclic vector.

There exist $\Lambda _i(w,D)\in R$, $i=1,\dots ,n$,
such that
$$
(-1)^nv\Lambda _n(w,D) + (-1)^{n-1}(Dw)\Lambda _{n-1}(w,D) + 
\dots -(D^{m-1}v)\Lambda _1(w,D) +D^nw=0.
$$
Denote also by $S_i(w,D)$, $i=1,2,\dots $, the corresponding
complete symmetric functions.

The matrix $C=\left(\matrix A&0\\0&D\endmatrix\right)$
acts on $R^{m+n}$. Suppose that the vector
$u=\left(\matrix v\\w\endmatrix\right)$ is a cyclic vector
for matrix $C$. We want to express $\Lambda_i(u,C)$,
$i=1,\dots , m+n$ in terms of $\Lambda_j(v,A)$, 
$S_k(v,A)$, $\Lambda_p(w,D)$, and $S_q(w,D)$.

Denote, for brevity, $\Lambda_j(v,A)=\Lambda_j$, $S_k(v,A)=S_k$, 
$\Lambda_p(w,D)=\Lambda_p'$, $S_q(w,D)=S_q'$.

For two sets of variables $\alpha=\{a_1,a_2, \dots , \}$ and
$\beta=\{b_1,b_2,\dots , \}$ introduce the following $(m+n)\times
(m+n)$-matrix $M(m,n;\alpha , \beta)$:
$$
\left(\matrix 
1&a_1&a_2&\dots &\dots  &\dots &a_{m-1} &\dots &a_{m+n-1}\\
0&1  &a_1&a_2   &\dots  &\dots &a_{m-1} &\dots &a_{m+n-2}\\
 &   &   &      &       &\dots &        &      &\\
0&0  & 0 &\dots &\dots  &\dots &1       &\dots &a_m\\
1&b_1&b_2&\dots &b_{n-1}&\dots &\dots   &\dots &b_{m+n-1}\\
0&1  &b_1&b_2  &\dots   &\dots &\dots   &\dots &b_{m+n-2}\\
 &   &   &      &       &\dots &        &      &\\
0&0  &0  &\dots &1      &b_1   &\dots   &\dots &b_n\\
\endmatrix\right).
$$

\proclaim{Proposition 8.3.1} For any $j=2, \dots , m+n$ we have
$$
|M(m,n;\alpha, \beta)|_{1j}=-|M(m,n;\alpha, \beta)|_{m+1,j}.
$$
\endproclaim

The elements $S_i(u,C)$, $i=1,2,\dots $, can be
computed as follows.
Denote by $N_k(m,n;\alpha , \beta)$ the matrix obtained
from $M$ by replacing its last column
by the following column: 
$$
(a_{m+n+k-1}, a_{m+n+k-2},\dots , a_{n+k-1}, b_{m+n+k-1}, b_{m+n+k-2},
\dots , b_{m+k-1})^T.
$$ 

Set $\alpha =\{-S_1,S_2,\dots , (-1)^kS_k, \dots \}$,
$\alpha '=\{-S_1',S_2',\dots , (-1)^kS_k', \dots \}$.
 
\proclaim{Theorem 8.3.2} For $k=1,2,\dots $ we have
$$ 
S_k(u,C)=|M(m,n;\alpha , \alpha ')|_{1m+n}^{-1}
\cdot |N_k(m,n;\alpha , \alpha ')|_{1m+n}.
$$
\endproclaim

\example{Example} For $m=3$, $n=2$ and $k=1,2,\dots $. Then
$$
\multline
S_k(u,C)=\\
=(-1)^{k-1}\left |\matrix
1&-S_1&S_2&-S_3&S_4\\
0&1&-S_1&S_2&-S_3\\
0&0&1&-S_1&S_2\\
1&-S_1'&S_2'&-S_3'&S_4'\\
0&1&-S_1'&S_2'&-S_3'
\endmatrix\right|_{15}^{-1}\times
\left |\matrix
1&-S_1&S_2&-S_3&S_{4+k}\\
0&1&-S_1&S_2&-S_{3+k}\\
0&0&1&-S_1&S_{2+k}\\
1&-S_1'&S_2'&-S_3'&S_{4+k}'\\
0&1&-S_1'&S_2'&-S_{3+k}'
\endmatrix\right|_{15}.
\endmultline
$$
\endexample

For $n=1$ denote $\Lambda _1(D)=S_1(D)$ by $\lambda '$. 

\proclaim{Corollary 8.3.3} If $n=1$, then for $k=1,2,\dots $ we have
$$
S_k(u,C)=S_k(v,A) +S_{k-1}(v,A)|M(m,n;\alpha , \alpha ')|_{1m+n}^{-1}
\lambda '|M(m,n;\alpha , \alpha ')|_{1m+n+1}. 
$$
\endproclaim

Note that
$$
\Lambda_{m+1}(u,C)=|M(m,n;\alpha , \alpha ')|_{1m+n}^{-1}
\lambda '|M(m,n;\alpha , \alpha ')|_{1m+n+1}\Lambda _m(v, A),
$$
i.e. the ``determinant" of the diagonal matrix equals
the product of two ``determinants".

\head 9. Some applications \endhead

In this section we mainly present some results from \cite {GR1, GR2, GR4}.

\subhead 9.1. Continued fractions and almost triangular
matrices\endsubhead
Consider an infinite matrix $A$ over a skew-field:
$$
A=\pmatrix &a_{11}&a_{12}&a_{13}&\dots &a_{1n}\dots\\
           &-1    &a_{22}&a_{23}&\dots &a_{2n}\dots\\
           &0     &-1    &a_{33}&\dots &a_{3n}\dots\\
           &0     &0  &-1      &\dots &\dotso\endpmatrix
$$

It was pointed out in [GR1], [GR2] that the quasideterminant $|A|_{11}$
can be written as a generalized continued fraction
$$
|A|_{11} = a_{11} + \sum_{j_1\neq 1} a_{1j_1}{1\over
a_{2j_1}+\sum \Sb j_2\neq 1,j_1\endSb a_{2j_2} {1\over a_{3j_2}+\dots}}.
$$
Let
$$
A_n=\pmatrix &a_{11} & a_{12} &\dots &a_{1n}\\
&-1 &a_{22} &\dots &a_{2n}\\
&0 &-1 &\dots  &a_{3n}\\
&{} &{} &\dots &{}\\
&\dots  &0 &-1 &a_{nn}\\
\endpmatrix.
$$
The following proposition was formulated in [GR1], [GR2].

\proclaim{Proposition 9.1.1}
$|A_n|_{11} = P_n Q^{-1}_n$, where
$$
\align
P_n&=\sum_{1\leq j_1<\dots < j_k< n} a_{1j_1}a_{j_1+1,j_2}
a_{j_2+1,j_3}\dots a_{j_k+1,n},\tag 9.1.1
\\
Q_n&=\sum_{2\leq j_1<\dots < j_k< n}
a_{2j_1}a_{j_1+1,j_2}a_{j_2+1,j_3}\dots a_{j_k+1,n}.\tag 9.1.2
\endalign
$$
\endproclaim
\demo{Proof} From the homological relations one has
$$
|A_n|_{11}|A^{1n}_n|^{-1}_{21}= - |A_n|_{1n}|A^{11}_n|^{-1}_{2n}.
$$
We will apply formula (1.2.2) to compute $|A_n|_{1n}$,
$|A^{11}_n|_{2n}$, and $|A^{1n}_n|_{21}$. It is easy to see that
$|A^{1n}_n|_{21}=-1$. To compute the two other quasideterminants, we
have to invert triangular matrices. Setting $P_n=|A_n|_{1n}$ and
$Q_n=|A^{11}_n|_{2n}$ we arrive at formulas (9.1.1), (9.1.2).
\qed\enddemo

\remark{Remark}In the commutative case Proposition 9.1.1 is well known.
In this case $P_n=|A_n|_{1n}=(-1)^n\det A_n$
and $Q_n=(-1)^{n-1}\det A^{11}_n$.
\endremark

Formulas (9.1.1), (9.1.2) imply the following result (see \cite{GR1, GR2}).
\proclaim {Corollary 9.1.2}
The polynomials $P_k$ for $k\geq 0$ and $Q_k$ for $k\geq 1$ are related
by the formulas
$$
\align
P_k&=\sum^{k-1}_{s=0} P_s a_{s+1,k},\qquad P_0 = 1,\tag 9.1.3
\\
Q_k&=\sum^{k-1}_{s=1} Q_s a_{s+1,k},\qquad Q_1 = 1.\tag 9.1.4
\endalign
$$
\endproclaim

\proclaim{Corollary 9.1.3}
Suppose that for any $i\neq j$ and
any $p,q$ the elements of the matrix $A$ satisfy the conditions
$$
\align
a_{ij}a_{pq}&=a_{pq}a_{ij}
\\
a_{jj}a_{ii}-a_{ii}a_{jj}&=a_{ij},\quad 1\leq i<j\leq n. 
\endalign
$$
Then
$$ 
P_n=|A_n|_{1n}=a_{nn}a_{n-1n-1}\dots a_{11}. \tag 9.1.5
$$
\endproclaim

The proof follows from (9.1.3).

\proclaim{Corollary 9.1.4 \rm(\cite{GR1, GR2})} For the Jacoby matrix
$$
A=\pmatrix &a_1 & 1 & 0 &\dots\\
&-1&a_2 & 1&{}\\
&0 &-1 & a_3 &\dots\endpmatrix
$$
we have
$$
|A|_{11} = a_1 + {1\over a_2+{1\over a_3 + \dots}},
$$
and
$$
\align
P_0 =1,\quad P_1=a_1,\quad & P_k = P_{k-1}a_k + P_{k-2},\ \text{for $k\geq 2;$}
\\
Q_1 = 1,\quad Q_2=a_2, \quad & Q_k=Q_{k+1}a_k + Q_{k-2},\ \text{for $k\geq 3;$}.
\endalign
$$
\endproclaim

In this case $P_k$ is a polynomial in $a_1,\dots ,a_k$ and
$Q_k$ is a polynomial in $a_2,\dots , a_k$.

\subhead 9.2. Continued fractions and formal series\endsubhead
In the notation of the previous subsection the infinite continued
fraction $|A|_{11}$ may be written as a ratio of formal series
in the  letters $a_{ij}$ and $a^{-1}_{ii}$.
Namely, set
$$
\align
P_\infty&=\sum\Sb 1\leq j_1
<j_2\dots<j_k<r-1\\ r=1,2,3,\dots\endSb a_{1j_1}a_{j_1+1j_2}\dots
a_{j_k+1r}a^{-1}_{rr}\cdot\dots\cdot a_{11}^{-1} 
\\
&=1+ a_{12}a^{-1}_{22}a^{-1}_{11} +
a_{13}a^{-1}_{33}a^{-1}_{22} a^{-1}_{11} +
a_{11}a_{23}a^{-1}_{33}a^{-1}_{22}a^{-1}_{11}+\dots ,
\endalign
$$
and
$$
\align
Q_\infty&=a^{-1}_{11}+\sum\Sb 2\leq j_1
<j_2\dots<j_k<r-1\\ r=2,3\dots\endSb a_{2j_1}a_{j_1+1j_2}\dots
a_{j_k+1r}a^{-1}_{rr}\cdot\dots\cdot a_{11}^{-1} 
\\
&=a^{-1}_{11} + a_{23}a^{-1}_{33}a^{-1}_{22}a^{-1}_{11} +
a_{24}a^{-1}_{44}a^{-1}_{33} a^{-1}_{22} a^{-1}_{11} + \dots .
\endalign
$$
Since each monomial appears in these sums at most once, these are
well-defined formal series.

The following theorem was proved in \cite{PPR}.
Another proof was given in \cite{GR4}.

\proclaim{Theorem 9.2.1} We have
$$
|A|_{11} = P_\infty\cdot Q^{-1}_\infty.
$$
\endproclaim

\demo{Proof} Set $b_{ij}=a_{ij}a_{jj}^{-1}$ and consider matrix
$B=(b_{ij})$,
$i,j=1, 2, 3, \dots $. According to a property of quasideterminants
$|A|_{11}=|B|_{11}a_{11}$.
Applying the noncommutative Sylvester theorem to $B$ with matrix
$(b_{ij}), i,j \geq 3$, as the pivot, we have
$$
|B|_{11}= 1+|B^{21}|_{12}|B^{11}|^{-1}_{22}a^{-1}_{11}.
$$
Therefore
$$
|A|_{11}= (a_{11}|B^{11}|_{22}a^{-1}_{11}
+ |B^{21}|_{12}a^{-1}_{11})
(|B^{11}|_{22}a^{-1}_{11})^{-1}.\tag 9.2.1
$$
By \cite{GKLLRT}, Proposition 2.4, the first
factor in (9.2.1) equals $P_{\infty}$, and the second
equals $Q_{\infty}^{-1}$. \qed\enddemo

\subhead 9.3. Noncommutative Rogers-Ramanujan continued fraction
\endsubhead
The following application of Theorem 9.2.1 to Rogers-Ramanujan continued fraction
was given in \cite{PPR}.
Consider a continued fraction with two formal variables $x$ and $y$:
$$
A(x,y) = {1\over 1+x{1\over 1+x{1\over 1+\dots}y}y}.
$$
It is easy to see that
$$
A(x,y)=\vmatrix &1 &x &{}&\cdot&{}&{}&{}\\
&-y&1&x&{}&\cdot &{}&0\\
&{}&-y&1&x&{}&\cdot&{}\\
&{}&{}&{}&1&{}&{}&\cdot\\
&{}&0&{}&\ddots&\ddots&{}&\cdot\endvmatrix^{-1}_{11} =\vmatrix
&1&x&0&{}&{}\\
&-1&y^{-1}&xy^{-1}&{}&{}\\
&0&-1&y^{-1} &xy^{-1} &{}\\
&{}&{}&-1 &y^{-1}&\ddots\endvmatrix_{11}
$$
Theorem 9.2.1 implies the following result.

\proclaim {Corollary 9.3.1}
$A(x,y) = P\cdot Q^{-1} $, where $Q=yPy^{-1}$ and
$$
P=1+\sum\Sb k\geq 1\\n_1,\dots,n_k\geq 1\endSb y^{-n_1}x
y^{-n_2}x\dots \, y^{-n_k} x y^{k+n_1+n_2 +\dots+n_k}.
$$
\endproclaim

Following \cite{PPR}, let us assume that $xy=qyx$, where $q$ commutes with
$x$ and $y$. Set $z=yx$.
Then Corollary 9.3.1 implies Rogers-Ramanujan continued fraction identity
$$
A(x,y)={1\over 1+{qz\over 1+{q^2 z\over 1+\dots}}}
={{1+\sum_{k\geq 1} {q^{k(k+1})\over
(1-q)\dots(1-q^k)}z^k}
\over 1+\sum_{k\geq 1}{q^{k^2}\over(1-q)\dots(1-q^k)}z^k}.
$$

\subhead 9.4. Quasideterminants and characterisric functions of graphs
\endsubhead
Let $A=(a_{ij})$,  $i,j=1,\dots , n$, where $a_{ij}$ are formal
noncommuting variables. Fix $p, q\in \{1,\dots , n\}$ and a set
$J\subset \{1,\dots , \hat p, \dots , n\}\times  \{1,\dots , \hat q, \dots , n\}$
such that $|J|=n-1$ and both projections of $J$ onto
$\{1,\dots , \hat p, \dots , n\}$ and   $\{1,\dots , \hat q, \dots , n\}$
are surjective. Introduce new variables $b_{kl}$, $k,l=1,\dots , n$,
by the formulas $b_{kl}=a_{kl}$ for $(l, k)\notin J$, $b_{kl}=a_{lk}^{-1}$ for
$(l, k)\in J$. 
Let $F_J$ be a ring of formal series in variables $b_{kl}$.

\proclaim{Proposition 9.4.1} The quasideterminant $|A|_{ij}$ is defined in
the ring $F_J$ and is given by the formula
$$
|A|_{ij}=b_{ij} -\sum (-1)^sb_{ii_1}b_{i_1i_2}\dots b_{i_sj}. \tag 9.4.1
$$
The sum is taken over all sequences $i_1,\dots , i_s$ such that
$i_k\neq i, j$ for $k=1,\dots ,s$.
\endproclaim

\proclaim{Proposition 9.4.2} The inverse to $|A|_{ij}$ is also defined in
the ring $F_J$ and is given by the following formula
$$
|A|_{ij}=b_{ij} -\sum (-1)^sb_{ii_1}b_{i_1i_2}\dots b_{i_sj}. \tag 9.4.2
$$
The sum is taken over all sequences $i_1,\dots , i_s$.
\endproclaim

All relations between quasideterminants, including the Sylvester identity,
can be deduced from formulas (9.4.1) and (9.4.2).

Formulas (9.4.1) and (9.4.2) can be interpreted in terms graph theory.
Let $\Gamma _n$ be a complete oriented graph with vertices
$1,\dots , n$ and edges $e_{kl}$, where $k, l=1,\dots , n$.
Introduce a bijective correspondence between edges of the
graph and elements $b_{kl}$ such that $e_{kl}\mapsto b_{kl}$.

Then there exist a bijective correspondence between
the monomials $b_{ii_1}b_{i_1i_2}\dots b_{i_sj}$ and the paths
from the vertex $i$ to the vertex $j$.

\subhead {9.5. Factorizations of differential operators and
noncommutative variation of constants}\endsubhead

Let $R$ be an algebra  with a derivation $D:R\rightarrow R$.
Denote $Dg$ by $g'$ and $D^kg$ by $g^{(k)}$.
Let $P(D)=D^n+a_1D^{n-1}+\dots +a_n$ be a differential operator
acting on $R$ and $\phi _i$, $i=1,\dots , n$, be solutions of
the homogeneous equation $P(D)\phi =0$, i.e.,
$P(D)\phi _i=0$ for all $i$.

For $k=1,\dots , n$ consider the Wronski matrix
$$W_k=\left (\matrix
\phi _1^{(k-1)}&\dots &\phi _k^{(k-1)}\\
 &\dots & \\
\phi _1&\dots &\phi _k\endmatrix \right ) $$ and
suppose that any square submatrix of $W_n$ is invertible.

Set $w_k=|W|_{1k}$ and $b_k=w_k'w_k^{-1}$, $k=1,\dots , n$.

\proclaim{Theorem 9.5.1 \rm \cite {EGR}}
$$P(D)=(D-b_n)(D-b_{n-1})\dots (D-b_1).$$
\endproclaim

\proclaim{Corollary 9.5.2} Operator $P(D)$ can be factorized as
$$P(D)=(w_n\cdot D\cdot w_n^{-1})(w_{n-1}\cdot D\cdot w_{n-1}^{-1})
\dots (w_1\cdot D\cdot w_1^{-1}).$$
\endproclaim
 
One can also construct solutions of the nonhomogeneous equation
$P(D)\psi =f$, $f\in R$, starting with solutions $\phi _1,\dots , \phi _n$
of the homogeneous equation. Suppose that any square submatrix of $W_n$
is invertible and that
there exist elements $u_j\in R$, $j=1,\dots , n$, such that
$$
u_j'=|W|_{1j}^{-1}f. \tag 9.5.1
$$
\proclaim {Theorem 9.5.3} The element $\psi=\sum _{j=1}^{j=n}\phi _ju_j$
satisfies the equation
$$
(D^n+a_1D^{n-1}+\dots +a_n)\psi =f.
$$
\endproclaim

In the case where $R$ is the algebra of complex valued functions $g(x)$,
$x\in \Bbb R$ the solution $\psi $ of the nonhomogeneous equation
is given by the classical formula
$$
\psi (x)=\sum _{j=1}^{j=n}\phi _j\int \frac {\det W_j}{\det W}dx \tag 9.5.2
$$
where matrix $W_j$ is obtained fron the Wronski matrix $W$ by
replacing the entries in the $j$-th column of $W$ by $f,0,\dots ,0$.
It is easy to see that formula (9.5.1) and Theorem 9.5.3 imply formula (9.5.2).

\subhead 9.6. Iterated Darboux transformations \endsubhead
Let $R$ be a differential algebra with a derivation $D:R\to R$
and $\phi \in R$ be an invertible element. Recall that we
denote $D(g)=g'$ and
$D^{k}(g)=g^{(k)}$. In particular $D^{(0)}(g)=g$.

For $f\in R$ define $\Cal D(\phi ;f)=f'-\phi '\phi ^{-1}f$. Following
\cite{Mat} we call $\Cal D(\phi ;f)$ the {\it Darboux transformation\/} of
$f$
defined by $\phi $. This definition was known for matrix functions $f(x)$
and $D=\partial _x$. Note that
$$
\Cal D(\phi ;f)=\left |\matrix \boxed {f'}&\phi '\\f&\phi \endmatrix \right|.
$$

Let $\phi _1,\dots , \phi _k$. Define the {\it iterated\/}
Darboux
transformation $\Cal D(\phi_k,\dots \phi _1;f)$ by induction as follows. 
For $k=1$, it coincides with the Darboux transformation defined above.
Assume that $k> 1$. The expression
$\Cal D(\phi_k,\dots ,\phi _1;f)$ is defined if $\Cal D(\phi_k,\dots ,\phi_2;f)$
is defined and invertible and $\Cal D(\phi _k;f)$ is defined. In this case,
$$
\Cal D(\phi_k,\dots \phi _1;f)=\Cal D(\Cal D(\phi_k,\dots \phi _2;f);
\Cal D(\phi _1;f).
$$

\proclaim{Theorem 9.6.1} If all square submatrices of matrix $(\phi _i^{(j)})$,
$i=1,\dots , k$; $j=k-1,\dots, 0$ are invertible, then
$$
\Cal D(\phi_k,\dots , \phi _1;f)=
\left |\matrix \boxed {f^{(k)}}&\phi _1^{(k)}&\dots &\phi _k^{(k)}\\
\dots &\dots &\dots &\dots \\
f&\phi _1&\dots &\phi _k\endmatrix \right |.
$$
\endproclaim

The proof follows from the noncommutative Sylvester theorem (Theorem 1.5.2).

\proclaim{Corollary 9.6.2} The iterated Darboux transformation $\Cal
D(\phi_k,\dots \phi _1;f)$
is symmetric in $\phi_1,\dots ,\phi _k$.
\endproclaim

The proof follows from the symmetricity of quasideterminants.

\proclaim{Corollary 9.6.3 \rm(\cite{Mat})} In commutative case, the iterated
Darboux transformation is a ratio of two Wronskians,
$$
\Cal D(\phi_k,\dots ,\phi _1;f)=\frac {W(\phi_1,\dots \phi
_k,f)}{W(\phi_1,\dots, \phi _k)}.
$$
\endproclaim

\subhead 9.7. Noncommutative Sylvester--Toda lattices \endsubhead
Let $R$ be a division ring with a derivation $D:R\rightarrow R$.
Let $\phi \in R$ and the quasideterminants
$$
T_n(\phi )=\left |\matrix
\phi &D\phi  &\dots &D^{n-1}\phi \\
D\phi &D^2\phi  &\dots &D^n\phi \\
\dots &\dots &\dots &\dots \\
D^{n-1}\phi &D^n\phi  &\dots &\boxed {D^{2n-2}\phi }
\endmatrix \right | \tag 9.7.1
$$
are defined and invertible. Set $\phi _1=\phi $ and $\phi _n=T_n(\phi )$,
$n=2,3,\dots $.

\proclaim{Theorem 9.7.1} Elements $\phi _n$, $n=1,2,\dots $,
satisfy the following  system of equations:
$$
\align
D((D\phi _1)\phi _1^{-1})&=\phi _2\phi _1^{-1},
\\
D((D\phi _n)\phi _n^{-1})&=\phi _{n+1}\phi _n^{-1}-
\phi _n\phi _{n-1}^{-1}, \ \ n\geq 2.
\endalign
$$
\endproclaim

If $R$ is commutative, the determinants of matrices used in formulas (9.7.1)
satisfy a nonlinear system of differential equations. In the modern
literature this system is called the Toda lattice (see, for example, \cite{Ok}
but in fact it was discovered by Sylvester in 1862 \cite {Syl} and, probably, 
should be called the Sylvester--Toda
lattice. Our system can be viewed as a
noncommutative generalization of the Sylvester--Toda lattice.
Theorem 9.7.1 appeared in \cite {GR1, GR2} and was generalized in \cite{RS}
and \cite {EGR}.

The following theorem is a noncommutative analog of the famous Hirota
identities.
\proclaim{Theorem 9.7.2} For $n\geq 2$
$$
T_{n+1}(\phi )=T_n(D^2\phi)-T_n(D\phi)\cdot ((T_{n-1}(D^2\phi )^{-1}-
T_n(\phi )^{-1})^{-1}\cdot T_n(D\phi ).
$$
\endproclaim
The proof follows from the Theorem 1.5.2.

\subhead 9.8. Noncommutative orthogonal polynomials \endsubhead
The results described in this subsection were obtained in \cite {GKLLRT}.
Let $S_0, S_1, S_2, \dots $ be elements of a skew-field
$R$ and $x$ be a commutative variable. Define a sequence of
elements $P_i(x)\in R[x]$, $i=0, 1, \dots $, by setting
$P_0=S_0$ and
$$
P_n(x)=\left |\matrix
S_n&\dots &S_{2n-1}&\boxed {x^n}\\
S_{n-1}&\dots &S_{2n-2}&x^{n-1}\\
\dots  &\dots &\dots   &\dots \\
S_0&\dots &S_{n-1}&1 \endmatrix \right |\tag 9.8.1
$$
for $n\geq 1$. We suppose here that quasideterminants
in (9.8.1) are defined. Proposition 1.5.1 implies that
$P_n(x)$ is a polynomial of degree $n$. If $R$ is
commutative, then $P_n, n\geq 0$, are orthogonal polynomials 
defined by the moments $S_n$, $n\geq 0$. We are going to show that if $R$ is a 
free division ring generated by $S_n$, $n\geq 0$, then polynomials
$P_n$ are indeed orthogonal with regard
a natural noncommutative $R$-valued product on $R[x]$.

Let $R$ be a free skew-field generated by $c_n$, $n\geq 0$.
Define on $R$ a natural anti-involution $a\mapsto a^*$ by
setting $c_n^*=c_n$ for all $n$. Extend the involution to
$R[x]$ by setting $(\sum a_ix^i)^*=\sum a_ix^i$. Define the
$R$-valued inner product on $R[x]$ by setting
$$
\Big\langle \sum a_ix^i, \sum b_jx^j\Big\rangle =\sum a_ic_{i+j}b_j^*.
$$

\proclaim{Theorem 9.8.1} For $n\neq m$ we have
$$
\langle P_n(x), P_m(x) \rangle =0.
$$
\endproclaim

The three term relation for noncommutative orhogonal polynomials
$P_n(x)$ can be expressed in terms of noncommutative quasi-Schur functions
$\check S_{i_1,\dots , i_N}$ defined in 6.4. We will use a notation
$\check S_{i^{N-1}j}$ if $i_1=\dots =i_{N-1}$ and $i_N=j$ and
write $\check S_{i^N}$ if $i_1=\dots =i_N$.

\proclaim{Theorem 9.8.2} The noncommutative orthogonal polynomials
$P_n(x)$ satisfy the three term recurrence relation
$$
P_{n+1}(x)-(x-\check S_{n^n(n+1)}^*\check S_{n^{n+1}}^{-1}
+\check S_{(n-1)^{n-1}n}^*\check S_{(n-1)^n}^{-1})P_n(x)
+\check S_{n^{n+1}}^*\check S_{(n-1)^n}^{-1}P_{n-1}(x)=0
$$
for $n\geq 1$.
\endproclaim


\Refs
\widestnumber\key{GKLLRT}

\ref\key A\by E. Artin  \book  Geometric algebra
\publ Wiley  \publaddr New York \yr 1988 
\endref

\ref\key Al \by S. Alesker \paper  Non-commutative linear algebra and
plurisubharmonic functions of quaternionic variables
\jour arXiv: math.CV/010429
\yr 2002
\endref

\ref\key As\by H. Aslaksen  \paper  Quaternionic determinants
\jour Math. Intelligencer \vol 18 \issue 3 \pages 57--65
\yr 1996
\endref


\ref\key B\by F. A. Berezin  \book  Introduction to algebra and analysis with
anticommuting variables 
\publ Izd. Moskov. Univ.  \publaddr Moscow \yr 1983 \lang Russian
\endref

\ref\key BW \by J. Birman and R. Williams \paper Knotted periodic
orbits I: Lorenz knots \jour Topology \vol 22 \pages 47--62
\yr 1983
\endref

\ref\key C\by A. Cayley \paper On certain results relating to
quaternions \jour Phil. Mag. \vol 26 \pages 141--145 \yr 1845
\endref

\ref\key CF\by  P. Cartier and D. Foata \paper Problemes
combinatoires de commutation et rearrangements
\inbook Lecture Notes Math. \vol 85 \yr 1969 \publ Springer-Verlag
\publaddr Berlin--Heidelberg
\endref

\ref\key CS\by A. Connes and A. Schwarz \paper Matrix Vi\`ete theorem
revisited \jour Lett. Math. Phys. \vol 39 \pages 349--353 \yr 1997
\endref

\ref\key Co \by P. M. Cohn \book Skew-field constructions
\bookinfo London Math. Soc. Lecture Notes \vol 27 \yr 1977 \publ Cambridge Univ. Press \publaddr Cambridge
\endref

\ref\key D\by J. Dieudonne \paper Les determinantes sur un corps
non-commutatif
\jour Bull. Soc. Math. France \vol 71 \yr 1943 \pages 27--45
\endref

\ref\key Dr\by P. K. Draxl \book Skew fields
\jour London Math. Soc. Lecture Note Series \vol 81 \yr 1983
\publ Cambridge Univ. Press \publaddr Cambridge
\endref

\ref\key DS \by V. G. Drinfeld and V.V. Sokolov
\paper Lie algebras and equations of Korteweg-de Vries type
\inbook jour Current problems in mathematics. Itogi Nauki i
Tekniki, Akad. Nauk. SSSR \vol 24 \pages 81--180 \yr 1984
\publ VINITI \publaddr Moscow \lang Russian
\endref
 
\ref\key Dy\by F. J. Dyson \paper Quaternion determinants
\jour Helv. Phys. Acta \vol 45 \yr 1972 \pages 289--302
\endref

\ref\key Dy1 \bysame \paper Correlations between eigenvalues of
random matrices
\jour Comm. Math. Phys. \vol 19 \yr 1970 \pages 235--250
\endref

\ref\key EGR\by P. Etingof, I. Gelfand, and V. Retakh\paper
Factorization of differential operators, quasideterminants, and
nonabelian Toda field equations\jour Math. Res. Letters \vol 4
\issue 2-3 \yr 1997\pages 413--425 \endref

\ref\key EGR1 \bysame \paper
Nonabelian integrable systems, quasideterminants, and Marchenko
lemma \jour Math. Res. Letters \vol 5 \issue 1-2 \yr 1998\pages 1--12
\endref

\ref\key ER \by P. Etingof and V. Retakh\paper Quantum determinants
and quasideterminants \jour Asian J. Math. \vol 3 \issue 2
\pages 345--352 \yr 1999
\endref


\ref\key F\by D. Foata \paper A noncommutative version of the
matrix inversion formula
\jour Advances in Math. vol 31 \yr 1979 \pages 330--349
\endref

\ref\key Fr\by G. Frobenius \paper \"Uber lineare Substituenen
und bilinear Formen \jour J. Reine Angew. Math. \vol 84 \pages 1--63
\yr 1877
\endref

\ref\key FIS\by S. H. Friedberg, A. J. Insel, and L. E. Spence
\book Linear algebra\bookinfo Third Edition \publ Prentice Hall,
Upper Saddle River, NJ \yr 1997
\endref

\ref\key G \by F.R. Gantmacher \book The theory of matrices
\publ Chelsea \publaddr New York \yr 1959 \endref

\ref\key GGR \by I. Gelfand, S. Gelfand, and V. Retakh \paper
Noncommutative algebras associated to complexes and graphs \jour
Selecta Math (N.S.) \vol 7\yr 2001 \pages 525--531 \endref
                          
\ref\key GGRSW \by I. Gelfand, S. Gelfand, V. Retakh, S. Serconek,
and R. Wilson\paper Hilbert series of quadratic algebras
associated with decompositions of noncommutative polynomials \jour
J. Algebra \vol 254 \yr 2002\pages 279--299  \endref

\ref\key GKLLRT \by I. Gelfand, D. Krob, A. Lascoux, B. Leclerc, V.
Retakh, and J.-Y. Thibon\paper Noncommutative symmetric functions
\jour Advances in Math.\vol 112\issue 2\yr 1995 \pages 218--348
\endref

\ref\key GR \by I. Gelfand and V. Retakh\paper Determinants of
matrices over moncommutative rings \jour Funct. Anal. Appl.\vol
25 \issue 2 \yr 1991 \pages 91--102 \endref

\ref\key GR1 \bysame\paper A theory of
noncommutative determinants and characteristic functions of graphs
\jour Funct. Anal. Appl.\vol 26 \issue 4 \yr 1992 \pages 1--20
\endref

\ref\key GR2 \bysame \paper A theory of
noncommutative determinants and characteristic functions of
graphs. \rm I \jour Publ. LACIM, UQAM, Montreal \vol 14 \yr 1993
\pages 1--26 \endref

\ref\key GR3 \bysame \paper Noncommutative Vieta
theorem and symmetric functions \inbook Gelfand Mathematical
Seminars 1993--95 \publ Birkh\"auser \publaddr Boston \yr 1996
\pages 93--100
\endref

\ref\key GR4 \bysame \paper Quasideterminants, \rm I
\jour Selecta Math. (N.S.) \vol 3 \issue 4\yr 1997 
\pages 517--546  \endref

\ref\key GRW \by I. Gelfand, V. Retakh, and R. Wilson \paper
Quadratic-linear algebras associated with decompositions of
noncommutative polynomials and Differential polynomials \jour
Selecta Math. (N.S.) \vol 7\yr 2001 \pages 493--523 \endref

\ref\key GRW1 \bysame \paper
Quaternionic quasideterminants and determinants 
\inbook Lie groups and symmetric spaces \bookinfo Amer. Math. Soc.
Transl. Series 2\vol 210 \pages 111--123 \yr 2003 \publ Amer. Math. Soc. 
\publaddr Providence, RI
\endref

\ref\key Gi \by A. Givental \paper Stationary phase integrals,
quantum Toda lattices, flag manifolds and the mirror conjecture
\inbook Topics in Singularity Theory \bookinfo Amer. Math. Soc.
Transl. Series 2\vol 180 \pages 103--115 \yr 1997 \publ Amer. Math. Soc. 
\publaddr Providence, RI
\endref

\ref\key GV \by V. Goncharenko and A. Veselov \paper Darboux
transformations for multidemensional Schr\"o\-dinger operators
\jour J. Phys. A: Math. Gen. \vol 31 \pages 5315--5326 \yr 1998
\endref

\ref\key H\by A. Heyting \paper Die Theorie der linearen Gleichungen
in einer Zahlenspezies mit nicht\-kommutativer Multiplication
\jour Math. Ann. \vol 98 \yr 1927 \pages 465--490
\endref

\ref\key Ho \by R. Howe \paper Perspectives in invariant theory:
Schur duality, multiplicity-free actions and beyond \inbook Israel Math.
Conf. Proc. \vol 8\yr 1995\pages 1--182 \publ Bar-Ilan Univ \publaddr Ramat Gan
\endref

\ref\key Ho1\bysame \paper Remarks on classical invariant theory
\jour Trans. Amer. Math. Soc.
\vol 313 \yr 1989 \pages 539--570
\endref

\ref\key J \by C. Joly \paper Supplement to oevres of Hamilton
\yr 1969 \inbook W. R. Hamilton, Elements of Quaternions \publ Chelsea
\publaddr New York
\endref

\ref\key Ja \by N. Jacobson \book Lectures in abstract algebra. \rm  Vol. II
\yr 1975 \publ Springer-Verlag \publaddr New York--Berlin \endref

\ref\key Ja2 \bysame \book Structure and representations of
Jordan algebras \bookinfo Colloquium Publications 
\yr 1968 \vol 39 \publ Amer. Math. Soc. \publaddr Providence, RI
\endref

\ref \key K \by M. Kapranov
\paper Noncommutative geometry based on commutator expansions
\jour J. Reine Angew. Math. \vol 505 \yr 1998 \pages 73--118
\endref


\ref\key Ki \by A. Kirillov \paper Introduction to family algebras
\jour  Moscow Math. J. \vol 1 \issue 1 \pages 49--63 \yr 2001
\endref

\ref \key KR \by  M. Kontsevich and A. Rosenberg
\paper Noncommutative smooth spaces
\inbook The Gelfand Mathematical Seminars, 1996--1999
\publ Birkh\"auser \publaddr Boston \yr 2000 \pages 85--108
\endref

\ref\key Ko \by B. Kostant \paper The solution to a generalized
Toda lattice and representation theory \jour Advances in Math. \vol 34
\pages 195--338 \yr 1979 \endref

\ref\key KL \by D. Krob and B. Leclerc\paper Minor identities for
quasi-determinants and quantum determinants\jour Comm. Math. Phys.
\vol 169 \issue 1\yr 1995 \pages 1--23 \endref

\ref\key KLT \by D. Krob, B. Leclerc, and J.-Y. Thibon \paper
Noncommutative symmetric functions. {\rm II.} Transformations of
alphabets \jour Internat. J. Algebra Comput. \vol 7 \pages 181--264 \yr 1997
\endref

\ref\key KS \by P. P. Kulish and E. K. Sklyanin \paper Quantum spectral
transform method. Recent developments \inbook Lecture Notes Phys.
\vol 151 \pages 61-119 \yr 1982 \publ Springer-Verlag \publaddr Belin--Heidelberg--New York
\endref

\ref\key L \by T. Y. Lam \paper Representations of finite groups: A
hundred years, \rm Part I \jour Notices Amer. Math. Soc. \vol 45 \pages
361--372 \issue 4 \yr 1998
\endref

\ref\key Le \by D. A. Leites \paper A certain analog of the determinant
\jour Russian Math. Surv. \vol 35 \issue 3 \page
156 \yr 1975 \lang Russian
\endref

\ref\key LST \by B. Leclerc, T. Scharf, and J.-Y. Thibon \paper
Noncommutative cyclic characters of symmetric groups \jour J. 
Combinatorial Theory, A \vol 75 \pages 55--69 \yr 1996
\endref

\ref\key L\"o \by L\"ofwall \paper On the subalgebra generated by the
one-dimensional elements in the Yoneda $Ext$-algebra \jour
Lecture Notes in Math. \vol 1183 \pages 291-338 \yr 1986 \endref

\ref\key Mac \by I. Macdonald
\book Symmetric functions and Hall polynomials
\bookinfo Second edition
\publ The Clarendon Press, Oxford University Press
\publaddr New York \yr  1995
\endref

\ref\key M \by P. A. MacMacon \book Combinatorial analysis
\publ Cambridge Univ. Press \publaddr Cambridge \yr 1915
\endref

\ref\key Ma \by Y. Manin \paper Multiparameter quantum deformations
of the linear supergroup \jour Comm. Math. Phys. \vol 113 \yr 1989
\pages 163--175
\endref

\ref\key Mat \by V. B. Matveev \paper Darboux transformations,
covariance theorems and integrable systems
\inbook L.~D.~Faddeev's Seminar on Mathematical Physics \bookinfo Amer. Math.
Soc. Transl. Series 2\vol 201 \yr 2000 
\pages 179--209 \publ Amer. Math. Soc. \publaddr Providence, RI
\endref

\ref\key Mol \by A. Molev \paper Noncommutative symmetric functions and
Laplace operators for classical Lie algebras\jour Lett. Math.
Phys. \vol 35\issue 2\yr 1995 \pages 135--143 \endref

\ref\key Mol1 \bysame \paper Gelfand-Tsetlin bases for representations of
Yangians \jour Lett. Math. Phys. \vol 30\issue 2\yr 1994 \pages 53--60 \endref

\ref\key MolR \by A. Molev and V. Retakh \paper Quasideterminants and
Casimir elements for the general Lie superalgebras\jour Intern. Math.
Res. Notes \yr 2004 \issue 13 \pages 611-619\endref

\ref\key MB \by E. H. Moore and R. W. Barnard \book General analysis pt. I
\yr 1935 \publ The American phylosophical society \publaddr Philadelphia
\endref

\ref\key Mo \by E. H. Moore \paper On the determinant of an Hermitian
matrix of quaternionic elements \jour Bull. Amer. Math. Soc.
\vol 28 \yr 1922 \pages 161--162
\endref

\ref\key O \by O. Ore \paper Liner equations in noncommutative rings
\jour Ann. Math.
\vol 32 \yr 1931 \pages 463--477
\endref

\ref\key Ok \by K. Okamoto \paper B\"acklund transformations of
classical orthogonal polynomials \inbook Algebraic
Analysis. \rm II \eds M. Kashiwara, T. Kawai
\pages 647--657 \yr 1988\publ Academic Press \publaddr Boston, MA
\endref

\ref\key Os \by B. L. Osofsky \paper Quasideterminants and roots of
polynomials over division rings \paperinfo
Proc. Intern. Lisboa Conf. on Algebras, Modules and Rings, 2003
\toappear \endref


\ref\key P \by A. Polischuk \paper Triple Massey products on curves,
Fay's trisecant identity and tangents to the canonical embedding
\jour math.AG/0109007 \yr 2001
\endref

\ref\key PPR \by I. Pak, A. Postnikov and V. Retakh \paper Noncommutative
Lagrange inversion  \paperinfo preprint \yr 1995
\endref

\ref\key RS \by  A. Razumov and M. Saveliev \paper Maximally nonabelian
Toda systems \jour Nuclear Phys. B\vol 494 \issue 3 \yr 1997 \pages 657--686
\endref

\ref\key RRV \by V. Retakh, C. Reutenauer, and A. Vaintrob \paper
Noncommutative rational functions and Farber's invariants of
boundary links \inbook Differential Topology,
Infinite-dimensional Lie Algebras and Applications \bookinfo Amer. Math. Soc.
Transl. Series 2\vol 194 \pages 237--246 \yr 1999 \publ Amer. Math. Soc. 
\publaddr Providence, RI
\endref


\ref\key RSh \by V. Retakh and V. Shander \paper Schwarz derivative for
noncommutative differential algebras \inbook Unconvential Lie Algebras 
\bookinfo Adv. Soviet Math. \pages 130--154 \vol 17 \yr 1993\publ Amer. Math.
Soc.  \publaddr Providence, RI \endref

\ref\key Re \by C. Reutenauer \paper Inversion height in free fields
\jour Selecta Math. (N.S.) \vol 2 \issue 1 \yr 1996 \pages 93--109
\endref

\ref\key Ri \by A. R. Richardson \paper Hypercomplex determinants
\jour Messenger of Math.
\vol 55 \yr 1926 \pages 145--152
\endref

\ref\key Ri1 \bysame \paper Simultaneous linear equations over
a division algebra
\jour Proc. London Math. Soc.
\vol 28 \yr 1928 \pages 395--420
\endref

\ref \key SvB \by J. T. Stafford and M. van den Bergh
\paper Noncommutative curves and noncommutative surfaces
\jour  Bull. Amer. Math. Soc. \vol 38 \yr 2001 \issue  2
\pages  171--216
\endref


\ref\key SW \by S. Serconek and R. L. Wilson \paper
Quadratic algebras
associated with decompositions of noncommutative polynomials
are Koszul algebras
\paperinfo in preparation
\endref

\ref\key Sch \by M. Schork \paper The bc-system of higher rank revisited
\jour J. Phys. A: Math. Gen. \vol 35 \yr 2001 \pages 2627--2637
\endref

\ref\key Schur \by I. Schur \paper \"Uber die charakteristischen Wurzeln einer
linearen Substitution mit einer Anwendung auf die Theorie der
Integralgleichungen
\jour Math. Ann.
\vol 66 \yr 1909 \pages 488--510
\endref


\ref\key St \by E. Study \paper Zur Theorie der lineare Gleichungen
\jour Acta Math.\vol 42
\pages 1--61 \yr 1920
\endref

\ref\key Syl\by  J. J. Sylvester \paper Sur une classe nouvelle d'equations
differ\'entielles et d'equations aux differences finies d'une forme
int\'egrable 
\jour Compt. Rend. Acad. Sci.
\vol 54 \yr 1862 \pages 129-132, 170 
\endref

\ref\key Thi \by J.-Y. Thibon \paper
Lectures on moncommutative symmetric functions
\jour MSJ Mem. \vol 11 \yr 2001 \pages 39--94
\endref

\ref\key VP \by P. Van Praag \paper Sur le norme reduite du determinant
de Dieudonne des matrices quaterniennes \jour J. Algebra \vol 136
\yr 1991 \pages 265--274
\endref


\ref\key W \by J. H. M. Wedderburn \paper On continued fractions in
noncommutative quantities
\jour Ann. Math.
\vol 15 \yr 1913 \pages 101--105
\endref

\ref\key We \by H. Weyl \book Classical groups, their invariants and representations
\publ Princeton Univ. Press \publaddr Princeton, NJ \yr 1946
\endref

\ref\key Wi \by R. L. Wilson \paper Invariant polynomials in the free skew
field
\jour Selecta Math.
\vol 7 \yr 2001 \pages 565--586
\endref

\endRefs
\enddocument